\documentclass[11pt,a4paper,oneside]{article} 
\usepackage[utf8]{inputenc}
\usepackage[T1]{fontenc}
\usepackage[affil-it]{authblk}

\usepackage{graphicx}
\usepackage{grffile}
\usepackage{caption}
\usepackage{latexsym}
\usepackage{graphpap}
\usepackage[english]{babel}
\usepackage[nodisplayskipstretch]{setspace}
\usepackage{amsmath,amsfonts,amssymb}
\usepackage{amsthm}
\usepackage{stmaryrd} 
\usepackage{caption}
\usepackage{appendix}
\usepackage{subcaption} 

\newtheorem{proposition}{Proposition}

\usepackage{amsmath}
\usepackage{tikz}
\usepackage{float}
\usepackage{lineno}
\usepackage{comment}
\usepackage{graphicx}
\usepackage[nomentbl]{nomencl}

\usepackage[colorlinks=true, citecolor=blue, linkcolor=blue,plainpages=false,pdfpagelabels]{hyperref}

    \newcommand{\incr}{\text{\raisebox{.7mm}{\scalebox{.6}{$\nearrow$\ }}}}
    \newcommand{\decr}{\text{\raisebox{.6mm}{\scalebox{.6}{$\searrow$\ }}}}
    
\usepackage{chngcntr}
\usepackage{lineno}
\usepackage{fancyhdr}

\usepackage{url}
\usepackage{lastpage}

\usepackage{anysize}
\marginsize{1.5cm}{1.5cm}{0.5cm}{1cm}

\usepackage[normalem]{ulem} 

\numberwithin{equation}{section}
\numberwithin{figure}{section}
\numberwithin{table}{section}

{\theoremstyle{plain}
\newtheorem{theorem}{Theorem}[section]

}

{\theoremstyle{definition}
\newtheorem{definition}[theorem]{Definition}

\newtheorem{restrict}[theorem]{Restriction}
}

{\theoremstyle{remark}
\newtheorem{remark}[theorem]{Remark}
}

\numberwithin{equation}{section}
\numberwithin{figure}{section}
\numberwithin{proposition}{section}

\begin{document}

\def\Box{ \framebox[5.5pt]}

\title{\centering \large \bf Solving Riemann Problems with a Topological Tool (Extended version)}

\author{Cesar S. Eschenazi
  \thanks{email: \texttt{cseschenazi\emph{{@}}gmail.com}}}
\affil{Universidade Federal de Minas Gerais -- UFMG\\ Belo Horizonte, MG, 31270-901}
\author{Wanderson J. Lambert
  \thanks{email: \texttt{wanderson.lambert\emph{{@}}unifal-mg.edu.br}  - Corresponding author.}}
\affil{Universidade Federal de Alfenas -- UNIFAL\\ Alfenas, MG, 37130-001}
\author{Marlon M. L\'{o}pez-Flores%
  \thanks{email: \texttt{mmlf\emph{{@}}impa.br}}}
\affil{Instituto Nacional de Matem\'{a}tica Pura e Aplicada -- IMPA\\ Rio de Janeiro, RJ, Brazil, 22460-320}
\author{Dan Marchesin%
  \thanks{email: \texttt{marchesi\emph{{@}}impa.br}}}
\affil{Instituto Nacional de Matem\'{a}tica Pura e Aplicada -- IMPA\\ Rio de Janeiro, RJ, Brazil, 22460-320}
\author{Carlos F.B. Palmeira
  \thanks{email: \texttt{fredpalm\emph{{@}}gmail.com}}}
\affil{Pontifícia Universidade Católica do Rio de Janeiro -- PUC-Rio\\
Rio de Janeiro, RJ - Brasil, 22451-900}
\author{Bradley J. Plohr
  \thanks{email: \texttt{bradley.j.plohr\emph{{@}}gmail.com}}}
\affil{Los Alamos, New Mexico, USA}                                                       \date{}             
\maketitle

\begin{abstract}
In previous work, we developed a topological framework
for solving Riemann initial-value problems
for a system of conservation laws.
Its core is a differentiable manifold,
called the wave manifold,
with points representing shock and rarefaction waves.
In the present paper,
we construct, in detail, the three-dimensional wave manifold
for a system of two conservation laws with quadratic flux functions.
Using adapted coordinates,
we derive explicit formulae for important surfaces and curves
within the wave manifold and display them graphically.
The surfaces subdivide the manifold into regions according to shock type,
such as ones corresponding to the Lax admissibility criterion.
The curves parametrize rarefaction, shock, and composite waves
appearing in contiguous wave patterns.
Whereas wave curves overlap in state space,
they are disentangled within the wave manifold.
We solve a Riemann problem by constructing a wave curve
associated with the slow characteristic speed family,
generating a surface from it using shock curves,
and intersecting this surface with a fast family wave curve.
This construction is applied to solve Riemann problems
for several illustrative cases.
\end{abstract}

\newpage
\tableofcontents
\section{Introduction}
\label{sec:Introduction}

Solutions of Riemann problems
for a hyperbolic system of conservation laws
are often complicated,
curtailing their usage as building blocks for general solutions.
This complication has motivated us to develop
a global perspective on constructing solutions.
Previous work proposed a topological framework
for representing the fundamental waves
appearing in Riemann solutions~\cite{isaacson92}
and for building solutions of Riemann problems~\cite{AEMP10}.
In the present work,
we apply this framework to solve Riemann problems
for a particular system of two conservation laws
with quadratic polynomial flux functions.
Focusing on this model allows us to derive explicit formulae for,
and to present graphical visualizations of,
the constructs in the topological framework.

\subsection{Conservation laws and Riemann problems}
Hyperbolic systems of conservation laws
are ubiquitous in applications, \emph{e.g.}, Fluid Dynamics.
Such a system of partial differential equations (PDE)
governs a \emph{state} $W(x, t) \in \mathbb{R}^n$ \nomenclature[W]{$W(x, t)$}{state}{}{}
as a function of \emph{time} $t \ge 0$ \nomenclature[t]{$t$}{time}{}{} and $\emph{position}$ $x$;\nomenclature[x]{$x$}{position}{}{}%
it takes the form
\begin{equation}
W_t + F(W)_x = 0,
\label{eq:cons-law}
\end{equation}
where $F(W) \in \mathbb{R}^n$\nomenclature[F]{$F(W)$}{flux}{}{} is the \emph{flux}.
Solutions comprise interacting nonlinear waves,
such as smooth \emph{rarefaction waves} that fan out
and discontinuous \emph{shock waves} that develop from compression.
The possible waves and their interactions are elucidated by
restricting to one spatial dimension, $x \in \mathbb{R}$,
and solving initial-value problems having an extremely simple form,
first studied by Riemann~\cite{rie}:
for a specified \emph{left state} $W_L$\nomenclature[W]{$W_L$}{left state}{}{} and \emph{right state} $W_R$\nomenclature[W]{$W_R$}{right state}{}{},
\begin{equation}
W(x, t=0) = \begin{cases}
W_L & \text{if $x < 0$}, \\
W_R & \text{if $x > 0$}.
\end{cases}
\label{eq:initial_condition}
\end{equation}
As we will see, a \emph{Riemann problem solution} (RPS)
is a sequence of rarefaction and shock waves
separated by regions in which the state is uniform.

\subsection{Importance of Riemann problems}
Riemann problems are crucial for understanding
hyperbolic systems of conservation laws,
as evidenced by Refs.~\cite{smoller94,serre,bressan1,lefloch02,dafermos}.
The constituent waves of a general solution
take an analytically tractable form within an RPS,
so that Riemann problems facilitate analysis of
nonlinear wave formation, propagation, and interaction.
Additionally,
solving Riemann problems is key to constructing
weak solutions for the Cauchy problem through
Glimm's random choice method~\cite{glimm,liu4,glimm_marshall_plohr84,hong}
and front tracking methods~\cite{bressan,holden}.
Furthermore,
advanced numerical methods for simulating hyperbolic systems
utilize Riemann solutions:
the primary step in Godunov's method
and its refinements~\cite{godunov,engquist1,harten,osher2,roe}
is finding an approximate RPS,
and the interface tracking method~\cite{glim2,glimFT,gliKliMcB85}
solves Riemann problems dominated by the tracked wave.
Riemann solutions are also invaluable as test cases
for verifying the accuracy of numerical simulations
of hyperbolic systems~\cite{toro}.

\subsection{Overview of the classical approach}
A distinctive feature in an RPS is
a jump discontinuity that propagates at constant speed:
\begin{equation}
\label{eq:shock}
W(x, t) = \begin{cases}
W & \text{for $x < \sigma\,t$}, \\
W' & \text{for $x > \sigma\,t$}.
\end{cases}
\end{equation}
The constant $\sigma$\nomenclature[s]{$\sigma$}{shock propagation speed}{}{} is the \emph{shock propagation speed},
whereas the constant states $W$ and $W'$,
with $W \ne W'$,
are the \emph{left} and \emph{right states} of the jump.
Because such a function is not differentiable,
Eq.~\eqref{eq:cons-law} does not apply.
Instead,
this PDE is interpreted as requiring conservation of $W$ in a weak sense,
which for a jump discontinuity reduces to
the Rankine-Hugoniot condition \cite{smoller94}:
\begin{equation}
\label{RH-condition}
F(W) - F(W') = \sigma\,(W - W').
\end{equation}
In the context of solving Riemann problems,
a \emph{shock wave} is a jump discontinuity~\eqref{eq:shock}
that satisfies Eq.~\eqref{RH-condition}.
A \emph{rarefaction wave}
is a one-parameter family of infinitesimal shock waves,
for which $\sigma$ is an eigenvalue of $D F(W)$.

The model examined in the present paper has $n = 2$:
a state has components $W = (u, v)$
and \emph{state space} is the $(u, v)$-plane.
The \emph{strictly hyperbolic region} in state space
is the subset where the derivative $D F(W)$ has real distinct eigenvalues,
$\lambda_s(W) < \lambda_f(W)$,
called the \emph{slow} and \emph{fast characteristic speeds}.

If we set $W = W_L$,
the states $W'$ that satisfy Eq.~\eqref{RH-condition}
for some $\sigma$ form the \emph{Hugoniot locus} of $W_L$ in state space,
which is a one-parameter family of shock waves.
Suppose that $W_L$ lies in the strictly hyperbolic region.
For $W'$ near $W_L$,
this locus consists of two curves that cross at $W_L$,
the \emph{slow} and \emph{fast branches}.
Rarefaction waves having $W_L$ at its left edge
also form \emph{slow} and \emph{fast rarefaction curves} in state space,
which cross at $W_L$.
These curves can be extended to form \emph{wave curves}
that parametrize the possible wave patterns with left state $W_L$,
including \emph{composite waves} comprising
adjacent shock and rarefaction waves.
Similarly, the solution set of Eq.~\eqref{RH-condition} with $W' = W_R$
forms the \emph{Hugoniot\,$'$ locus},
rarefaction waves with right edge state $W_R$ form curves,
and wave curves with right state $W_R$ parametrize wave patterns.
Such wave curves are essential
constructs~\cite{bet42,wen72,wen72a,liu1,liu}.
The classical approach to solving a Riemann problem
is to intersect the slow wave curve for $W_L$
with the fast wave curve for $W_R$.

\subsection{Overview of the topological approach}
A topological perspective on waves appearing in Riemann solutions
was proposed for quadratic models in \cite{marpal94b}
and generalized in \cite{isaacson92}.
Its framework is a space $\mathcal{W}$\nomenclature[W]{$\mathcal{W}$}{wave manifold}{}{},
termed the \emph{wave manifold},
in which fundamental constructs
(\emph{e.g.}, waves, wave curves, and bifurcation loci)
are represented as points, curves, and surfaces, respectively.
This approach was extended in \cite{AEMP10}
to construct solutions of Riemann problems.

A point in the wave manifold represents a shock wave,
\emph{i.e.}, a solution of the Rankine-Hugoniot condition~\eqref{RH-condition}.
The speed $\sigma$ can be eliminated
between the two components of this equation,
so that a shock wave has coordinates
$(W, W') = (u, v, u', v') \in \mathbb{R}^4$
constrained by a single equation.
Although this equation has a singularity when $W = W'$,
a blowing-up construction yields
the smooth three-dimensional manifold $\mathcal{W}$.

Many concepts from the theory of conservation laws
are represented by features within the wave manifold.
A Hugoniot locus in state space corresponds
to the curve within $\mathcal{W}$ along which $(u, v)$ is fixed,
called a \emph{Hugoniot curve},
and holding $(u', v')$ fixed gives a Hugoniot$'$ curve.
Viewing a Hugoniot locus as a curve in three-dimensions has an advantage:
whereas each Hugoniot locus in state space has a self-intersection
and the loci for distinct left states can overlap,
Hugoniot curves within $\mathcal{W}$ do not intersect.
In other words, Hugoniot curves form
a one-dimensional \emph{foliation}~\cite{camlin85} in $\mathcal{W}$.
Similarly, the two families of rarefaction curves that overlap
in state space correspond to a foliation
on a two-dimensional submanifold $\mathcal{C} \subseteq \mathcal{W}$,\nomenclature[C]{$\mathcal{C}$}{characteristic manifold}{}{}
called the \emph{characteristic manifold},
which is a two-fold cover of the strictly hyperbolic region in state space.
In this sense, the wave manifold disentangles shock and rarefaction curves.

Some jump discontinuities that satisfy Eq.~\eqref{RH-condition}
are unsuitable in a Riemann solution because they violate physical principles.
Accordingly, shock waves are required to satisfy
an \emph{admissibility criterion}.
For the model studied in this work,
we adopt the \emph{Lax admissibility criterion} \cite{lax}
that the shock speed and the characteristic speeds obey
certain strict inequalities.
These inequalities define open regions within $\mathcal{W}$,
which are bounded by $\mathcal{C}$ together with
the \emph{sonic manifolds} $Son$ and $Son'$
where the shock speed equals a characteristic speed.
Like rarefaction waves within $\mathcal{C}$,
composite waves form a foliation within $Son'$.

By combining portions of Hugoniot, rarefaction, and composite curves
we obtain the counterparts within $\mathcal{W}$ of wave curves in state space.
Solving a Riemann problem
by intersecting slow and fast wave curves in state space
corresponds to a construction in the wave manifold
identified in \cite{AEMP10}:
define the \emph{intermediate surface} to be
the union of Hugoniot$'$ loci starting at points on the slow wave curve;
a Riemann solution corresponds to an intersection
of the reflected fast wave curve with the intermediate surface.
(The \emph{reflection} operation interchanges left and right states states.
Its role was overlooked in \cite{AEMP10};
we correct this omission here.)

Several authors \cite{palmeira88,marpal94b,isaacson92,canplo95,eschenazi02,%
bastos05,wenplo07,AEMP10,eschenazi13,lopez20}
have investigated the wave manifold
when the two components of $F$ are quadratic polynomials.
For such models,
$\mathcal{W}$ is $M_2 \times \mathbb{R}$,
where $M_2$\nomenclature[M]{$M_2$}{M\"obius strip}{}{} is the M\"obius strip,
while $\mathcal{C}$, $Son$, and $Son'$ are ruled surfaces
embedded as cylinders into $\mathcal{W}$.
The shock, rarefaction, and composite foliations have been characterized,
and the question of shock admissibility has been explored.
We intend to add to this wealth of abstract information about quadratic models
by providing details for a specific example
and presenting graphical visualizations of the results.

\subsection{Present work}
In this paper, we apply the topological framework
to a specific quadratic model
belonging to symmetric Case~IV of the Schaeffer-Shearer
classification \cite{schaeffer87}, as modified in \cite{palmeira88}.
This model was solved in \cite{isaacson88}
when it is strictly hyperbolic everywhere except at one point;
here the strictly hyperbolic region is the exterior of an ellipse.
We introduce a coordinate system on $\mathcal{W}$
that is convenient for this model
and develop explicit formulae,
involving polynomial and rational functions,
for the primary structures in $\mathcal{W}$,
including its subdivision according to shock type
and the separation of sonic surfaces into slow and fast portions.
In particular,
we construct examples of wave curves
and use them to solve representative Riemann problems.
Central to this work is visualizing the topological components
of the three-dimensional wave manifold.

\subsection{Organization of this work}
This paper is organized as follows.
In Sec.~\ref{sec:Wave_Manifold_and_Hugoniot_Curves},
we review basic facts and definitions from \cite{isaacson92},
introduce adapted coordinates,
and describe the characteristic, sonic, and sonic\,$'$ surfaces
for the model being studied.
In Sec.~\ref{sec:Decomposition_of_the_Wave_Manifold},
we describe how these surfaces subdivide the wave manifold
into twelve regions. We locate the regions where Lax inequalities are satisfied.)
We also characterize a new surface,
formed by Hugoniot curves through points
for which the characteristic speeds coincide,
and demonstrate its tangency to the sonic surfaces.
This tangency curve separates slow and fast portions,
as demonstrated in Sec.~\ref{sec:sonlifs}.

Section~\ref{subsec:raref} addresses rarefaction solutions,
while Lax shock waves are the focus of Sec.~\ref{sec:lax}. In Sec.~\ref{subsec:NEW4.3}, we incorporate Lax conditions into the decomposition of $\mathcal{W}$. Analysis of composite curves is presented in Sec.~\ref{subsec:Compos}.
Then, in Sec.~\ref{sec:Wave_Curves},
we define wave curves in $\mathcal{W}$,
which parametrize contiguous groups
of rarefaction and shock waves having the same family, slow or fast.
In Sec.~\ref{sec:Wave_Groups_and_the_Intermediate_Surface},
we introduce and give an example of the intermediate surface
for a slow wave curve.
This surface allows us to obtain the solution
for a Riemann problem in the wave manifold.
In Sec.~\ref{sec:Riemann_Solutions},
the elements we have defined and constructed
are used to solve Riemann problems.
Representative examples are presented in Sec.~\ref{sec:examples}.
Finally, some technical results have been relegated to Appendices.

\section{Wave Manifold and Hugoniot Curves}
\label{sec:Wave_Manifold_and_Hugoniot_Curves}

This section reviews the construction of the wave manifold
and of Hugoniot curves within it.
We refer the reader to Sec.~2 of \cite{eschenazi02} for further details.

\subsection{Quadratic model}
The model treated in this paper
is a two-component systems of conservation laws
with quadratic polynomial flux functions.
It falls in the symmetric Case~IV \cite{schaeffer87,palmeira88}
with the strictly hyperbolic being the exterior of an ellipse.
As shown in Appendix~\ref{subsec:Normal_form_for_the_flux} and \cite{palmeira88},
it suffices to examine Eq.~\eqref{eq:cons-law} with $F=(f,g)$
having the normal form
\begin{align}
f(u,v) &= (b_1 + 1)\,u^2/2 + v^2/2 + a_1\,u + a_2\,v,
\label{eq:f} \\
g(u,v) &= u\,v - b_2\,v^2/2 + a_3\,u + a_4\,v.
\label{eq:g}
\end{align}
with $b_1 > 1$ and $b_2 = 0$.
We choose this case for simplicity and specificity,
but we expect that the treatment presented here can be extended to other cases,
such as those described in the Appendix of \cite{schaeffer87}.
The general solution of Riemann problems for a symmetric Case~IV model
with homogeneous quadratic flux functions
(so that strict hyperbolicity fails only at the origin)
is found in \cite{isaacson88};
see also \cite{matos2015}.

\subsection{Wave manifold}\label{sec:2.2}
To study Hugoniot curves,
we consider $\mathbb{R}^4$ with coordinates $(u,v,u',v')$
and solve the equation obtained from the pair~\eqref{RH-condition}
by eliminating $\sigma$, i.e.,
\begin{equation}
[f(u,v)-f(u',v')](v-v')=[g(u,v)-g(u',v')](u-u').
\end{equation}
Since this  equation
is singular along the diagonal $W=W'$,
we perform a blowing-up,
which in this case consists of
introducing the following coordinate transformation:
\begin{equation}\label{transf1}
\text{$U=(u+u')/2$, $V=(v+v')/2$, $X=u-u'$, $Y=v-v'$, and $Z=Y/X$}.
\end{equation}
Using \eqref{eq:f} and \eqref {eq:g}, we get
\begin{equation}
X^{2}[(1-Z^{2})\left(V+a_{2}\right)-Z\left(b_{1}\,U+a_{1}-a_{4}\right)+c]=0
\quad\text{and}\quad Y=ZX,
\label{eqrh2}
\end{equation}
and we complete the blowing-up by discarding the factor $X^{2}$
from the first equation.
Here $c=a_3-a_2$,
which we take to be positive without loss of generality.
These two equations in $(U, V, X, Y, Z)$-space define a smooth,
three-dimensional manifold,
referred to as the \emph{wave manifold} and denoted by $\mathcal{W}$.
We often use the notation $\mathcal{U}$ for a point in $\mathcal{W}$.

\begin{remark}
One can also regard the wave manifold as
a three-dimensional submanifold of $\mathbb{R}^5$
with coordinates $(W,\sigma,W')$.
The solution set of Eq.~\eqref{RH-condition}
is the union of this submanifold and the diagonal, $W' = W$,
which is eliminated by the blowing up.
This approach was used in~\cite{isaacson92}
for systems of conservation laws with more than two equations.
\end{remark}

A fundamental operation in $\mathcal{W}$ is defined as follows.

\begin{definition}
The \emph{reflection} of a point $\mathcal{U} \in\mathcal{W}$
replaces $Y$ by $-Y$ and $X$ by $-X$.
The point in $\mathcal{W}$ corresponding to $(u,v,u',v')$
is mapped to the point corresponding to $(u',v',u,v)$.
\label{def:reflection}
\end{definition}

Defining
\begin{equation}\label{transf2}
\widetilde{U}=b_1\,U+a_1-a_4\qquad \text{and} \qquad \widetilde{V}=V+a_2,
\end{equation}
and substituting in Eq.~\eqref{eqrh2}, we obtain the equation $G_{1}=0$, where
\begin{equation}\label{G1}
G_1=(1-Z^2)\,\widetilde{V}-Z\,\widetilde{U}+c,
\end{equation}
which was formulated in \cite{marpal94b,eschenazi02,eschenazi13}.

\subsection{Hugoniot curves}
Given a state $W_0 = (u_0, v_0)$,
the \emph{Hugoniot curve} with left state $W_0$
comprises points of $\mathcal{W}$ for which $W = W_0$,
\emph{i.e.}, $U + X/2 = u_0$ and $V + Z\,X/2 = v_0$.
In terms of $\widetilde{U}$ and $\widetilde{V}$ from~\eqref{transf2},
these equations become
\begin{equation}\label{eq:KL}
\text{$\widetilde{U}+b_{1}\,X/2=\widetilde{U}_0$
and $\widetilde{V}+Z\,X/2=\widetilde{V}_0$},
\end{equation}
where $\widetilde{U}_0=b_{1}\,u_{0}+a_{1}-a_{4}$
and $\widetilde{V}_0=v_{0}+a_{2}$.

Equivalently, a Hugoniot curve is a solution
of $dG_1 = 0$, $du = 0$, and $dv = 0$.
As shown in \cite{eschenazi02},
the differential forms $dG_1$, $du$, $dv$ are linearly dependent,
\emph{i.e.}, have a singularity,
when $Z = 0$, $\widetilde{U} + X/2 = 0$, and $\widetilde{V} = - c$.
These equations define the \emph{secondary bifurcation line} $\mathcal{B}$.
Therefore, Hugoniot curves foliate $\mathcal{W} \setminus \mathcal{B}$.

The \emph{Hugoniot\,$'$ curve} with right state $W_0$
is the reflection of the Hugoniot curve with left state $W_0$;
its defining equations, obtained simply by replacing $X$ with $-X$, are
\begin{equation}\label{eq:KL'}
\text{$\widetilde{U}-b_{1}\,X/2=\widetilde{U}_0$
and $\widetilde{V}-Z\,X/2=\widetilde{V}_0$}.
\end{equation}
With $\mathcal{B}'$ denoting the reflection of $\mathcal{B}$,
along which $Z = 0$, $\widetilde{U} - X/2 = 0$, and $\widetilde{V} = - c$,
Hugoniot$'$ curves foliate $\mathcal{W} \setminus \mathcal{B}'$.

If $\mathcal{U}_0 \in \mathcal{W}$,
then the Hugoniot curve that contains $\mathcal{U}_0$
is denoted $\mathcal{H}(\mathcal{U}_0)$.
It is the Hugoniot curve with the same left state as $\mathcal{U}_0$.
Similarly, the Hugoniot$'$ curve containing $\mathcal{U}_0$,
which is the Hugoniot$'$ curve with same right state as $\mathcal{U}_0$,
is denoted $\mathcal{H}'(\mathcal{U}_0)$.

\subsection{Boundaries of shock admissibility regions}
\label{subsec:Boundaries-admiss-regions}
We describe three important surfaces in the study of the wave manifold:
the characteristic surface $\mathcal{C}$,
the sonic surface $Son$,
and the sonic\,$'$ surface $Son'$.
These surfaces are \emph{boundaries of shock admissibility regions};
see Sec.~\ref{sec:lax}.

\begin{definition}
\label{carateristica}
The \emph{characteristic surface} $\mathcal{C}$
is the set of points in $\mathcal{W}$ with $X = 0$.
For such points, $(u, v) = (u', v')$.
\end{definition}

This surface is a natural place for defining rarefaction curves,
as discussed in Sec.~\ref{subsec:raref}.

To define $Son$ and $Son'$,
we express the \emph{shock speed} $\sigma$ as the real-valued function
defined on $\mathcal{W}$ by $[f(u,v)-f(u',v')]/(u-u')$.
In $(\widetilde{U}, \widetilde{V}, X, Z)$-coordinates,
we find that
\begin{equation}
\sigma=(b_1 + 1)\,\widetilde{U}/b_1 + Z\,\widetilde{V} + \sigma_0,
\label{eq:sp1}
\end{equation}
where $\sigma_0 = [(b_1 + 1)\,a_4 - a_1]/b_1$.
For $\mathcal{U} \in \mathcal{W}$, this speed is denoted $\sigma(\mathcal{U})$.

\begin{definition}
The \emph{sonic surface} $Son$ is the set of critical points
for $\sigma$ as restricted to Hugoniot curves.
Similarly, the \emph{sonic\,$'$ surface} $Son'$
is the set of the critical points of $\sigma$
restricted to Hugoniot$'$ curves.
\end{definition}
The equations for $Son$ and for $Son'$ are derived in \cite{marpal94b}:
\begin{align}
\left[Z^2 + b_1 + 1\right] \widetilde{U}
+ Z \left[Z^2 + b_1 + 3\right] \widetilde{V}
- \left[Z^2 - (b_1 + 1)\right] \frac{X}{2} &= 0,
\label{eq:eqson} \\
\left[Z^2 + b_1 + 1\right] \widetilde{U}
+ Z \left[Z^2 + b_1 + 3\right] \widetilde{V}
+ \left[Z^2 - (b_1 + 1)\right] \frac{X}{2} &= 0,
\label{eq:eqson'}
\end{align}
Notice that Eq.~\eqref{eq:eqson'} is the reflection of Eq.~\eqref{eq:eqson}.

\subsection{Previous results}
We list previous results about $\mathcal{C}$, $Son$, and $Son'$
for the quadratic model considered here.

\begin{itemize}
\item Topologically, $\mathcal{C}$ is a cylinder;
see \cite{eschenazi02}.

\item A Hugoniot curve either
intersects $\mathcal{C}$ transversely at two points,
intersects once where it is tangent to $\mathcal{C}$,
or does not intersect $\mathcal{C}$.
See \cite{eschenazi02}
and Sec.~\ref{subsec:coincidence_curve}.

\item The set of Hugoniot tangency points
forms a simple closed curve within $\mathcal{C}$,
denoted by $\mathcal{E}$
and called the \emph{coincidence curve}.
It divides $\mathcal{C}$ into two components,
denoted by $\mathcal{C}_{s}$ and $\mathcal{C}_{f}$.
See~\cite{AEMP10} and Sec.~\ref{subsec:coincidence_curve}.

\item If the Hugoniot curve $\mathcal{H}(\mathcal{U})$
of a point $\mathcal{U} \in \mathcal{W}$ intersects $\mathcal{C}$
at two distinct points,
one belongs to $\mathcal{C}_s$ and the other belongs to $\mathcal{C}_f$.
Let $\mathcal{U}_{s} = \mathcal{H}(\mathcal{U}) \cap \mathcal{C}_{s}$
and $\mathcal{U}_{f} = \mathcal{H}(\mathcal{U}) \cap \mathcal{C}_{f}$.
Then $\sigma(\mathcal{U}_{s}) < \sigma(\mathcal{U}_{f})$.
Moreover, $\sigma(\mathcal{U})$ is distinct from
$\sigma(\mathcal{U}_{s})$ and $\sigma(\mathcal{U}_{f})$
unless $\mathcal{U}$ belongs to $\mathcal{C}$ or $Son'$.
Analogous statements hold for a Hugoniot$'$ curve.
See~\cite{AEMP10}, Prop.~\ref{th:slow_fast}, and Sec.~\ref{sec:lax}.

\item Topologically,
$Son$ and $Son'$ are cylinders that intersect
$\mathcal{C}$ along a curve called the \emph{inflection locus},
denoted by $\mathcal{I}$,
which contains the point $\mathcal{B} \cap \mathcal{C}$.
Additionally, $Son$ and $Son'$ intersect along two straight lines,
called the \emph{double sonic locus} and denoted $\mathcal{D}$,
which intersect $\mathcal{C}$ transversely at two points
in the inflection locus.
See \cite{isaacson92}
and Sec.~\ref{sec:Subdivision_of_the_wave_manifold}.

\item Generically, $\mathcal{H}(\mathcal{U})$ intersects $Son'$ at
0, 2, or 4 points. If there are two intersection points,
$\mathcal{U}_{1}$ and $\mathcal{U}_{2}$,
then either
$\sigma(\mathcal{U}_{1})=\sigma(\mathcal{U}_{2})=\sigma(\mathcal{U}_{s})$ or
$\sigma(\mathcal{U}_{1})=\sigma(\mathcal{U}_{2})=\sigma(\mathcal{U}_{f})$,
whereas if there are four intersection points,
$\mathcal{U}_{1}$, $\mathcal{U}_{2}$, $\mathcal{U}_{3}$, and $\mathcal{U}_{4}$,
then $\sigma(\mathcal{U}_{1})=\sigma(\mathcal{U}_{2})
=\sigma(\mathcal{U}_{s})<\sigma(\mathcal{U}_{3})
=\sigma(\mathcal{U}_{4})=\sigma(\mathcal{U}_{f})$;
see \cite{eschenazi02}.
\end{itemize}

\subsection{Adapted coordinates}
\label{subsec:adapted_coordinates}
Let $\mathbb{R}P^1$ be the space of lines
through the origin in the $(X, Y)$-plane.
The coordinate $Z = Y/X$
is the slope of the line through $(X, Y) \ne (0, 0)$,
provided that $X \ne 0$, \emph{i.e.}, the line is not vertical.
To include the vertical line,
we use the coordinate $z = X/Y$,
the inverse slope of a line that is not horizontal.
Away from the vertical and horizontal lines,
the coordinate change $z = 1/Z$ is a diffeomorphism,
so that $\mathbb{R}P^1$ is a manifold.

Substituting $Z = 1/z$ into the defining equation $G_1 = 0$ for $\mathcal{W}$,
introducing $V_{1}=V+a_{3}=\widetilde{V}+c$,
and omitting an overall factor of $1/z^2$,
we arrive at the equation $G=0$, where
\begin{equation}\label{G}
G=(z^{2}-1)\,V_{1}-z\,\widetilde{U}+c.
\end{equation}
In addition,
the formula $\sigma=[g(u,v)-g(u',v')]/(v-v')$ reads
\begin{equation}
\sigma= \widetilde{U}/b_1 + z\,V_1 + \sigma_0.
\label{eq:sp2}
\end{equation}
where, again, $\sigma_0 = [(b_1 + 1)\,a_4 - a_1]/b_1$.

As a result, $\mathcal{W}$ is a submanifold of
$\mathbb{R}^4 \times \mathbb{R}P^1$,
being the solution set of $G_1 = 0$ and $Y = ZX$ away from the plane $z = 0$
and of $G = 0$ and $X = zY$ away from the plane $Z = 0$.
Recall, however, that the secondary bifurcation lines
$\mathcal{B}$ and $\mathcal{B}'$ belong to the plane $Z = 0$,
which we denote $\Pi$.
Therefore,
the Hugoniot and Hugoniot$'$ foliations are non-singular
on $\mathcal{W} \backslash \Pi$.
For this reason, $z$ and $Y$ are advantageous as coordinates.

Furthermore, we can solve $G = 0$ by introducing another coordinate.
For each $z$, the solution set of $G=0$ is a line,
which we parametrize by $\tau \in \mathbb{R}$:
\begin{equation}\label{eq:tdef}
\widetilde{U}=2\,c\,z/(z^{2}+1)+c\,\tau\,(z^{2}-1)
\qquad\text{and}\qquad V_{1}=c/(z^{2}+1)+c\,\tau\,z.
\end{equation}
For details, see \ref{subsec:Coordinates}.
Therefore, $\mathcal{W}\setminus\Pi$ has global coordinates
$(z, \tau, Y) \in \mathbb{R}^3$.

Formulae involving $(\widetilde{U}, \widetilde{V}, X, Z)$-coordinates
are translated into $(z, \tau, Y)$-coordinates
by replacing $Z$ by $1/z$ and $X$ by $z\,Y$,
substituting for $\widetilde{U}$ and $\widetilde V$ using \eqref{eq:tdef},
and omitting overall factors of $z$.
Thus, $\mathcal{C}$ is given by $Y=0$,
whereas the surfaces $Son$ and $Son'$ of Eqs.~\eqref{eq:eqson}
and~\eqref{eq:eqson'} are given by
\begin{align}
c\,z \left[(b_{1}+1)\,z^2 + 3\right] \tau
+ \left[(b_{1}+1)\,z^2 - 1\right] \frac{Y}{2}
+ \,\frac{c \left[(b_{1}-1)\,z^{2} + 1\right]}{z^2+1} &= 0,
\label{eq:son} \\
c\,z \left[(b_{1}+1)\,z^2 + 3\right] \tau
- \left[(b_{1}+1)\,z^2 - 1\right] \frac{Y}{2}
+ \frac{c \left[(b_{1}-1)\,z^{2} + 1\right]}{z^2+1} &= 0.
\label{eq:son'}
\end{align}
Likewise, Eq.~\eqref{eq:sp2} shows that the shock speed $\sigma$ is
\begin{equation}
\sigma=\frac{c}{b_1} \left[(b_1 + 1)\,z^2 - 1\right] \tau
+ \frac{c}{b_1}\,\frac{(b_1 + 2)\,z}{z^2 + 1} + \sigma_0.
\label{eq:speed}
\end{equation}

In the remainder of this paper,
we use $(z, \tau, Y)$-coordinates exclusively.
In effect, we excise the plane $\Pi$.
Consequently, some topological features are lost.
For example,
$\mathcal{W} \simeq M_2 \times \mathbb{R}$,
where $M_2$ is the M\"obius strip,
but $\mathcal{W} \backslash \Pi \simeq \mathbb{R}^3$;
similarly, $\mathcal{C}$ is a cylinder
but $\mathcal{C} \backslash \Pi$ is a plane.
For simplicity,
we will often write $\mathcal{W}$ rather than $\mathcal{W}\setminus\Pi$
and do likewise for surfaces and curves within $\mathcal{W}$.\break
XXXXX
\subsection{Hugoniot curves}
\label{subsec:hugoniot_curves}
Parametric equations for the Hugoniot curve with left state $W_0$
derive from \eqref{eq:KL},
which become linear equations for $\tau$ and $Y$ having solution
\begin{align}
\label{eq:parahug}
&z=z,\nonumber\\
&\tau= \dfrac{(\widetilde{V}_0+c)\,b_{1}z^{3}-\widetilde{U}_0\,z^{2}
+[b_{1}\,\widetilde{V}_0+2\,c]z-\widetilde{U}_0}
{c\,(z^{2}+1)[(b_{1}-1)\,z^{2}+1]},\\
&Y=\frac{-2\,(\widetilde{V}_0+c)\,z^{2}+2\,\widetilde{U}_0\,z
+2\,\widetilde{V}_0}{(b_{1}-1)\,z^{2}+1}.\nonumber
\end{align}
The corresponding parametrization of
the Hugoniot$'$ curve with right state $W_0$
is obtained from \eqref{eq:parahug} by replacing $Y$ with $-Y$.

More generally, the equations for a Hugoniot curve containing
the point $\mathcal{U}_0 = (z_0,\tau_0, Y_0)$ are derived as follows.
First, we set $(z, \tau, Y) = (z_0, \tau_0, Y_0)$
in Eqs.~\eqref{eq:parahug}:
\begin{align}
\label{eq:parahugpont}
&\tau_0= \dfrac{(\widetilde{V}_0+c)\,b_1\,z_0^{3}-\widetilde{U}_0\,z_0^2
+[b_1\,\widetilde{V}_0+2\,c]\,z_0-\widetilde{U}_0}
{c\,(z_0^2+1)[(b_1-1)\,z_0^2+1]},\\
&Y_0=\frac{-2\,(\widetilde{V}_0+c)\,z_0^2+2\,\widetilde{U}_0\,z_0
+2\,\widetilde{V}_0}{(b_1-1)z_0^2+1}.\nonumber
\end{align}
Solving \eqref{eq:parahugpont}
for $\widetilde{U}_0$ and $\widetilde{V}_0$
and substituting into Eqs.~\eqref{eq:parahug} yields
\begin{equation}
\begin{array}{l}
z=z,\quad
\tau_{hug(z)}=\dfrac{Dz^{3}+Ez^2+Fz+G}
{\vartheta_0 c[(b_1-1)z^{4}+b_1z^2+1]},\quad
Y_{hug(z)}=\displaystyle \frac{Az^2+Bz+C}{\vartheta_0[(b_1-1)z^2+1]},
\end{array}
\label{eq:hugponto}
\end{equation}
where $A=-[2\,c+\left(2\,c\,\tau_0\,z_0+Y_0\right)\vartheta_0]$,
$B=4\,c\,z_0+2\,c\,\tau_0\,(z_0^2-1)\,\vartheta_0^2+b_1\,Y_0\,z_0\vartheta_0$,
$C=-2\,c\,z_0^2+(2\,c\,\tau_0\,z_0+Y_0)\vartheta_0$,
$D=2\,c\,b_1+\left(2\,c\,b_1\,\tau_0\,z_0+b_1\,Y_0\right)
\vartheta_0$,
$E=-2\,c\,\tau_0\,(z_0^2-1)\,\vartheta_0^2
-b_1\,z_0\,Y_0\,\vartheta_0-4\,c\,z_0$,
$F=\left(4\,c+2\,c\,b_1\,\tau_0\,z_0+b_1\,Y_0\right)\vartheta_0
-2\,c\,b_1\,z_0^2$,
$G=-4\,c\,z_0-2\,c\,\tau_0\,(z_0^2-1)\,\vartheta_0^2
-b_1\,z_0\,Y_0\,\vartheta_0$,
and $\vartheta_0=z_0^2+1$.

\subsection{Coincidence curve}
\label{subsec:coincidence_curve}
Because a Hugoniot curve intersects $\mathcal{C}$ when $Y = 0$,
the numerator in the expression for $Y$ in~\eqref{eq:parahug}
determines whether such an intersection exists.
This quadratic polynomial in $z$ has discriminant
$\Delta_0 = \widetilde{U}_0^{2}+4\,\widetilde{V}_0\,(\widetilde{V}_0+c)$.
Therefore,
the Hugoniot curve with left state $W_0$ intersects $\mathcal{C}$
at two points if $\Delta_0 > 0$ and a single point if $\Delta_0 = 0$,
whereas it does not intersect $\mathcal{C}$ if $\Delta_0 < 0$.
As $(u_0, v_0)$ is related to $(\widetilde{U}_0, \widetilde{V}_0)$
by an affine linear transformation,
the solution set of $\Delta_0 = 0$
is an ellipse in the $(u_0, v_0)$-plane.
(We have assumed that $c > 0$.)

\begin{definition}
The \emph{coincidence curve} $E_0$ in state space
is the set of states such that $\Delta_0 = 0$.
The \emph{strictly hyperbolic region} comprises
states for which $\Delta_0 > 0$,
and the \emph{elliptic region} is where $\Delta_0 \le 0$.
\end{definition}

Similarly, the numerator in the expression for $Y$ in \eqref{eq:hugponto}
determines the intersections
of $\mathcal{H}(\mathcal{U}_0)$ with $\mathcal{C}$:
it consists of two distinct points,
a single point, or no points according to whether the discriminant
$B^{2} - 4\,A\,C$ is positive, zero, or negative.
To find the coordinates $(z, \tau, 0)$ of an intersection point,
we solve $Az^{2}+Bz+C=0$ for $z$
and calculate $\tau$ using the second equation in \eqref{eq:hugponto}.

Suppose that $Y_0 = 0$,
\emph{i.e.}, $\mathcal{U}_0 \in \mathcal{C}$;
let $W_0$ be the corresponding state.
Then Eqs.~\eqref{eq:hugponto} coincide with Eqs.~\eqref{eq:parahug},
and comparing the respective formulae for $Y$ shows that
$B^2 - 4\,A\,C = \Delta_0\,\vartheta_0^2$.
Also, straightforward computation shows that
$B^{2} - 4\,A\,C = 4\,c^{2}\,\tau_0^{2}\,\vartheta_0^{4}$.
Therefore,
$\tau_0 = 0$ if and only if $W_0$ belongs to
the coincidence curve $E_0$ in state space,
and otherwise $W_0$ lies in the strictly hyperbolic region.

\begin{definition}
\label{coincidence}
The \emph{coincidence curve} $\mathcal{E}$
is the line where $\tau = 0$ and $Y = 0$.
\end{definition}

If $\tau_0 = 0$,
then $\mathcal{H}(\mathcal{U}_0)$ intersects $\mathcal{C}$
only at $\mathcal{U}_0$, which belongs to $\mathcal{E}$.
By computing $dY/dz$ from \eqref{eq:hugponto}
and setting $\tau_0 = 0$ and $Y_0 = 0$,
we find that $dY/dz = 0$,
confirming that if $\mathcal{U}_0 \in \mathcal{E}$,
then $\mathcal{H}(\mathcal{U}_0)$
is tangent to $\mathcal{C}$ at $\mathcal{U}_0$.
If, on the other hand, $\tau_0 \ne 0$,
$\mathcal{H}(\mathcal{U}_0)$ intersects $\mathcal{C}$ at a second point.
A simple calculation shows that $dY/dz \ne 0$ at these intersections,
\emph{i.e.}, $\mathcal{H}(\mathcal{U}_0)$
intersects $\mathcal{C}$ transversally.

We denote the region where $\tau < 0$ by $\mathcal{C}_s$
and the region where $\tau > 0$ by $\mathcal{C}_f$.
The subscripts $s$ and $f$ stands for ``slow'' and ``fast'',
in accordance with the following result,
which is proved in \ref{sec:Proofs_of_Results}.

\begin{proposition}
If $(z_0,\tau_0,0)$ and $(z_1,\tau_1,0)$
are distinct intersections of a Hugoniot curve with $\mathcal{C}$,
then one of the intersection points belongs to $\mathcal{C}_s$
and the other belongs to $\mathcal{C}_f$.
Furthermore, the value of $\sigma$ at the point in $\mathcal{C}_s$
is smaller than the value at the point in $\mathcal{C}_f$.
\label{th:slow_fast}
\end{proposition}

\section{Decomposition of the Wave Manifold}
\label{sec:Decomposition_of_the_Wave_Manifold}
In this section, we describe how the \emph{characteristic},
\emph{sonic} and \emph{sonic}$'$ surfaces subdivide $\mathcal{W}$.
We also decompose \emph{sonic}$'$ into $Son'_s$ and $Son'_f$.

\subsection{Subdivision of the wave manifold}
\label{sec:Subdivision_of_the_wave_manifold}
\begin{proposition}
The characteristic, sonic, and sonic\,$'$ surfaces
subdivide $\mathcal{W} \backslash \Pi$ into twelve regions.
\label{proposition:divideM}
\end{proposition}

\begin{figure}[!ht]
\begin{center}
\includegraphics[scale=1,width=0.75\linewidth]
{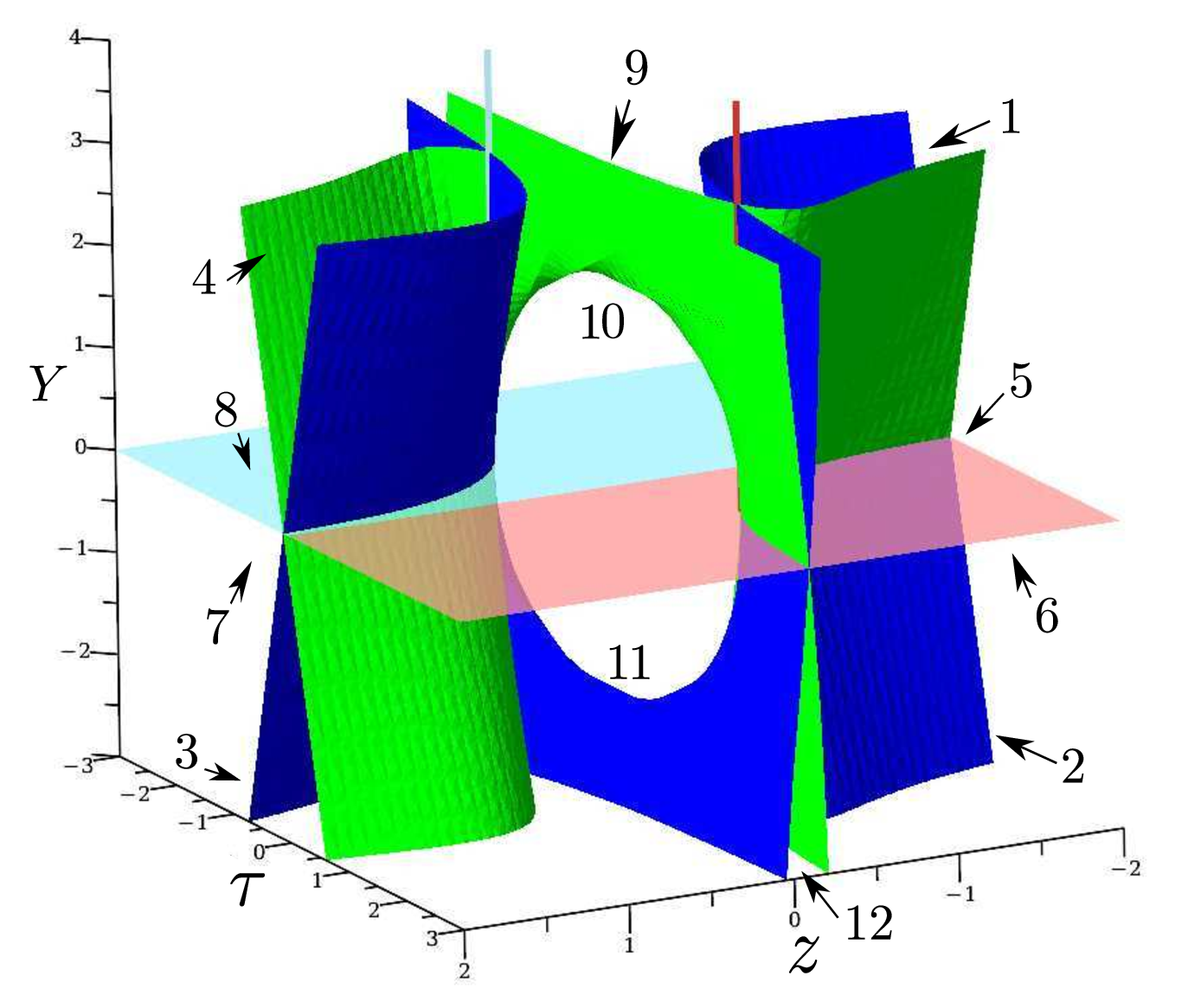}
\caption[]{
Decomposition of $\mathcal{W}$ by the surfaces
$\mathcal{C} = \mathcal{C}_{s} \cup \mathcal{C}_{f}$
(the plane colored light blue in $\mathcal{C}_{s}$
and light coral in $\mathcal{C}_{f}$),
$Son$ (the green surface),
and $Son'$ (the blue surface).
There are twelve regions.
All three surfaces, intersect at the inflection locus,
a hyperbola-like curve on $\mathcal{C}$
also shown in Fig.~\ref{fig:02raref}.
In this and all subsequent figures, $b_1 = 8$ and $c = 1$.
However, the decomposition is qualitatively the same
for any $b_1 > 1$ and $c > 0$.
}
\label{fig:02}
\end{center}
\end{figure}

The proof of Proposition \ref{proposition:divideM}
is in \ref{sec:Proofs_of_Results},
and the subdivision is depicted in Fig.~\ref{fig:02}.

Let us describe this subdivision.
Examining Fig.~\ref{fig:02},
we notice symmetry between half-spaces $Y>0$ and $Y<0$.
We focus on the upper half-space.
Here, the blue surface ($Son'$) has two connected components,
subdividing the half-space into three regions: left, middle, and right.
The middle region is subdivided further into two by the green surface ($Son$).
One is Region~9, and the other is Region~10,
bounded by parts of both $Son$ and $Son'$.

From Eqs.~\eqref{eq:son} and~\eqref{eq:son'} we deduce that
the intersection of the green and blue surfaces ($Son$ and $Son'$)
is as follows:
\begin{itemize}
\item two branches of a curve in the plane $Y=0$,
given by $\tau=\tau_{\text{infl}}(z)$, where
\begin{equation}
\tau_{\text{infl}}(z)=\frac{-\left[(b_{1}-1)z^{2}+1\right]}
{z\,(z^2 + 1) \left[(b_{1}+1)\,z^2 + 3\right]}
\label{eqinfle}
\end{equation}
for $z \ne 0$,
which form the \emph{inflection locus} $\mathcal{I}$,
illustrated as well in Fig.~\ref{fig:02raref} below;
\item two lines where
\begin{equation}
\text{$(b_1 + 1)\,z^2 - 1 = 0$,
\quad $\tau = \tau_{\text{infl}}(z)$,
\quad and \quad $Y$ is arbitrary,}
\label{eq:double_sonic}
\end{equation}
which constitute the \emph{double sonic locus} $\mathcal{D}$.
\end{itemize}
Parts of these intersecting surfaces bound Regions~1 (right) and~4 (left),
while other parts bound Regions~5 and 8.

In the half-plane $Y<0$
we have an analogous decomposition obtained by reflection,
which interchanges $Son$ and $Son'$.

Recall that the coincidence curve $\mathcal{E}$ in $\mathcal{C}$
separates $\mathcal{C}_{s}$ and $\mathcal{C}_{f}$.
As shown at the end of Sec.~\ref{sec:Wave_Manifold_and_Hugoniot_Curves},
Hugoniot curves are tangent to $\mathcal{C}$
along the coincidence curve $\mathcal{E}$.
Such Hugoniot curves are also tangent to $Son'$;
see \cite{eschenazi02}.

\begin{definition}
Hugoniot curves emanating from points along $\mathcal{E}$
form a two-dimensional submanifold of $\mathcal{W}$
called the \emph{saturation of the coincidence curve},
denoted $SCC$.
Similarly, $SCC\,'$ is generated by Hugoniot$'$ curves
of points along $\mathcal{E}$.
\end{definition}

In other words, a point $\mathcal{U}$ belongs to $SCC$
if and only if $\mathcal{H}(\mathcal{U})$ intersects $\mathcal{C}$ once,
namely in $\mathcal{E}$.
Consequently:

\begin{proposition}
The $SCC$ separates $\mathcal{W}$ into two regions:
\begin{itemize}
  \item points $\mathcal{U}$ for which $\mathcal{H}(\mathcal{U})$
  does not intersect $\mathcal{C}$; and
  \item points $\mathcal{U}$ for which $\mathcal{H}(\mathcal{U})$
  intersects $\mathcal{C}$ at two distinct points.
\end{itemize}
\end{proposition}

\begin{remark}
Topologically, the surface $SCC$
is a cylinder within $\mathcal{W}\backslash\Pi$,
each constant-$z$ slice being an ellipse
tangent to the $\tau$-axis at $(0,0)$
and lying within the half-plane $Y\leq 0$.
As a result,
$SCC$ is entirely contained in the half-space $Y\leq 0$,
specifically Region~11 in Fig.~\ref{fig:02}.
A Hugoniot curve that does not intersect $\mathcal{C}$
lies inside the open three-dimensional region bounded by the $SCC$,
whereas a Hugoniot curve that intersects $\mathcal{C}$ twice
lies outside.
\end{remark}
\par
\begin{figure}[!ht]
\begin{center}
\includegraphics[scale=1,width=0.75\linewidth]
{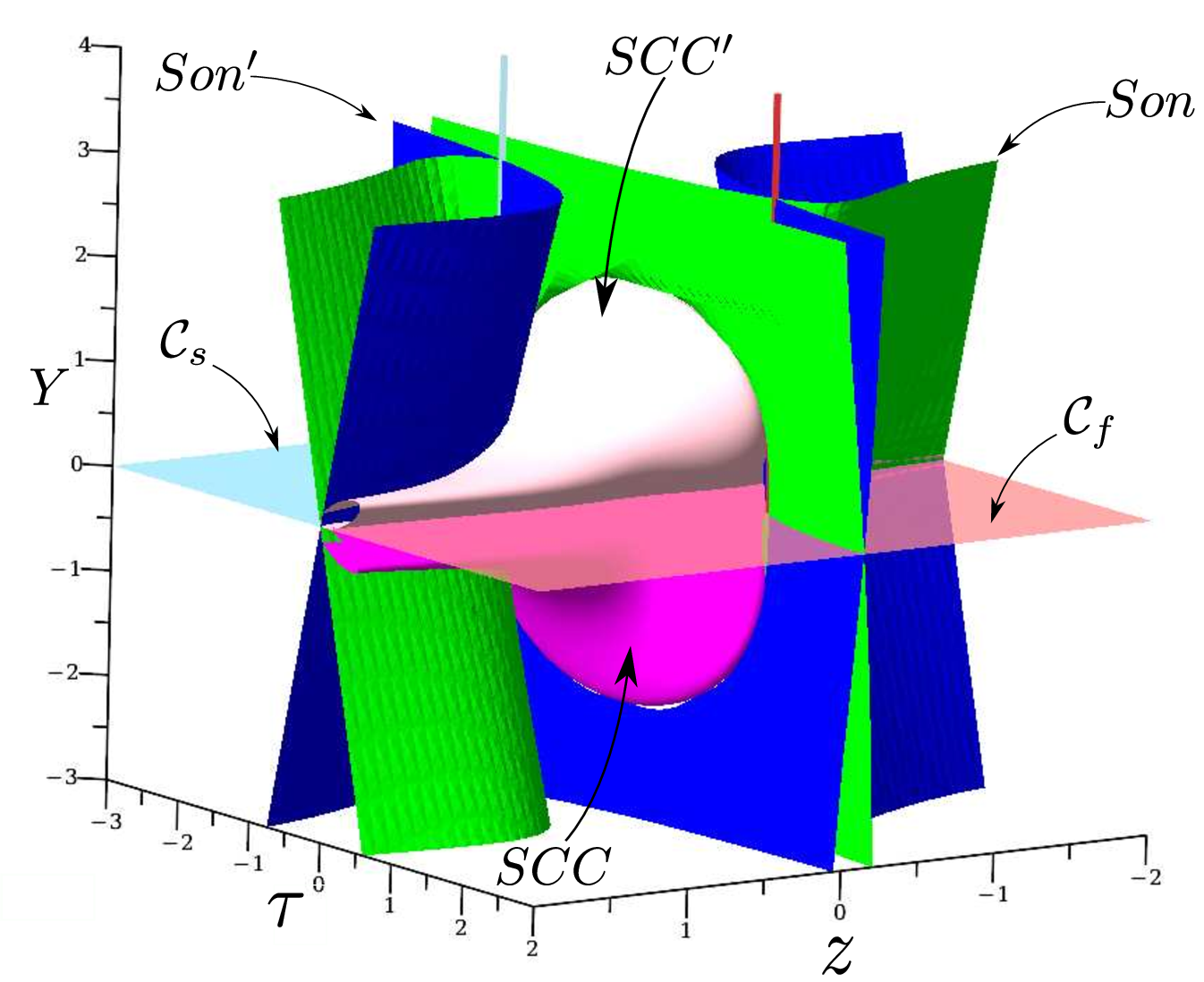}
\caption[]{
The surfaces $\mathcal{C}$ (light blue and light coral),
$Son$ (green), $Son'$ (blue),
$SCC$ (magenta), and $SCC\,'$ (pink).
}
\label{fig:02SCC}
\end{center}
\end{figure}

\noindent
The $SCC$ and $SCC'$ are depicted as the magenta and pink surfaces
in Fig.~\ref{fig:02SCC}.

\begin{remark}
\label{rem:SCC}
Adding the $SCC$ and $SCC\,'$ surfaces to Fig.~\ref{fig:02}, we obtain
Fig.~\ref{fig:02SCC}, where region 10 is divided into 2 sub-regions, one with a
part of $\mathcal{C}_s$ as lower boundary and the other with a part of
$\mathcal{C}_f$ as lower boundary. The same is true for Region 11. There are
now fourteen regions, these regions are studied in Sec.~4.3.There it is shown that each of these fourteen regions has uniform
\emph{shock type} (as defined in Sec.~\ref{sec:lax}),
which changes when a boundary is crossed.
\end{remark}

\begin{definition}\label{def:ECC'}
The intersection of $SCC$ and $Son'$ is called
the \emph{extension of the coincidence curve} by Hugoniot curves
and is denoted by $ECC\,'$.
\end{definition}
\par\noindent
The curve $ECC\,'$ was referred to as
the \emph{sonic fold} in \cite{eschenazi02}.
The following result is proved in \ref{sec:Proofs_of_Results}.

\begin{proposition}
The surface $SCC$ is tangent to $\mathcal{C}$ along $\mathcal{E}$
and to $Son'$ along $ECC\,'$.
\label{th:Tf}
\end{proposition}

\subsection{Finding the surfaces \texorpdfstring{$Son'_{s}$}{Son'\_s}
and \texorpdfstring{$Son'_{f}$}{Son'\_f}\label{sec:sonlifs}}
In this subsection,
we describe how the surface $Son'$ is subdivided into two parts,
called the slow sonic$'$ surface and the fast sonic$'$ surface.

\begin{proposition}
\label{propson}
Suppose that $\mathcal{U}_0=(z_0,\tau^\prime_0,Y_0) \in Son'$.
Then $\mathcal{H}(\mathcal{U}_0)$ intersects $\mathcal{C}$
at a point for which the shock speed equals $\sigma(\mathcal{U}_0)$.
\end{proposition}

In the proof, we take the following steps.
\begin{itemize}
\item Calculate $\tau'_0$ in terms of $z_0$ and $Y_0$
using \eqref{eq:son'}.
\item Find the shock speed of $\mathcal{U}_0$ using \eqref{eq:speed}.
\item Parametrize the Hugoniot curve of $\mathcal{U}_0$
by $z$ using \eqref{eq:parahugpont}.
\item Derive the quadratic equation for $z$ such that $Y = 0$
and explicitly find its roots.
\item Calculate the values of $\tau$
corresponding to these roots using \eqref{eq:hugponto}.
\item Calculate the shock speeds at these intersection points.
\item Observe that one of these shock speeds equals $\sigma(\mathcal{U}_0)$.
\end{itemize}

\begin{definition}
The \emph{slow sonic\,$'$} surface, denoted by $Son'_{s}$,
is the set of points $\mathcal{U} \in Son'$ such that
$\sigma(\mathcal{U})=\sigma(\mathcal{U}_{s})$.
Similarly, the \emph{fast sonic\,$'$} surface,
denoted by $Son'_{f}$, is the set of points $\mathcal{U} \in Son'$
such that $\sigma(\mathcal{U})=\sigma(\mathcal{U}_{f})$.
\end{definition}

\begin{proposition}
\label{proposition:son'sson'f}
The curve $ECC\,'$ of Def.~\textup{\ref{def:ECC'}}
is the boundary between $Son'_{s}$ and $Son'_{f}$.
\end{proposition}

\begin{figure}[!ht]
\begin{center}
\begin{subfigure}{0.45\linewidth}
\includegraphics[width=\linewidth]{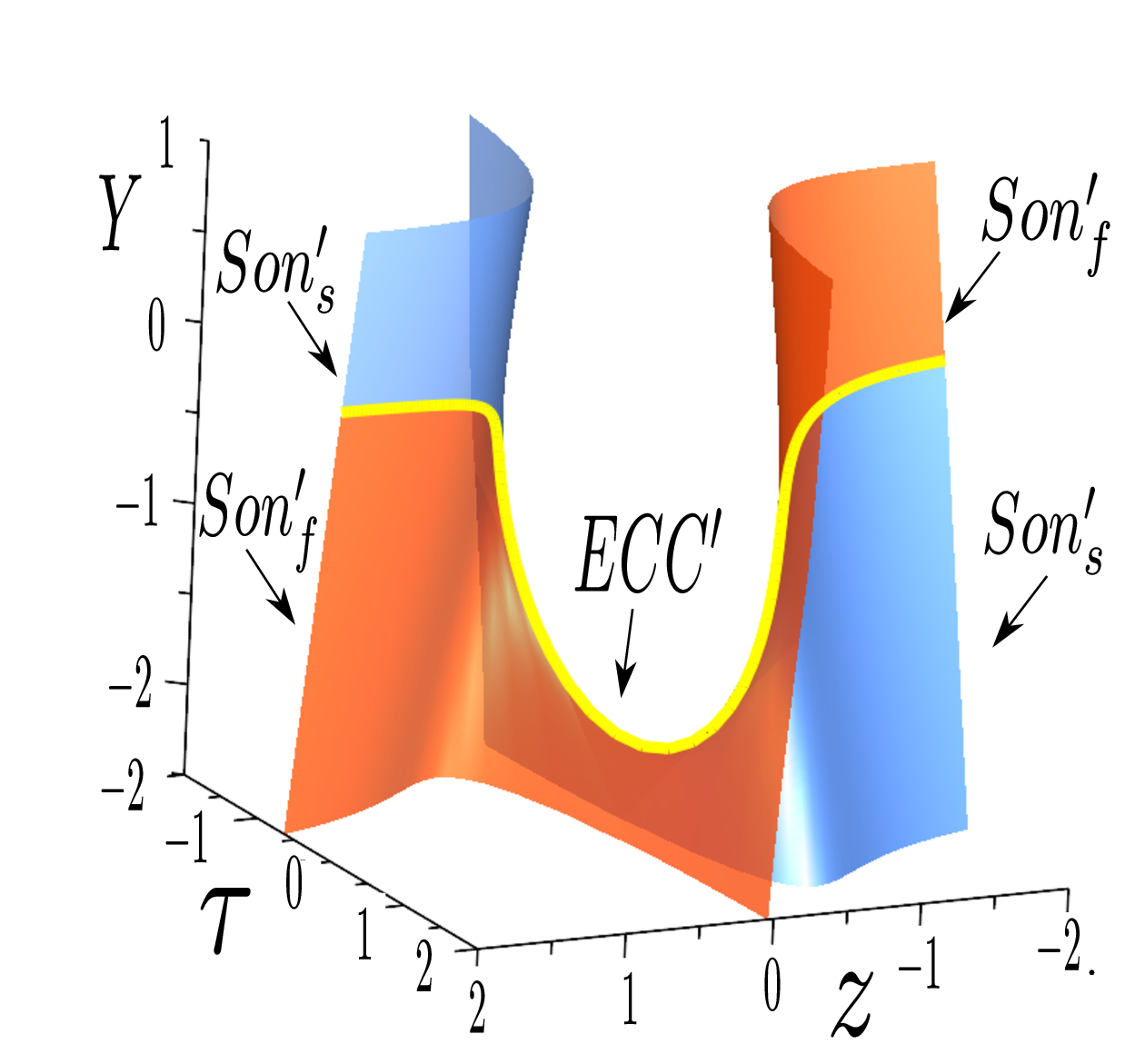}
\caption{}
\end{subfigure}
\hfil
\begin{subfigure}{0.45\linewidth}
\includegraphics[width=\linewidth]{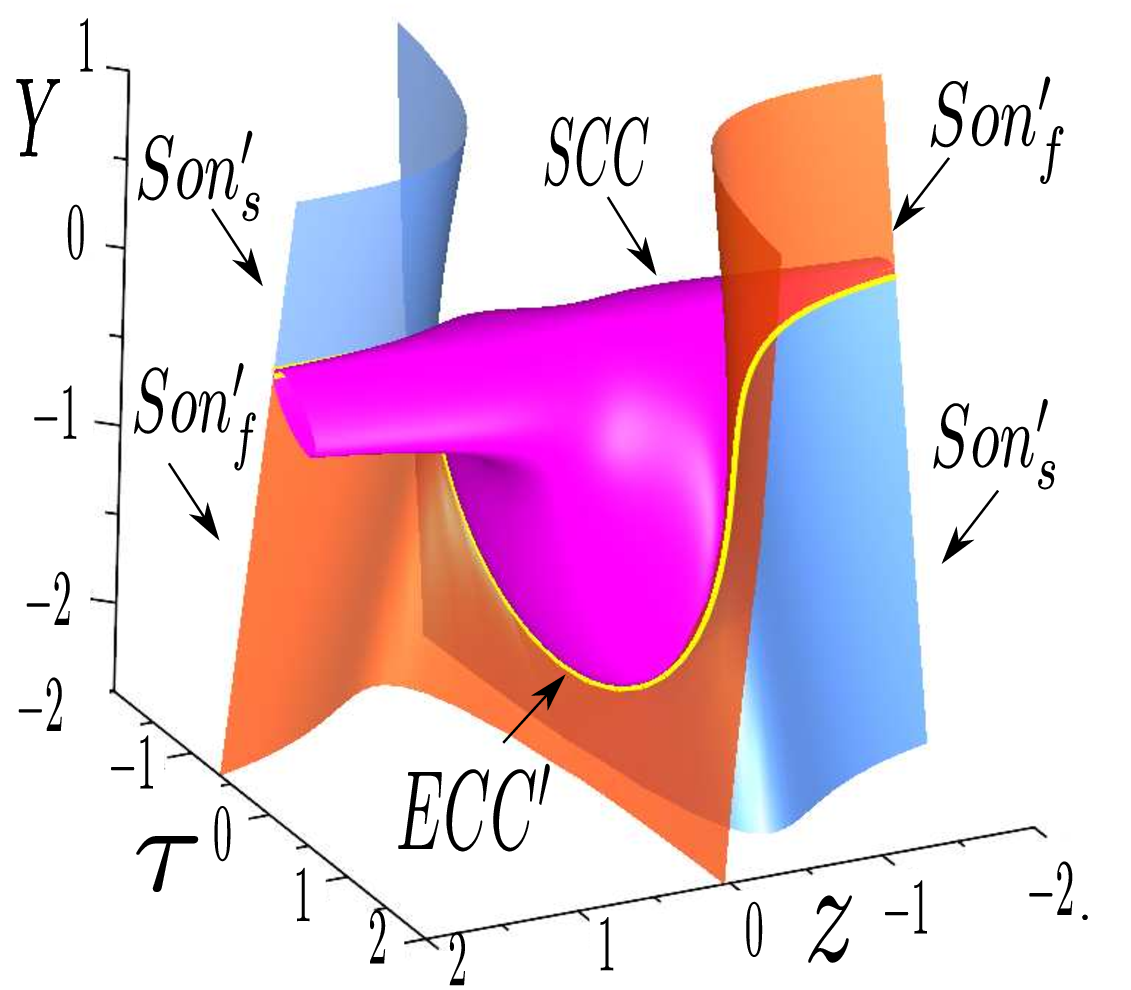}
\caption{}
\end{subfigure}
\caption[]{
(a) The $Son'$ surface subdivided in two surfaces,
$Son'_{s}$ (light blue) and $Son'_{f}$ (light orange),
by the $ECC\,'$ curve (yellow).
(b) The $ECC\,'$ curve is the intersection of $SCC$ (magenta)
and $Son'$, which are tangent.
}
\label{fig:son-linha-s-f}
\end{center}
\end{figure}

Proofs of Propositions \ref{propson} and \ref{proposition:son'sson'f}
can be found in \ref{sec:Proofs_of_Results}.
Refer to Fig.~\ref{fig:son-linha-s-f} for an illustration.

\section{Rarefaction, Shock, and Composite Curves}
\label{sec:Rarefaction_Shock_and_Composite_Curves}

This section describes three fundamental types
of curves for constructing the RPS:
rarefaction, shock, and composite curves.

\subsection{Rarefaction curves}
\label{subsec:raref}
We now examine smooth solutions of $W_t + F(W)_x = 0$
in the form $\widetilde{W}(x/t)$,
which are called \emph{rarefaction waves}.
Such solutions satisfy
\begin{equation}
(-x/t^{2})\,d\widetilde{W}/d\lambda
+(1/t)\,DF(\widetilde{W})\cdot d\widetilde{W}/d\lambda=0,
\label{rareffor}
\end{equation}
where $\lambda=x/t$, $\widetilde{W} = \widetilde{W}(\lambda)$,
and $DF(\widetilde{W})$ is the Jacobian derivative of the flux $F$
\cite{smoller94}.

Multiplying \eqref{rareffor} by $t$ and replacing $x/t$ by $\lambda$,
we deduce the requirement
$DF(\widetilde{W})\cdot d\widetilde{W}/d \lambda
= \lambda\,d \widetilde{W}/d\lambda$,
which is the differential equation for the line field
of eigenspaces of $DF$.
Because $\widetilde{W}(x/t)$ is constant
along space-time trajectory $x = \lambda\,t$,
the eigenvalue $\lambda$ is a wave propagation speed,
known as a \emph{characteristic speed}
(which is why $\mathcal{C}$ is called the characteristic manifold).
Suitable arcs of integral curves of this line field
are called \emph{rarefaction curves}.
A comprehensive study of rarefaction curves in quadratic models
can be found in \cite{palmeira88} and \cite{bastos05}.
A general definition of rarefaction curves in the wave manifold
can be found in \cite{isaacson92}.

For the function $F$ defined by Eqs.~\eqref{eq:f}--\eqref{eq:g} with $b_2 = 0$,
$DF$ is
\begin{align}
DF=\begin{bmatrix}
(b_{1}+1)\,u+a_{1}& v+a_{2}\\
v+a_{3}&u+a_{4}\\
\end{bmatrix}.
\end{align}
Let $\vec{r}=(r_{1},r_{2})^T$ be
the right eigenvector associated to the eigenvalue $\lambda$:
$(\lambda\mathbb{I}-DF)\,\vec{r}=0$.
Eliminating $\lambda$ and replacing $r_{1}/r_{2}$ by $du/dv$,
we arrive at the following differential equation:
\begin{equation}
(v+a_{3})\,(du/dv)^{2} - (b_{1}\,u+a_{1}-a_{4})\,du/dv - (v+a_{2}) = 0.
\label{eq:raref}
\end{equation}

Let $\Delta$ be the discriminant of this quadratic equation in $du/dv$:
\begin{equation}
\Delta=(b_{1}\,u+a_{1}-a_{4})^{2} +4(v+a_{2})\,(v+a_{3}).
\end{equation}
Notice that $\Delta$ coincides with the expression $\Delta_0$
defined in Sec.~\ref{subsec:coincidence_curve}
if we replace $u_0$ and $v_0$ by $u$ and $v$.
Therefore, the curve $\Delta=0$ is
the coincidence locus $E_0$ in state space,
which is an ellipse.
We see that within the elliptic region (the interior of the ellipse),
there are no rarefaction curves,
whereas at each point of the strictly hyperbolic region
(the exterior of the ellipse),
two rarefaction curves cross transversely.

The differential equation \eqref{eq:raref}
is not of the usual form in which $du/dv$ is specified
as a function of $u$ and $v$.
Instead of solving the quadratic equation for $du/dv$,
we use the technique due to Lagrange,
of replacing $du/dv$ by a variable $z$,
so that \eqref{eq:raref} becomes a simple differential equation
on an implicitly-defined surface in $(u, v, z)$-space:
\begin{equation*}
\text{$z\,dv - du = 0$ on $G_2 = 0$},
\end{equation*}
where $G_2 = (v + a_3)\,z^{2} - (b_1\,u + a_1 - a_4)\,z - (v + a_2)$.

When expressed in terms of the variables
$\widetilde{U} = b_1\,u + a_1 - a_4$ and $V_1 = v + a_3$,
$G_2$ is identical to the function $G$ of Eq.~\eqref{G}.
Therefore, the surface $G_2 = 0$ can be identified
as the $Y = 0$ plane within $\mathcal{W}$,
\emph{i.e.}, as $\mathcal{C}$.
In other words,
rarefaction curves are integral curves of
$b_1\,z\,dV_1 - d\widetilde{U} = 0$ and $dG = 0$.
These equations are linearly independent away from $z = \infty$.
Thus, the family of rarefaction curves foliate $\mathcal{C} \setminus \Pi$.

We can solve $(\lambda\mathbb{I}-DF)\,\vec{r}=0$
to obtain the eigenvalue $\lambda$ as a function of $u$, $v$ and $z=r_{1}/r_{2}$.
If we replace $u$ and $v$ by the coordinates $\widetilde{U}$ and $V_{1}$,
we see that $\lambda$ is given by the
same formula \eqref{eq:sp2} as is $\sigma$.
Thus, the shock speed $\sigma$ and the characteristic speed $\lambda$ coincide
on the characteristic manifold.
Specifically, if $\mathcal{U}_s\in\mathcal{C}_{s}$
and $W$ is its associated state,
then the smaller eigenvalue of $DF(W)$
is $\lambda_s(W)=\sigma(\mathcal{U}_{s})$.
Similarly, if $\mathcal{U}_f\in\mathcal{C}_{f}$ has associated state $W$,
then the larger eigenvalue of $DF(W)$ is $\lambda_f(W)=\sigma(\mathcal{U}_{f})$.

Let us derive the differential equation for rarefaction curves
in $\tau$ and $z$ coordinates.
By taking the expressions for $\widetilde{U}$ and $V_1$ from \eqref{eq:tdef}
and calculating $d\widetilde{U}$ and $dV_1$ in terms of $d\tau$ and $dz$,
$b_{1}\,z\,dV_{1}-d\widetilde{U}=0$ becomes
\begin{equation}
\frac{d\tau}{dz}
= \frac{-(b_1 - 2)\,z}{(b_1 - 1)\,z^2 + 1}\,\tau + \frac{2}{(z^2 + 1)^2}.
\label{eq:edoraref}
\end{equation}

Rarefaction waves are constructed from arcs of rarefaction curves
along which $\sigma$
(\emph{i.e.}, the characteristic speed $\lambda$)
varies monotonically.
(By \emph{arc} we mean a homeomorphic image
of an interval within $\mathbb{R}$.)

\begin{definition}\label{forwardr}
A \emph{slow rarefaction curve},
denoted by $\mathcal{R}_s$,
is an arc of a rarefaction curve within $\mathcal{C}_{s}$
along which the speed $\sigma$ is strictly monotone.
Similarly, a \emph {fast rarefaction curve}, denoted by $\mathcal{R}_f$,
is an arc of a rarefaction curve within $\mathcal{C}_{f}$
along which $\sigma$ is strictly monotone.
\end{definition}

\begin{figure}[!ht]
\begin{center}
\includegraphics[scale=0.9,width=0.85\linewidth]{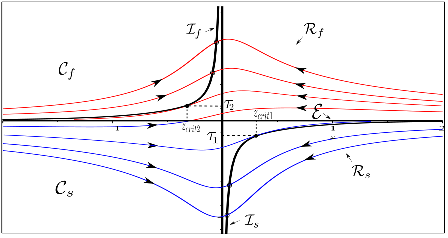}
\caption[]{
Rarefaction curves $\mathcal{R}_{s}$ and $\mathcal{R}_{f}$
in the characteristic plane $\mathcal{C}$.
The arrows indicate the direction
in which the speed $\sigma$ increases for $\mathcal{R}_s$
and decreases for $\mathcal{R}_f$.
Also shown are the inflection loci,
$\mathcal{I}_{s}$ in the quadrant where $z>0$ and $\tau<0$,
and $\mathcal{I}_{f}$ in the quadrant where $z<0$ and $\tau>0$.
At the points $z_{crit1}$ and $z_{crit2}$,
which are defined in Eqs.~\eqref{eq:zcrit1} and~\eqref{eq:zcrit2},
the rarefaction curve is tangent to $\mathcal{I}$.
}
\label{fig:02raref}
\end{center}
\end{figure}

The expression for $\sigma$ in terms of $z$ and $\tau$
is given by Eq.~\eqref{eq:speed}.
Calculating
$d\sigma/dz=\partial \sigma/\partial z
+(\partial \sigma/\partial \tau)\,(d\tau/dz)$,
where $d\tau/dz$ is given by \eqref{eq:edoraref},
we find that
\begin{equation}
\frac{d\sigma}{dz}
= \dfrac{c\,z \left[(b_1 + 1)\,z^2 + 3\right]}{(b_1 - 1)\,z^2 + 1}\,\tau
+ \frac{c}{z^2 + 1}.
\end{equation}
Observe that $d\sigma/dz$ takes the form
$z\,r(z) \left[\tau - \tau_{\text{infl}}(z)\right]$,
where $\tau_{\text{infl}}$ is given by \eqref{eqinfle}
and $r(z)$ is a positive rational function of $z$.
Therefore, the inflection locus $\mathcal{I}$,
defined above as $\mathcal{C} \cap Son = \mathcal{C} \cap Son'$,
is also where $\sigma$ has a critical point along a rarefaction curve,
as is true generally \cite{isaacson92}.
For the model studied in this paper,
it is the hyperbola-like curve shown in Fig.~\ref{fig:02raref}.
We denote by $\mathcal{I}_{s}$ and $\mathcal{I}_{f}$
the branches of $\mathcal{I}$ within $\mathcal{C}_{s}$ and $\mathcal{C}_{f}$,
respectively.
The sign of $d\sigma/dz$ is positive between the two branches
and negative outside them.

By differentiating the formula for $d\sigma/dz$
and evaluating at the inflection locus,
we calculate that
\begin{equation}
\frac{d^2\sigma}{dz^2}
= z\,r(z) \left(\frac{d\tau}{dz} - \frac{d\tau_{\text{infl}}}{dz}\right)
= \frac{3\,c \left[(b_1 + 1)\,z^2 - 1\right]}
{z \left[(b_1 - 1)\,z^2 + 1\right] \left[(b_1 + 1)\,z^2 + 3\right]}.
\end{equation}
Consequently, when $z > 0$,
$\sigma$ has a strict minimum along a rarefaction curve
at an intersection with $\mathcal{I}$
if and only if the slope of the rarefaction curve
exceeds that of the inflection locus,
and vice versa when $z < 0$.
Moreover, a rarefaction curve and $\mathcal{I}$ are tangent
if and only if $z$ satisfies $(b_1 + 1)\,z^2 - 1 = 0$
at their intersection,
which is where the double sonic locus $\mathcal{D}$,
given by~\eqref{eq:double_sonic},
meets the characteristic manifold $\mathcal{C}$.
(This feature also holds more generally~\cite{isaacson92}.)

\begin{proposition}\label{th:rarefp}
The slow rarefaction curve
of $(z_{\text{crit}_1}, \tau_1)\in \mathcal{C}_{s}$
is tangent to $\mathcal{I}_{s}$,
where
\begin{equation}
\text{$z_{\text{crit}_1} = 1/\sqrt{b_{1}+1}$
and $\tau_1 = -b_1\,\sqrt{b_1 + 1} / [2\,(b_1 + 2)]$.}
\label{eq:zcrit1}
\end{equation}
Similarly, the fast rarefaction curve
of $(z_{\text{crit}_2}, \tau_{2})\in \mathcal{C}_{f}$
is tangent to $\mathcal{I}_{f}$,
where
\begin{equation}
\text{$z_{\text{crit}_2} = -1/\sqrt{b_{1}+1}$
and $\tau_1 = b_1\,\sqrt{b_1 + 1} / [2\,(b_1 + 2)]$.}
\label{eq:zcrit2}
\end{equation}
\end{proposition}

\subsection{Lax shock curves in
\texorpdfstring{$\mathcal{W}$}{W} }
\label{sec:lax}
This subsection discusses the shock wave admissibility condition
we adopt for arcs of Hugoniot and Hugoniot$'$ curves within $\mathcal{W}$.

Beyond satisfying the Rankine-Hugoniot condition,
shock waves that appear in solutions of Riemann problems
must obey an \emph{admissibility criterion}
that reflects their dissipative nature.
Lax~\cite{lax} abstracted this criterion for the shock wave~(\ref{eq:shock})
as inequalities relating the shock speed $\sigma$
to the characteristic speeds $\lambda_s(W)$, $\lambda_f(W)$,
$\lambda_s(W')$, and $\lambda_f(W')$ of its left and right states.
We refer to the particular ordering of these speeds
as the \emph{type} of the shock wave.

For a system of two conservation laws,
two types of shock waves
are admissible according to the Lax criterion
(inequalities~\eqref{lax1} and~\eqref{lax2} explained below);
we call them \emph{slow} and \emph{fast Lax shock waves}.
For the specific model examined in this paper,
these types of shock waves
are deemed suitable for use in Riemann solutions,
whereas other types of shock waves are excluded.

\begin{remark}
We adopt the Lax admissibility criterion primarily for simplicity.
As shown in~\cite{isaacson88}
solutions of Riemann problems exist and are unique
under the Lax admissibility criterion
for symmetric Case~IV quadratic models that are homogeneous,
for which the elliptic region is a point.
However, for Case~I--III homogeneous quadratic models,
existence and uniqueness of the general RPS
fails for the Lax admissibility criterion
but holds for the \emph{viscous profile admissibility criterion}%
~\cite{schshe87a,isaacson88,isatem88,isatem88a,isamarplo90},
for which some Lax shock waves are inadmissible
and certain shock waves of another type are admissible.
Existence and uniqueness has not been established
for quadratic models with elliptic regions that are not a point.
\end{remark}

In order for the Lax admissibility criterion to be applicable,
we adopt the following restriction.

\begin{restrict}
\label{restrict:hyperbolic}
We only treat points $\mathcal{U} \in \mathcal{W}$
for which $\mathcal{H}(\mathcal{U})$ intersects $\mathcal{C}$ at two points
and $\mathcal{H}'(\mathcal{U})$ likewise intersects $\mathcal{C}$ at two points.
As explained in Remark~\ref{rem:SCC},
$\mathcal{U}$ lies outside both the $SCC$ and the $SCC'$;see Fig.\ref{fig:02SCC}. 
\end{restrict}
An arc of a Hugoniot curve along which the Lax admissibility criterion
is satisfied is referred to as a \emph{slow} or \emph{fast shock curve}
and is denoted by $\mathcal{S}_s$ or $\mathcal{S}_f$.
For a point $\mathcal{U} \in \mathcal{W}$
satisfying Restriction~\ref{restrict:hyperbolic}, we define
\begin{alignat}{2}
&\mathcal{U}_{s}=\mathcal{H}(\mathcal{U}) \cap \mathcal{C}_{s},
\qquad &&\mathcal{U}_{f}=\mathcal{H}(\mathcal{U}) \cap \mathcal{C}_{f},
\label{usuli} \\
&\mathcal{U}_{s}^{\prime} = \mathcal{H}'(\mathcal{U}) \cap \mathcal{C}_{s},
\qquad &&\mathcal{U}_{f}^{\prime}
= \mathcal{H}'(\mathcal{U}) \cap \mathcal{C}_{f}.
\label{usuli'}
\end{alignat}
Recall that a point $\mathcal{U}\in \mathcal{W}$
has an associated left state $W=(u,v)$ and right state $W'=(u',v')$;
moreover, for a point in $\mathcal{C}$,
the left and right states coincide (see Def.~\ref{carateristica}).
Because $\mathcal{U}_s$ and $\mathcal{U}_f$
both belong to $\mathcal{\mathcal{U}}$,
they have the same state $W$.
Likewise, $\mathcal{U}'_s$ and $\mathcal{U}'_f$
belong to $\mathcal{H}'(\mathcal{U})$,
so they have the same state $W'$.
In particular,
$\sigma(\mathcal{U}_s) = \lambda_s(W)$,
$\sigma(\mathcal{U}_f) = \lambda_f(W)$,
$\sigma(\mathcal{U}'_s) = \lambda_s(W')$,
and $\sigma(\mathcal{U}'_f) = \lambda_f(W')$.

\begin{figure}[htpb]
\begin{center}
\includegraphics[scale=0.43]{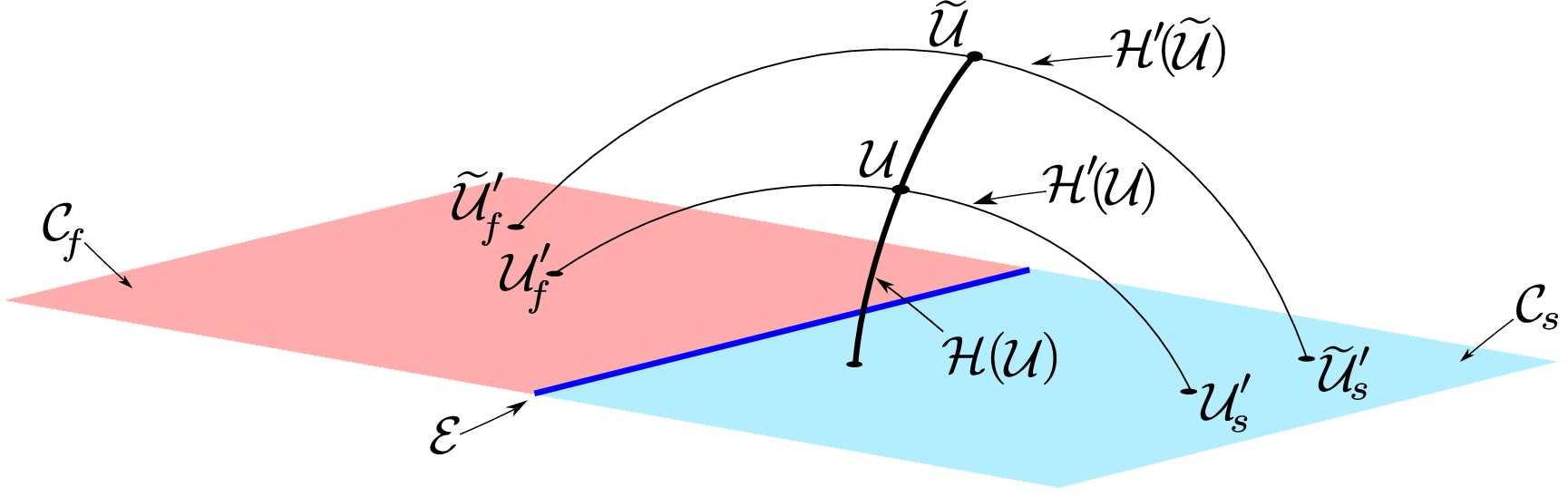}
\caption[]{
For a point $\mathcal{U}\in\mathcal{W}$,
the projections $\mathcal{U}'_{s}$ and $\mathcal{U}'_{f}$ of $\mathcal{U}$
onto $\mathcal{C}_{s}$ and $\mathcal{C}_{f}$, respectively,
are constructed using the Hugoniot$'$ curve $\mathcal{H}'(\mathcal{U})$.
In the same manner, for $\widetilde{\mathcal{U}} \in \mathcal{H}(\mathcal{U})$,
we project $\widetilde{\mathcal{U}}$
onto $\mathcal{C}_{s}$ and $\mathcal{C}_{f}$,
obtaining $\widetilde{\mathcal{U}}'_{s}$ and $\widetilde{\mathcal{U}}'_{f}$.
}
\label{fig:shock}
\end{center}
\end{figure}

Figure~\ref{fig:shock} illustrates these definitions.
This drawing is justified because we use $(z, \tau, Y)$-coordinates,
so that the plane $\Pi$ containing
the secondary bifurcation lines $\mathcal{B}$ and $\mathcal{B}'$
is excluded.
In effect, we adopt another restriction:

\begin{restrict}
\label{restrict:secbif}
We treat only Hugoniot and Hugoniot$'$ curves that do not meet
the secondary bifurcation lines $\mathcal{B}$ and $\mathcal{B}'$,
\emph{i.e.}, do not self-intersect.
\end{restrict}
\par\noindent
For the model studied in this paper,
such a Hugoniot or Hugoniot$'$ curve is homeomorphic to $\mathbb{R}$.

Fix $\mathcal{U}\in\mathcal{W}$ and define the points
in \eqref{usuli}--\eqref{usuli'}. We define a slow Lax shock curve
to be a maximal of $\mathcal{H}(\mathcal{U})$
such that, for all $\mathcal{U}$ in its interior,
\begin{equation}
\text{$\sigma(\mathcal{U}) < \sigma({\mathcal{U}}_{s})$
and $\sigma(\mathcal{U}'_{s}) < \sigma(\mathcal{U})
< \sigma(\mathcal{U}'_{f})$.}
\label{lax1}
\end{equation}
Similarly, we define a fast Lax shock curve to be
a maximal arc of $\mathcal{H}'(\mathcal{U})$
such that, for all $\mathcal{U}$ in its interior,
\begin{equation}
\text{$\sigma(\mathcal{U}) > \sigma(\mathcal{U}'_{f})$
and $\sigma(\mathcal{U}_{s}) < \sigma(\mathcal{U})
< \sigma(\mathcal{U}_{f})$.}
\label{lax2}
\end{equation}
Because shock speeds equal characteristic speeds on $\mathcal{C}$,
inequalities~\eqref{lax1} and \eqref{lax2}
are equivalent to the shock inequalities of Lax~\cite{lax}.

Suppose that a slow Lax shock curve
begins or ends at the point $\mathcal{U}$.
One possibility is that $\mathcal{U}$ belongs to $\mathcal{C}_s$.
Otherwise, if $\mathcal{U}\notin\mathcal{C}_s$,
then $\sigma(\mathcal{U})$ equals either $\sigma(\mathcal{U}_s)$,
$\sigma(\mathcal{U}'_s)$, or $\sigma(\mathcal{U}'_f)$,
meaning, respectively,
that $\mathcal{U}$ belongs to $Son'_s$, $Son_s$, or $Son_f$.
Examining inequalities~(\ref{lax1}),
we deduce that,
if a slow shock curve is oriented in the direction of decreasing speed,
it can only begin at $\mathcal{C}_s$ or $Son'_s$
and can only end at $Son_s$ or $Son_f$.
A slow shock curve starting at $\mathcal{C}_s$ is called \emph{local},
whereas one starting at $Son'_s$ is called \emph{non-local}.

Analogous statements hold:
if a fast Lax shock curve is oriented in the direction of increasing speed,
it can only begin at $\mathcal{C}_f$ or $Son_f$
and can only end at $Son'_s$ or $Son'_f$;
a fast shock curve is local if it starts at $\mathcal{C}_f$
and non-local if it starts at $Son_f$.
\subsection{Incorporating Lax conditions into the decomposition of $\mathcal{W}$ \label{subsec:NEW4.3}}

Since we are considering only Hugoniot curves that cross $\mathcal C$, we must add the $SCC$ surface to Fig.~\ref{fig:02}. Doing so, Region $11$ is subdivided into 3 regions: the region under $\mathcal {C}_s$ which we call $11_s$; the region under $\mathcal {C}_f$, $11_f$; and the region inside the $SCC$ surface (magenta surface in Fig.~\ref{fig:02SCC}), which is
discarded because it is formed by Hugoniot curves that do not intersect $\mathcal C$, In the same way, Region $10$ is subdivided into 3 regions: inside of $SCC’$ (light pink surface in Fig.~\ref{fig:02SCC}); the region above $\mathcal{C}_s$, $10_s$, and the region above $\mathcal C_f$, $10_f$.

\begin{proposition}
All points in Regions $8$ and $11_s$ belong to slow Lax shock curves, and all points in Regions $6$ and $10_f$ belong to fast Lax shock curves. There are no Lax shock curves in any other region.    
\end{proposition}

\begin{proof} We will focus on slow Lax shock curves. The result for fast Lax shock curves follows by using the reflection map defined in Section \ref{sec:2.2}.
Let us denote by $Lax\,1$ the condition that the arc is oriented in the direction of decreasing $\sigma$, by $Lax\,2$ the condition $\sigma(\mathcal U) \leq \sigma(\mathcal {U}_s)$, and by $Lax\,3$ the condition $\sigma(\mathcal {U’}_s) \leq \sigma(\mathcal {U})\leq \sigma(\mathcal {U’}_f )$.

As mentioned in the previous section, Lax shock curves can only begin and end at $\mathcal C$, $Son$, and $Son’$. Therefore, it is sufficient to check Lax conditions 1, 2, and 3 at the boundaries of the regions.

Let us start with the regions having $\mathcal C_s$ in their boundaries. In this case, it is sufficient to check $Lax\,1$ (direction of decreasing speed), since $Lax\,2$ and $Lax\,3$ are satisfied in $\mathcal C_s$.


Given a point in\ \ $\mathcal{C}_{s}$,\ \ we know that the Hugoniot curve is transversal to\ \ $\mathcal{C}$\ \ at this point. Let us determine to which side of $\mathcal{C}_{_{s}}$ $\sigma$ decreases. We have $\sigma(z_{0}, \tau_{0},Y_{0},z)$ given by \eqref{eq:velhug}. We call it $\sigma_{hug}$. Differentiating $\sigma$ with respect to $z$ and setting $Y_{0}=0$ and $z=z_{0}$, we can see that $d\sigma/dz$  is positive between the branches of the inflection locus and negative out of them. It remains to check in which direction the Hugoniot curve crosses $\mathcal{C}$ as $z$ increases. This is done by computing $\frac{dY_{hug}}{dz}$ from \eqref{eq:hugponto} and again setting $Y_{0}=0$ and $z=z_{0}$. We get $\frac{dY_{hug}}{dz}=-2c\tau_{0}(1+z_{0}^2)/[1+(b_{1}-1)z_{0}^2]$ which has the opposite sign of $\tau_{0}$. So, in $\mathcal{C}_{_{s}}(\tau_{0}<0)$, $\sigma$ decreases as the Hugoniot curve passes from Region 7 (below $\mathcal{C}$) to Region 8 (above $\mathcal{C}$) and also from Region 10 to Region $11_{s}$, indicating that Regions 8 and $11_{s}$ contain local arcs  which satisfy the Lax conditions.

Passing to regions that have $Son’_s$ in their boundaries:\\
Comparing Figs.~\ref{fig:cor}, \ref{fig:02}, \ref{fig:02SCC}, and \ref{fig:son-linha-s-f},
we can see that $Son'_s$ separates Region 4 from Region 10 and Region 8 from Region 9. In each case, we must check all three Lax conditions. We will show that a point in $Son’$ satisfies $Lax\,3$ if and only if  $z^2 < 1/(b1+1)$. This will be used to simplify the verification of $Lax\,1$ and $Lax\,2$ conditions.

  Parametric equations for the Hugoniot$'$ curve through a point $(\tau_0,z_0,Y_0)$ are given by 

\begin{equation}
\left\{
\begin{array}{l}
Y= -\displaystyle \frac{A'z^2+B'z +C'}{(z_0^2+1)[(b_1-1)z^2+1]}\\
\tau =\displaystyle \frac{D'z^3+E'z^2+F'z+G'}{c(z_0^2+1)[(b_1-1)z^4+b_1z^2+1]},
\end{array}
\right.
\label{eq:hugpontoli}
\end{equation}

\noindent obtained from Equations in \eqref{eq:hugponto} by changing $Y$ into $-Y$ and $A'$, $B'$, $C'$, $D'$, $E'$, $F'$, $G'$ are obtained from the same equations by changing  $Y_0$ into $-Y_0$.
Parametric equations for Hugoniot$'$ curve through a point of $Son'$ are obtained from the Equations above by substituting $\tau_0$ into $\tau'_0$, where $\tau'_0$ is given in Equation \eqref{eq:velponto}, getting $(\tau_{hug'}(z), Y_{hug'}(z))$.



We need to calculate the speed $\sigma$ at the intersection points of the Hugoniot$'$ curve with $\mathcal{C}$ and then write the condition $Lax\,3$. As the inequalities involved change only in $Son$, it is enough to verify them at a point of $Son'$.


Let $z_{1}$ and $z_{2}$ , $z_{1}< z_{2}$, be the solutions of  the equation $Y_{hug'}=0$. The roots $z_{1}$ and $z_{2}$ are the $z$ coordinates of the intersection points of the Hugoniot$'$ curve with $\mathcal {C}$. The speed, $\sigma_{hug'}(z)$, along the Hugoniot$'$ curve through $(z_0,\tau'_0,Y_0)$ is obtained by changing $\tau$ into $\tau_{hug'}(z)$ in Equation \eqref{eq:speed}. We must compare the values $\sigma_{hug'}(z_1)$ and $\sigma_{hug'}(z_2)$ with $\sigma_{son'}$. Denoting the numerator of $\sigma_{hug'}(z)$, by $\sigma_{hug'}^{num}(z)$ and the denominator by $\sigma_{hug'}^{den}(z)$, and denoting  the numerator of $Y_{hug'}(z)$ by $Y_{hug'}^{num}(z)$, we can write

\begin{equation}
\sigma_{hug'}^{num}(z)=p_1(z)Y_{hug'}^{num}(z)+r_1(z)
\end{equation}
and
\begin{equation} 
\sigma_{hug'}^{den}(z)=p_2(z)Y_{hug'}^{num}(z)+r_2(z).
\end{equation}
It follows that, $\sigma_{hug'}^{\mathcal {C}}(z)$, the speed at points of $\mathcal C$, is given by
$$\sigma_{hug'}^{\mathcal {C}}(z)= \frac{r_1(z)}{r_2(z)},$$ which has the same value as $\sigma_{hug'}(z)$ at the intersection points of the Hugoniot$'$ curve through $(z_0,\tau'_0,Y_0)$ with $\mathcal {C}$. The expression of $\sigma_{hug'}^{\mathcal {C}}(z)$ is of the form $\displaystyle \frac{az+b}{cz+d}$. Our goal is to get a condition such that $\sigma_{son'}(z_0.\tau'_0,Y_0)$ satisfies $\sigma_{hug'}^{\mathcal {C}}(z_1)< \sigma_{son'}(z_0.\tau'_0,Y_0)<\sigma_{hug'}^{\mathcal {C}}(z_2)$. In a more general way, we have the following problem: given a number $\sigma$ we want a condition such that 
$$\frac{az_1+b}{cz_1+d}<\sigma<\frac{az_2+b}{cz_2+d},$$ 
where $z_1$ and $z_2$ are the roots of the polynomial $fz^2+gz+h$, equivalently the condition is 
$$(\frac{az_1+b}{cz_1+d}-\sigma)(\frac{az_2+b}{cz_2+d}-\sigma)<0.$$ 
Taking into account that $z_1+z_2=-g/f$ and $z_1z_2=h/f$, straightforward computations give that the condition is 
$$\frac{(c^2h-dcg+d^2f)\sigma^2+(agd+bcg-2(ach+bdf))\sigma+a^2h-abg+b^2f}{c^2h-dcg+d^2f}<0.$$

Changing $a$, $b$ to the coefficients of the numerator of $\sigma_{hug'}^{\mathcal {C}}(z)$, $c$ and $d$ to the coefficients of the denominator of $\sigma_{hug'}^{\mathcal {C}}(z)$ and $f$, $g$, $h$ to the coefficients of $Y_{hug'}^{num}(z)$, we get that the condition is $$cond= \frac{Y_0^2((b_1+1)z_0^2-1)}{2}<0,$$
it follows that the condition is satisfied if and only if $\displaystyle \frac{-1}{\sqrt{b_1+1}}<z< \displaystyle \frac{1}{\sqrt{b_1+1}}.$

 For $Lax\,1$ and $Lax\,2$ we must determine in which direction is $\sigma$ decreasing along a Hugoniot curve through a point $\mathcal U_0=(\tau_0',z_0, Y_0)$ in $Son'_s$. Let $V_0$ be the tangent vector of such a Hugoniot curve at $\mathcal U_0$. Let $W_0$ be the gradient of $sonli$ at $\mathcal U_0$, where $sonli$ is the expression in Equation \eqref{eq:son'}. Since $sonli(0,0,0)<0$, we see that $W_0$ points away from $(0,0,0)$. Since $Son'$ divides $\mathcal W$ into two regions let us call $O$ the region which contains $(0,0,0)$ and $nonO$ the region which does not contain $(0,0,0)$. Let $sca$ be the scalar product of $V_0$ and $W_0$. We will show that for any $\mathcal U_0$ in $Son'_s$, $sca$ is positive, so the Hugoniot curve crosses $Son’_s$ entering $nonO$ as $z$ increases. 

 We will restrict ourselves to points in $Son'_s$ such that $z^2<1/(b_1+1)$, i.e., points in the region between the straight lines that form the double sonic locus, since this is the $Lax\,3$ condition, as we have seen.

 As before, using Eq. \eqref{eq:velponto}, writing that the point is in
$Son’$ (changing $\tau_0$ into $\tau_0'$, with $\tau'_0$ as in \eqref{eq:velponto}) and differentiating with respect to $z$, we get $$d\sigma_{Son'}= \frac{Y_0(-1+(b1+1)z_0^2)}{(1+(b1-1)z_0^2)},$$ which has the sign of $-Y_0$, since $z_0^2<1/(b1+1)$. So, if $Y_0>0$, $\sigma$ decreases as $z$ increases, i.e., as the Hugoniot curve enters $nonO$, and if $Y_0<0$ $\sigma$ decreases as the Hugoniot curve enters $O$.

Let us compute $sca$. Solving \eqref{eq:son'} for $\tau$ and replacing it in Eq. \eqref{eq:hugponto}, we get the expression of the Hugoniot curve from a point in $Son'_s$ parametrized by $z$ and $Y$. Differentiating with respect to $z$, we get $V_0$. In the same way, we compute the gradient of $sonli$ and replace $\tau$ by its value from Eq. \eqref{eq:son'} to get $W_0$.

Computing the scalar product, we find that $sca$ equals a fraction with a numerator
$$-4(z^2+1)(z^2(b1+1)-1)(A(z)Y+B(z))$$ and denominator $$z(1+(b1-1)z^2)(3+(b1+1)z^2).$$
Since we are restricting ourselves to $z^2<1/(b1+1)$, the denominator has the sign of $-z$.
So, the sign of $sca$ is the same as the sign of $z(A(z)Y+B(z))$, where  $A(z)= A_{\tau_{\mathcal C 2}} /2$ and $B(z)=B_{\tau_{\mathcal C 2}}$ from \eqref{eq:A.16} and \eqref{eq:A.17}.
It is shown in Appendix \ref{sec:propson} that a point $(z,\tau,Y)$ is in $Son’_s$ if and only if $z(A(z)Y+B(z))>0$, so $sca>0$.

Returning to Fig.~\ref{fig:02}, we see that the part of $Son’_s$ which separates Regions $4$ and $10$ does not satisfy condition $Lax\,3$, and in the part which separates Regions $8$ and $9$, all three conditions are satisfied by arcs that enter Region $8$ (contained in $nonO$), so in the Region $8$ we have both, local and non-local  slow Lax shock arcs.
It is easy to see that in the Region $8$, the saturation of the inflection locus by Hugoniot curves separates local from non-local slow Lax shock arcs.

Let us now consider $Y<0$. Fig.~\ref{fig:son-linha-s-f}, shows that $Son’_s$ is contained in the $z<0$ quadrant and separates Regions $6$ and $11_s$. We have seen that $\sigma$ is decreasing as the Hugoniot curve enters $O$, i.e., enters Region $11_s$. Therefore, Region $11_s$ contains non-local slow Lax shock arcs, which are separated from local slow Lax shock arcs by the saturation of the inflection locus by Hugoniot curves.

To apply the reflection map and obtain fast Lax shock curves, we must change $Y$ into $-Y$ and $z$ into $-z$ (remembering that $\sigma$ is increasing instead of decreasing), so Region $8$ becomes Region $6$ and Region $11_s$ becomes Region $10_f$.

\end{proof}




\begin{remark} Throughout this paper, we have used $z$, $\tau$, and $Y$ coordinates, which do not cover the plane $z=\infty$, formed by Hugoniot and Hugoniot’ curves. Since lax conditions do not change when we cross the $z=\infty$ plane it is interesting for this division in regions to notice which pairs of regions connect at infinity. For that, we will use coordinates $Z$. $T$ , and $X$ defined by  $Z=1/z$, $X=zY$, $T=(1+z^2).\tau$, which are valid near the $z=\infty$ plane.We can use these coordinates to obtain a figure analogous to Fig.~\ref{fig:02} and then check that we have the following pairs of regions connecting at  $z=\infty$ $1$ and $3$, $2$ and $4$, $5$ and 7, $6$ and $10$, $8$ and $11$.  Regions $9$ and $12$ do not extend to  $z=\infty$, so they do not connect.
\end{remark}



\subsection{Composite curves\label{subsec:Compos}}
Given a point $\mathcal{U}$ in $Son' \backslash ECC'$,
the Hugoniot curve of $\mathcal{U}$
intersects $\mathcal{C}$ at two distinct points 
$\mathcal{U}_{s}$ and $\mathcal{U}_{f}$
and either $\sigma(\mathcal{U}_{s})=\sigma(\mathcal{U})$
or $\sigma(\mathcal{U}_{f})=\sigma(\mathcal{U})$;
let $T(\mathcal{U})$ denote this point.
Then $T$ is a smooth map from $Son'$ to $\mathcal{C}$
satisfying $\sigma(T(\mathcal{U}))= \sigma(\mathcal{U})$.
The curves obtained by pulling
back rarefaction curves from $\mathcal{C}$ to $Son'$
are called \emph{composite curves}.
They are described in \cite{eschenazi02}.

Let us derive the differential equation for composite curves
in coordinates $z$ and $Y$.
By Prop.~\ref{propson},
a point $\mathcal{U}_{0}=(z_{0},\tau'(z_{0},Y_{0}), Y_{0})$ in $Son'$
has an associated point $(z_{\mathcal{C}},\tau_{\mathcal{C}},0)$
having the same speed as $\sigma(\mathcal{U}_0)$.
Explicitly (see \ref{sec:propson}),
\begin{align}
z_\mathcal{C}
&=-\frac{2\,(Y_{0}-c)\,z_{0}}{(b_{1}+1)\,Y_{0}\,z_{0}^{2}+Y_{0}+2c},
\label{eq:z_C} \\ 
\tau_\mathcal{C}&=\frac{
(A_\mathcal{C}\,Y_0 + B_\mathcal{C})[(\theta_+\,z_0^2+1)\,Y_0 + 2\,c]^2
}{2\,c\,z_0\,(z_0^2\,\theta_+ + 3)
\,(C_\mathcal{C}\,Y_{0}^2 + D_\mathcal{C}\,Y_{0} +E_\mathcal{C})},
\label{eq:tau_C}
\end{align}
with $\theta_+=b_1+1$, $\theta_-=b_1-1$,
$\vartheta_1=z_{0}^{2}+1$,
$A_\mathcal{C}=\theta_+\,z_{0}^{4} + 2\,(b_{1}+3)\,z_{0}^{2}+1$,
$B_\mathcal{C}=2\,c\,(\theta_-\,z_{0}^{2} + 1)$,
$C_\mathcal{C}=\theta_+\,z_{0}^{4} + 2\,(b_{1}+3)\,z_{0}^{2}+1$,
$D_\mathcal{C}=4\,c\,(\theta_-\,z_{0}^{2} + 1)$,
and $E_\mathcal{C}=4\,c^2\,\vartheta_0$.

We perform the pullback via the map
$(z_{0},Y_{0}) \mapsto (z_\mathcal{C},\tau_\mathcal{C})$
of the rarefaction differential equation \eqref{eq:edoraref},
obtaining
\begin{equation}
\frac{\mu_{2}(z)\,Y^2+\mu_1(z)\,Y + \mu_0(z)}{z\,(z^2\,\theta_+ + 3)}\,dz
+ [\nu_1(z)\,Y + \nu_0(z)]\,dY=0,
\label{eq:edocomp}
\end{equation}
where $\mu_{0}(z)=-12c^{2}[\theta_+z_{1}^{2}-1)]
[\theta_-z_{1}^{2}+1]$,
$\mu_{1}(z)=4c[\theta_+^{2}z_{1}^{6}
+\theta_+(5b_{1}-7)z^{4}+3\vartheta_1]$,
$\mu_{2}(z)=-(z^{^{8}}\theta_+^4-4(2b_{1}+3)\theta_+^{2}z^{^{6}}
-2\theta_+(13b_{1}+1)z^{4}+12z^{2}+3$,
$\nu_{0}(z)=2c [\theta_+z^{2}-1][\theta_-z^{2}+1]$,
$\nu_{1}(z)=-[\theta_+^{3}z^{^{6}}
+\theta_+(7b_{1}-1)z^{4}+(7b_{1}-1)z^{2}+1]$,
and $\theta_+=b_1+1$, $\theta_-=b_1-1$, $\vartheta_1=z^2+1$.

Equation~\eqref{eq:edocomp} is singular at the points
$(z=-1/\sqrt{\theta_+},Y=0)$ and $(z=1/\sqrt{\theta_+},Y=0)$.
These are the intersection points of the double sonic locus
and the characteristic surface.
Furthermore, the differential equation is not defined at $z=0$.
Two other singularities of saddle type at $z=0$ and $z=\infty$
are discussed in \cite{eschenazi02}.

For this paper, we define a \emph{slow composite curve}
to be a composite arc in $Son'_{s}$ that starts at $\mathcal{I}_{s}$
and is followed in the direction of decreasing speed $\sigma$.
Similarly, we define a \emph{fast composite curve}
to be a composite arc in $Son_{f}$ that starts at $\mathcal{I}_{f}$
and is followed in the direction of increasing speed $\sigma$.
A composite wave curve is denoted by $\mathcal{C}o$.

\section{Wave Curves}
\label{sec:Wave_Curves}

Wave curves in $\mathcal{W}$ are constructed in this section.
A \emph{wave curve} in $\mathcal{W}$ is a sequence of
rarefaction, shock, and composite curves satisfying certain conditions.
As we do throughout the paper,
we require shock waves to satisfy Lax admissibility criterion.
See~\cite{schecter96} for a discussion of wave curves for diffusive systems
of conservation laws.
Our list of wave curves is not exhaustive:
we treat only wave curves that occur in the quadratic model being studied.
We follow the methodology outlined in Sec.~4 of \cite{AEMP10}.

\subsection{Notation for wave curves\label{sec:notation}}
To facilitate the description of wave curves,
we introduce the following notation.
\begin{itemize}
\item Slow rarefaction, composite, and shock curves
are denoted $\mathcal{R}_{s}$, $\mathcal{C}o_{s}$, and $\mathcal{S}_{s}$,
respectively.
To indicate the beginning and end points
$\mathcal{U}^{b}$ and $\mathcal{U}^{e}$,
we sometimes amend the notation as in
$\mathcal{R}_s(\mathcal{U}^{b}\incr\mathcal{U}^{e})$,
$\mathcal{S}_s(\mathcal{U}^{b}\decr\mathcal{U}^{e})$,
and $\mathcal{C}o_s(\mathcal{U}^{b}\decr\mathcal{U}^{e})$;
here the arrows indicate that the speed $\sigma$
increases ($\incr$) or decreases ($\decr$) when the wave curve is traversed.

\item Similarly, fast wave curves are indicated by the subscript $f$,
as in $\mathcal{R}_f(\mathcal{U}^{b}\decr\mathcal{U}^{e})$,
$\mathcal{S}_f(\mathcal{U}^{b}\incr\mathcal{U}^{e})$,
and $\mathcal{C}o_f(\mathcal{U}^{b}\incr\mathcal{U}^{e})$.
\end{itemize}

\subsection{Slow and fast wave curves\label{sec:swc}}
We limit ourselves to certain examples of Riemann solutions
for the model treated in this paper.
Here, we describe the wave curves that occur in these examples.

A slow wave curve beginning from the point
$\mathcal{U}_0\in \mathcal{C}_{s}$
is a succession of slow rarefaction,
shock curves, and composite curves.
There are rules for how each type of curve ends
and what type of curve follows it.
For the model under consideration,
the steps for constructing wave curves are as follows.

\begin{itemize}
\item \textbf{Step 1}:
A slow wave curve has an initial point $\mathcal{U}_0 \in \mathcal{C}_{s}$.
The speed at $\mathcal{U}_0$
(\emph{i.e.}, the slow characteristic speed at its corresponding state)
is denoted $\sigma_0$.

\item \textbf{Step 2}:
The first arc is a slow shock curve
that starts at $\mathcal{U}_0$ and extends in the direction of
decreasing speed $\sigma$.
There are two possibilities:
either (a)~the arc extends to infinity as $z \to \infty$
or (b)~it ends at $Son_f$,
where the speed attains a minimum.
We denote this endpoint as $\mathcal{U}_1$,
which can be located at $z = \infty$ (\emph{i.e.}, $Z = 0$).
When diagramming a wave curve,
our convention is to place this shock curve
on the left side of $\mathcal{U}_0$,
as in $\mathcal{S}_s(\mathcal{U}_1\incr\mathcal{U}_0)\,\mathcal{U}_0$.

\item \textbf{Step 3}:
The second arc is a slow rarefaction curve contained in $\mathcal{C}_{s}$,
which is traversed in the direction of increasing speed.
(See Fig.~\ref{fig:02raref}.)
This rarefaction curve begins at $\mathcal{U}_0$
and ends at a point $\mathcal{U}_2$;
it is denoted by $\mathcal{R}_s(\mathcal{U}_0\incr\mathcal{U}_2)$,
and the speed at $\mathcal{U}_2$ is denoted $\sigma_2$.
There are two possibilities for $\mathcal{U}_2$:
(a)~it belongs to the slow inflection locus $\mathcal{I}_{s}$,
where the speed attains a maximum or
(b)~it belongs to the coincidence locus $\mathcal{E}$.
Within a wave curve diagram, we place this rarefaction curve
on the right side of the initial point $\mathcal{U}_0$,
as in $\mathcal{U}_0\,\mathcal{R}_s(\mathcal{U}_0\incr\mathcal{U}_2)$.

\item \textbf{Step 4}:
If the endpoint $\mathcal{U}_2$ of the rarefaction curve
lies on the coincidence locus $\mathcal{E}$,
as in~2\,(b), the wave curve stops.
Otherwise, $\mathcal{U}_2$ lies on
the slow inflection locus $\mathcal{I}_{s}$,
and we add a third arc,
namely the slow composite curve starting at $\mathcal{U}_2$
and stopping at a point $\mathcal{U}_3$.
Because the composite arc corresponding to
$\mathcal{R}_s(\mathcal{U}_0\incr\mathcal{U}_2)$
can cross the double sonic line,
there are two possibilities:
either (a)~$\mathcal{U}_3 \in \mathcal{D} = Son_s \cap Son'_f$,
so that $\sigma(\mathcal{U}_3) = \sigma\left((\mathcal{U}_3)_f\right)$, or
(b)~$\mathcal{U}_3 \in \mathcal{H}(\mathcal{U}_0)$,
\emph{i.e.}, the left state of $\mathcal{U}_3$
coincides with the state of $\mathcal{U}_0$.

\item \textbf{Step 5}:
If $\mathcal{U}_3$ belongs to $\mathcal{D}$, and therefore $Son'_f$,
then the wave curve stops.
(It cannot be continued with a slow rarefaction arc.)
Otherwise, we add a fourth arc, which is a non-local shock curve
starting at $\mathcal{U}_3$
and stopping at $\mathcal{U}_{4} \in Son_f$,
denoted $\mathcal{S}_{s}(\mathcal{U}_3\decr\mathcal{U}_{4})$.
\end{itemize}

\begin{figure}[htpb]
\begin{center}
\begin{subfigure}{0.48\linewidth}
\includegraphics[width=\linewidth]{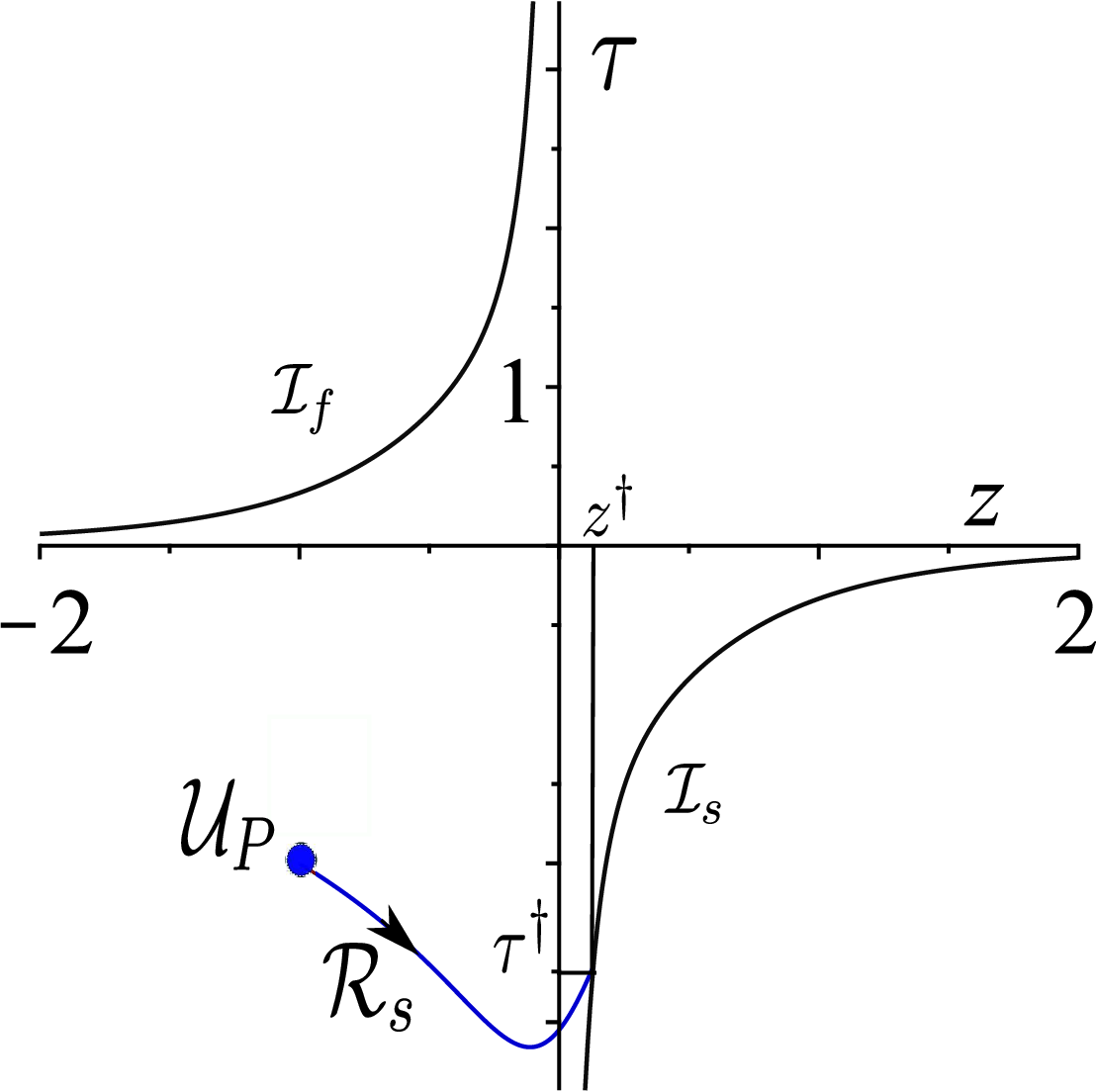}
\caption{}
\end{subfigure}
\hfil
\begin{subfigure}{0.48\linewidth}
\includegraphics[width=\linewidth]{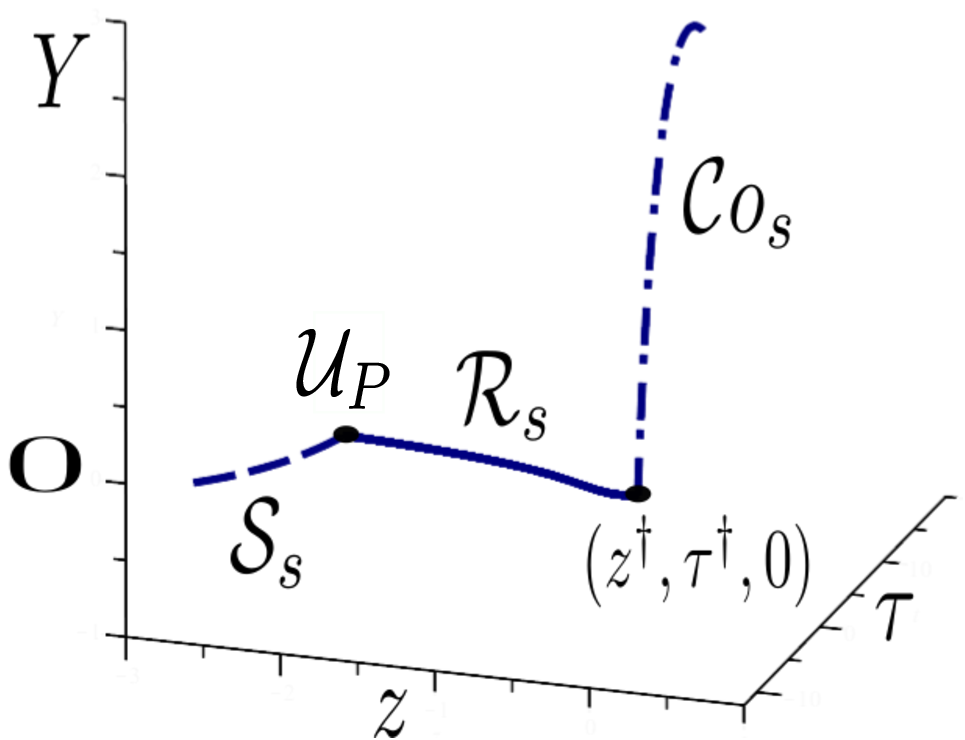}
\caption{}
\end{subfigure}
\caption[]{
The wave curve starting at the point
$\mathcal{U}_0 = \mathcal{U}_{P}$ for Case~2.1.
(a)~The rarefaction curve $\mathcal{R}_{s}$ in $\mathcal{C}_{s}$
reaches the slow inflection locus $\mathcal{I}_{s}$.
(See also Fig.~\ref{fig:subdivision} below.)
(b)~The wave curve is
$\mathcal{S}_{s}^\incr\mathcal{U}_{P}\mathcal{R}_{s}^\incr
\mathcal{C}o_{s}^\decr$ in $\mathcal{W}$.
}
\label{fig:wave_curve_example}
\end{center}
\end{figure}

We summarize the wave curve constructions as follows.
In all cases, the initial point of the wave curve
is denoted $\mathcal{U}_0\in\mathcal{C}_{s}$.
\begin{itemize}
\item \textbf{Case~1}:
$\mathcal{U}_2 \in \mathcal{E}$.
According to the constructions in steps~1--3,
the slow wave curve is
\begin{equation}
\text{$\mathcal{S}_s(\mathcal{U}_1\incr\mathcal{U}_0)
\,\mathcal{U}_0\,\mathcal{R}_s(\mathcal{U}_0\incr\mathcal{U}_2)$
or $\mathcal{S}_s^\incr\,\mathcal{U}_0\,\mathcal{R}_s^\incr$}.
\label{eq:structI}
\end{equation}
Either $\mathcal{U}_1 \in Son_f$ or it is located at $z = \infty$.

\item \textbf{Case~2.1}:
$\mathcal{U}_2 \in \mathcal{I}_{s}$
and $\mathcal{U}_3 \in \mathcal{D}$.
According to the constructions in steps~1--4,
the slow wave curve is
\begin{equation}
\text{$\mathcal{S}_s(\mathcal{U}_1\incr\mathcal{U}_0)
\,\mathcal{U}_0\,\mathcal{R}_s(\mathcal{U}_0\incr\mathcal{U}_2)
\,\mathcal{C}o_{s}(\mathcal{U}_2\decr\mathcal{U}_3)$
or $\mathcal{S}_s^\incr\,\mathcal{U}_0
\,\mathcal{R}_s^\incr\,\mathcal{C}o_s^\decr$.}
\label{eq:structII}
\end{equation}
Either $\mathcal{U}_1 \in Son_f$ or it is located at $z = \infty$.
This wave curve is illustrated
in Fig.~\ref{fig:wave_curve_example}.

\item \textbf{Case~2.2}:
$\mathcal{U}_2 \in \mathcal{I}_s$
but $\mathcal{U}_3 \notin \mathcal{D}$.
According to the constructions in steps~1--5,
the slow wave curve is
\begin{equation}
\text{$\mathcal{S}_s(\mathcal{U}_1\incr\mathcal{U}_0)
\,\mathcal{U}_0\,\mathcal{R}_s(\mathcal{U}_0\incr\mathcal{U}_2)
\,\mathcal{C}o_s(\mathcal{U}_2\decr\mathcal{U}_3)
\,\mathcal{S}_s(\mathcal{U}_3\decr\mathcal{U}_4)$
or $\mathcal{S}_s^\incr\,\mathcal{U}_0\,\mathcal{R}_s^\incr
\,\mathcal{C}o_s^\decr\mathcal{S}_s^\decr$.}
\label{eq:structIII}
\end{equation}
Here either $\mathcal{U}_1 \in Son_f$ or it is located at $z = \infty$,
and either $\mathcal{U}_4 \in Son_f$ or it is located at $z = \infty$;
moreover, it is not possible for both $\mathcal{U}_1$
and $\mathcal{U}_4$ to be located at $z = \infty$.
\end{itemize}

We refer to the sequence of rarefaction, shock, and composite curves
within a wave curve as its \emph{structure}.
The wave curve structures
that occur in the examples of Riemann solutions
in Sec.~\ref{sec:Riemann_Solutions}
are \eqref{eq:structI}, \eqref{eq:structII}, and \eqref{eq:structIII},
corresponding to Cases~1, 2.1, and 2.2.

A fast wave curve,
which starts at $\mathcal{U}_{0}\in \mathcal{C}_{f}$,
is analogous.
In fact, because our model is symmetric under the map taking
$x$ and $u$ to $-x$ and $-u$,
the possible fast wave curves are simply obtained from slow wave curves
by replacing $s$ by $f$ and interchanging the direction arrows
$\incr$ and $\decr$.

\subsection{Decomposing \texorpdfstring{$\mathcal{C}_{s}$}{C\_s}
according to wave curve structure}
\label{subsec:structure}
To determine the correct sequence of arcs
in a slow wave curve,
it is necessary to subdivide $\mathcal{C}_{s}$
into regions such that
wave curves starting in the same region
have similar properties.
Here, we subdivide $\mathcal{C}_{s}$ according to
properties of the rarefaction arc within
the slow wave curve.
(We do not take into account the properties
of the shock arcs.)

Along a rarefaction curve in $\mathcal{C}_{s}$,
the speed $\sigma$ increases towards $\mathcal{I}_{s}$,
as in Fig.~\ref{fig:02raref}.
With $(z,\tau,0)$ shortened to $(z,\tau)$,
the following statements hold.

\begin{figure}[htpb]
\begin{center}
\includegraphics[scale=1.1]{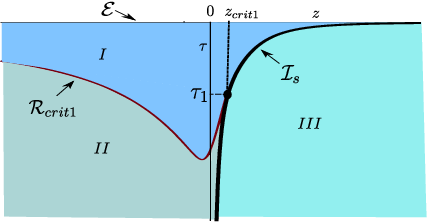}
\caption[]{
The subdivision of $\mathcal{C}_{s}$ into
the three regions~$I$, $II$, and $III$.
The point ($z_{crit1}$, $\tau_1$)
is defined in Eq.~\eqref{eq:zcrit1}. Here $z_{crit1}=1/3$ and $\tau_{1}=-6/5$.
}
\label{fig:subdivision}
\end{center}
\end{figure}

\begin{itemize}
\item[(a)]
If $(z,\tau)$ lies on the {left side} of $\mathcal{I}_{s}$
in Fig.~\ref{fig:subdivision} (Regions $I$ and $II$),
then the speed increases along the slow rarefaction curve through this point
as $z$ increases.

\item[(b)]
If a point $(z,\tau)$ lies on the {right side} of $\mathcal{I}_{s}$
in Fig.~\ref{fig:subdivision} (Region $III$),
then, the speed increases along the slow rarefaction curve through this point
as $z$ decreases.
\end{itemize}

The slow rarefaction curve of a point $(z,\tau)$ either
intersects $\mathcal{I}_{s}$
or ends at the coincidence curve $\mathcal{E}$.
The curve that separates these possibilities is
the slow rarefaction curve that ends at the point
$(z_{crit1},\tau_1)$ of Prop.~\ref{th:rarefp}.
This rarefaction curve,
denoted $\mathcal{R}_{crit1}$,
is parametrized by $z \le z_{crit1}$
and is tangent to $\mathcal{I}_s$ at $(z_{crit1},\tau_1)$.

Therefore, $\mathcal{C}_s$ is subdivided by
$\mathcal{I}_s$ and $\mathcal{R}_{crit1}$
into the three regions shown in Fig.~\ref{fig:subdivision}:
Region~$I$ (light blue) is bounded by
the coincidence curve $\mathcal{E}$,
by $\mathcal{R}_{crit1}$,
and by the $z \ge z_{crit1}$ portion of $\mathcal{I}_{s}$;
Region~$II$ (gray) is bounded by $\mathcal{R}_{crit1}$
and the $z \le z_{crit1}$ portion of $\mathcal{I}_{s}$;
and Region~$III$ (cyan) lies to the right of $\mathcal{I}_{s}$.

Let $\mathcal{R}_s$ denote the slow rarefaction curve
starting at $\mathcal{U}_0$.
If $\mathcal{U}_0$ belongs to Region~$I$,
then $\mathcal{R}_s$
stops at the coincidence curve $\mathcal{E}$;
in fact, the wave curve has the structure of Case~1.
If, instead, $\mathcal{U}_0$ belongs to Region~$II$ or~$III$,
then $\mathcal{R}_s$ intersects $\mathcal{I}_s$
and consequently is followed by a slow composite curve;
in fact, the wave curve has the structure of either Case~2.1 or~2.2.

Because the model being considered is symmetric ($b_2 = 0$),
it is preserved under the transformation $(x, u) \mapsto (-x, -u)$.
Therefore, fast wave curves are obtained from slow ones
by replacing $u$ by $-u$.
Accordingly,
the subdivision of $\mathcal{C}_f$ into regions
corresponding to the properties of the fast rarefaction curve
is analogous to the subdivision of $\mathcal{C}_s$.

\section{Wave Groups and the Intermediate Surface}
\label{sec:Wave_Groups_and_the_Intermediate_Surface}

In this section,
we describe how to use the wave curves of Sec.~\ref{sec:Wave_Curves}
to solve the Riemann problem~\eqref{eq:cons-law}
for the left state $W_L$ and right state $W_R$.

\subsection{Wave groups}
\label{subsec:wave_groups}
Each point $\mathcal{U}_M$ on the slow wave curve
from $\mathcal{U}_L \in \mathcal{C}_s$
determines a contiguous sequence of
slow rarefaction, shock, and composite waves
between $\mathcal{U}_L$ and $\mathcal{U}_M$,
called a \emph{slow wave group}.
The left-most state in this wave group is
the left state $W_L$ of $\mathcal{U}_L$,
and the right-most state in this wave group is
the right state $W_M$ of $\mathcal{U}_M$.

For example, suppose that the local portion of the wave curve is
\begin{equation}
\mathcal{S}_s(\mathcal{U}_1 \incr \mathcal{U}_L)
\,\mathcal{U}_L\,\mathcal{R}_s(\mathcal{U}_L \incr \mathcal{U}_1)
\,\mathcal{C}o_s(\mathcal{U}_1 \incr \mathcal{U}_3).
\label{eq:local}
\end{equation}
If $\mathcal{U}_M$ lies on the rarefaction curve
$\mathcal{R}_s(\mathcal{U}_L \incr \mathcal{U}_1)$,
then the corresponding wave group consists of
a rarefaction wave extending from $W_L$ to $W_M$.
The solution of $W(x, t)$ of Eq.~\eqref{eq:cons-law} is given by
\begin{equation}
W(x, t) = \begin{cases}
W_L & \text{if $x/t \le \sigma(\mathcal{U}_L)$,} \\
W^* & \text{if $\sigma(\mathcal{U}_L) \le x/t \le \sigma(\mathcal{U}_M)$,} \\
W_M & \text{if $\sigma(\mathcal{U}_M) \le x/t$,}
\end{cases}
\end{equation}
where $W^*$ is the state corresponding to
the point $\mathcal{U}^*$
along $\mathcal{R}_s(\mathcal{U}_L \incr \mathcal{U}_M)$
such that $\sigma(\mathcal{U}^*) = x/t$
(which exists and is unique because
$\sigma$ varies monotonically along a rarefaction curve).

If, on the other hand, $\mathcal{U}_M$ lies on the shock curve
$\mathcal{S}_s(\mathcal{U}_2 \incr \mathcal{U}_L)$,
then the corresponding wave group consists of a shock wave
with left state $W_L$ and right state $W_M$.
The solution of $W(x, t)$ of Eq.~\eqref{eq:cons-law} is given by
\begin{equation}
W(x, t) = \begin{cases}
W_L & \text{if $x/t <\sigma(\mathcal{U}_M)$,} \\
W_M & \text{if $\sigma(\mathcal{U}_M) < x/t$.}
\end{cases}
\end{equation}

Finally, if $\mathcal{U}_M$ lies on the composite curve
$\mathcal{C}o_s(\mathcal{U}_1 \incr \mathcal{U}_3)$,
and the left state of $\mathcal{U}_M$ is denoted $W_\ell$,
then the corresponding wave group consists of
a rarefaction wave extending from $W_L$ to $W_\ell$ adjoined, on its right side,
by the left-sonic shock wave from $W_\ell$ to $W_M$.
The solution of $W(x, t)$ of Eq.~\eqref{eq:cons-law} is given by
\begin{equation}
W(x, t) = \begin{cases}
W_L & \text{if $x/t \le \sigma(\mathcal{U}_L)$,} \\
W^* & \text{if $\sigma(\mathcal{U}_L) \le x/t < \sigma(\mathcal{U}_M)$,} \\
W_M & \text{if $\sigma(\mathcal{U}_M) < x/t$.}
\end{cases}
\end{equation}
Here $W^*$ is the state corresponding to
the point $\mathcal{U}^*$
along $\mathcal{R}_s(\mathcal{U}_L \incr \mathcal{U}_\ell)$
such that $\sigma(\mathcal{U}^*) = x/t$,
with $\mathcal{U}_\ell = (\mathcal{U}_M)_s \in \mathcal{C}_s$
corresponding to $W_\ell$.

\subsection{Intermediate surface}
\label{subsec:intermediate_surface}
In the classical approach to solving a Riemann problem,
one works with classical wave curves in state space,
which coincide with the projections of the wave curves in $\mathcal{W}$.
One seeks an intersection $W_M$ of
the slow wave curve for $W_L$
and the fast wave curve for $W_R$.
Then the solution comprises
the slow wave group between $W_L$ and $W_M$
and the fast wave group between $W_M$ and $W_R$.

This procedure can be carried out directly in the wave manifold~\cite{AEMP10}.
Let $\mathcal{U}_L$ denote the point in $\mathcal{C}_s$
corresponding to $W_L$,
and let $\mathcal{U}_R$ denote the point in $\mathcal{C}_f$
corresponding to $W_R$.
Then construct the slow wave curve for $\mathcal{U}_L$
and the fast wave curve for $\mathcal{U}_R$.
Unlike in the classical procedure,
1-dimensional wave curves in the 3-dimensional space $\mathcal{W}$
generally do not intersect.

For example, suppose that the slow and fast wave groups are both shock waves,
the slow shock wave having left state $W_L$ and right state $W_M$
and the fast shock wave having left state $W_M$ and right state $W_R$.
Represented within the wave manifold,
the slow shock wave $(W_L, \sigma_L, W_M)$
is a point $\mathcal{U}_{M,s} \in \mathcal{W}$,
whereas the fast shock wave $(W_M, \sigma_R, W_R)$
is $\mathcal{U}_{M,f} \in \mathcal{W}$.
Here, $\mathcal{U}_{M,s}$ lies on the slow wave curve of $\mathcal{U}_L$
and $\mathcal{U}_{M,f}$ lies on the fast wave curve of $\mathcal{U}_R$.
In general, $\mathcal{U}_{M,s}$ is distinct from $\mathcal{U}_{M,f}$:
$(W_L, W_M) \ne (W_M, W_R)$.
However, the right state of $\mathcal{U}_{M,s}$, namely $W_M$,
does coincide with left state of $\mathcal{U}_{M,f}$.

In general,
the requirement for obtaining a solution of the Riemann problem
is that some point $\mathcal{U}_{M,s}$
on the slow wave curve of $\mathcal{U}_L$
has the same right state as
the left state of a point $\mathcal{U}_{M,f}$
on the fast wave curve of $\mathcal{U}_R$.
To formulate this requirement within the wave manifold,
we apply the reflection operation of Def.~\ref{def:reflection}
to the fast wave curve.
Then the right state of $\mathcal{U}_{M,s}$
equals the right state of the reflection of $\mathcal{U}_{M,f}$,
(In the example of the previous paragraph,
the reflection of the fast shock wave $(W_R, \sigma_R, W_M)$,
which has the same right state as
the slow shock wave $(W_L, \sigma_L, W_M)$.)
Equivalently,
$\mathcal{U}_{M,s}$ and the reflection of $\mathcal{U}_{M,f}$
lie on the same Hugoniot$'$ curve.

Accordingly, we construct the Riemann solution
for $\mathcal{U}_L \in \mathcal{C}_s$
and $\mathcal{U}_R \in \mathcal{C}_f$
as follows.
\begin{itemize}
\item Form the slow wave curve of $\mathcal{U}_L$.
\item Draw the Hugoniot$'$ curve
from each point of this slow wave curve,
thereby forming a surface
called the \emph{intermediate surface of $\mathcal{U}_L$}.
\item Form the fast wave curve of $\mathcal{U}_R$ and reflect it.
\item Intersect the reflected fast wave curve with the intermediate surface,
obtaining a point $\mathcal{U}_{M,f}$.
\item Find a point $\mathcal{U}_{M,s}$
within the slow wave curve for which the Hugoniot$'$
curve contains $\mathcal{U}_{M,f}$.
\end{itemize}
Then the slow wave group corresponding to $\mathcal{U}_{M,s}$
and the fast wave group corresponding to $\mathcal{U}_{M,f}$
constitute a solution of the Riemann problem.

\begin{remark}
It is conceivable that the reflected fast wave curve of $\mathcal{U}_R$
fails to intersect the intermediate surface of $\mathcal{U}_L$,
in which case there is no solution of this particular Riemann problem.
Likewise,
the reflected fast wave curve
might have multiple intersections with the intermediate surface;
then there are multiple solutions of the Riemann problem.
\end{remark}

\begin{remark}\label{rem:snugly}
The Bethe-Wendroff theorem~\cite{bet42,wen72,wen72a,menPlo89,furtado}
identifies conditions under which a wave curve in state space is a smooth curve.
For the present model, all state-space wave curves are smooth.
It follows (see \cite{AEMP10}) that the intermediate surface
is a 2-dimensional submanifold of $\mathcal{W}$,
\emph{i.e.},
it is a differentiable surface,
with all pieces joining smoothly.
\end{remark}

\begin{figure}[htpb]
\begin{center}
\begin{subfigure}{0.48\linewidth}
\includegraphics[width=\linewidth]{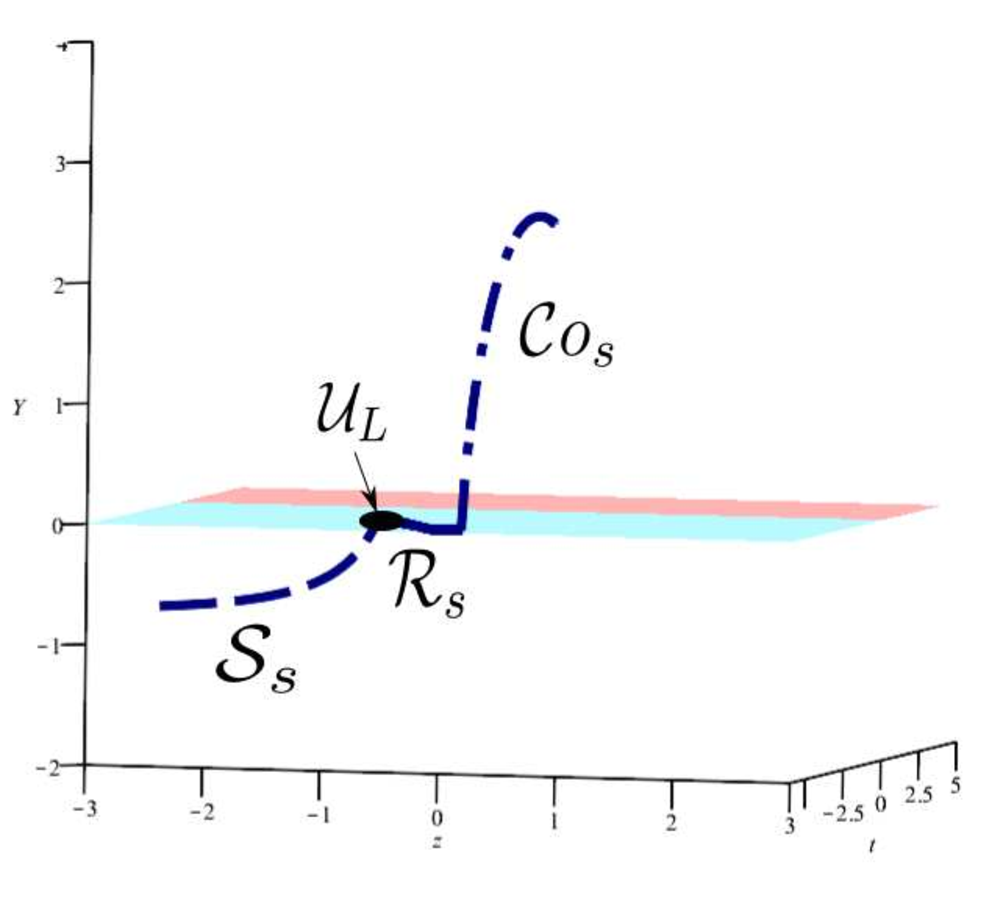}
\caption{}
\end{subfigure}
\hfil
\begin{subfigure}{0.48\linewidth}
\includegraphics[width=\linewidth]{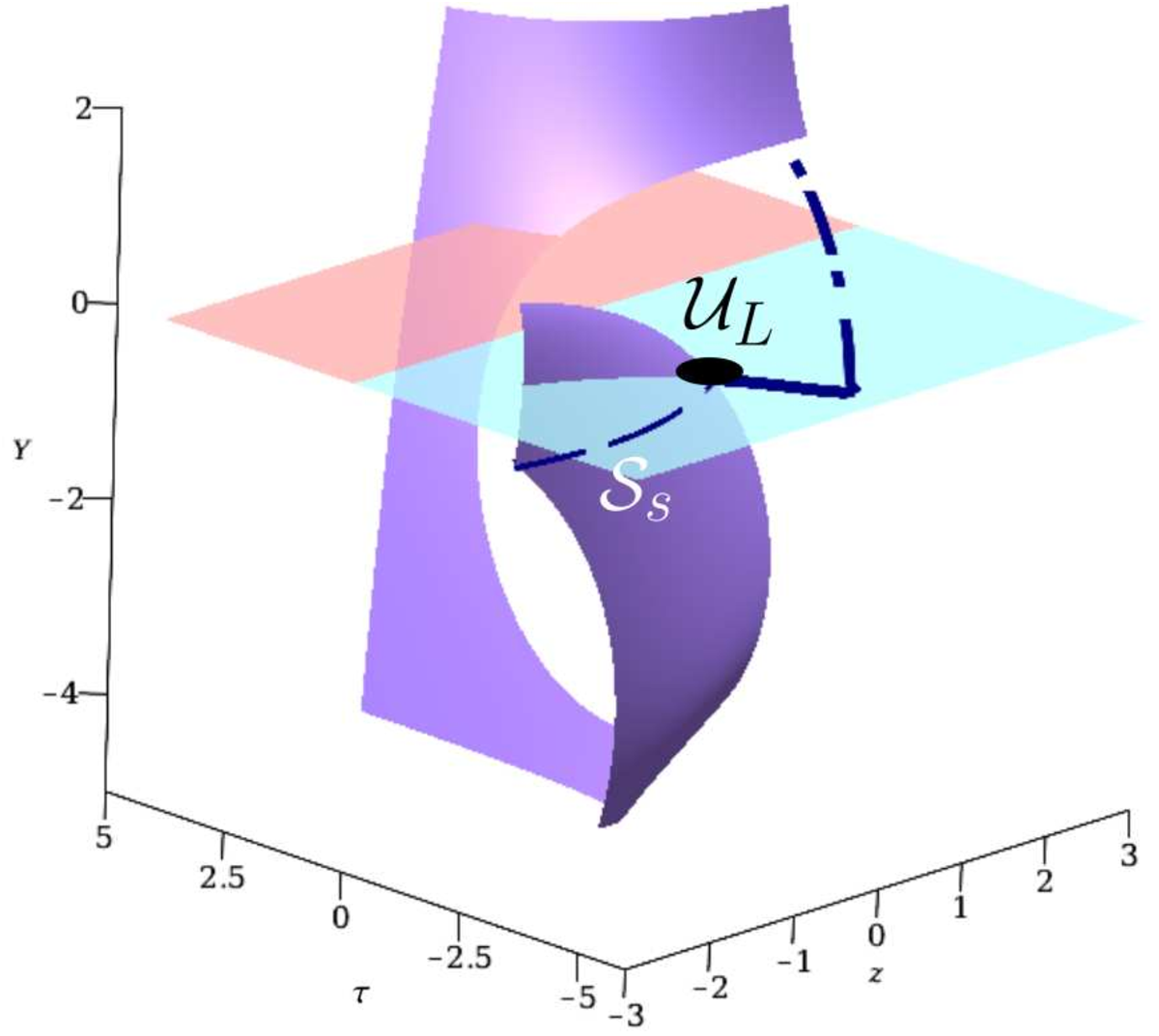}
\caption{}
\end{subfigure}
\caption[]{We recall from Fig.~\ref{fig:02} that
$\mathcal{C}_{s}$ is the light blue half plane
and $\mathcal{C}_{f}$ is the light coral half plane.
(a)~The dark blue slow wave curve from
a point $\mathcal{U}_{L}\in\mathcal{C}_{s}$.
The dashed curve is a slow shock curve ($\mathcal{S}_{s}$),
the solid curve is a rarefaction curve ($\mathcal{R}_{s}$),
and the dash-dotted curve is a composite curve ($\mathcal{C}o_{s}$).
The wave curve is
$\mathcal{S}_{s}^\incr\mathcal{U}_{L}\mathcal{R}_{s}^\incr\mathcal{C}o_s^\decr$.
(b)~The light purple intermediate surface of
the slow shock curve $\mathcal{S}_{s}$.
The wave curve is the same as in~(a),
viewed from a different perspective.
}
\label{fig:lenta1}
\end{center}
\end{figure}

\begin{figure}[htpb]
\begin{center}
\begin{subfigure}{0.48\linewidth}
\includegraphics[width=\linewidth]{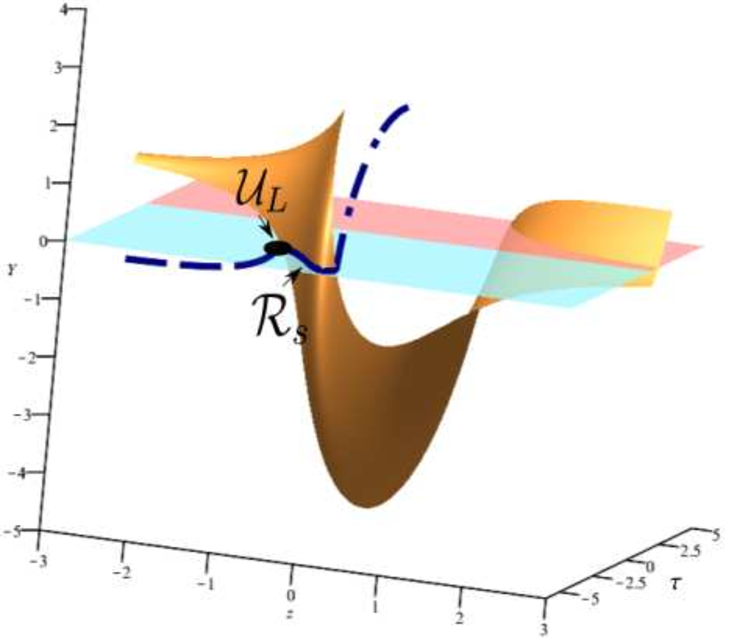}
\caption{}
\end{subfigure}
\hfil
\begin{subfigure}{0.48\linewidth}
\includegraphics[width=\linewidth]{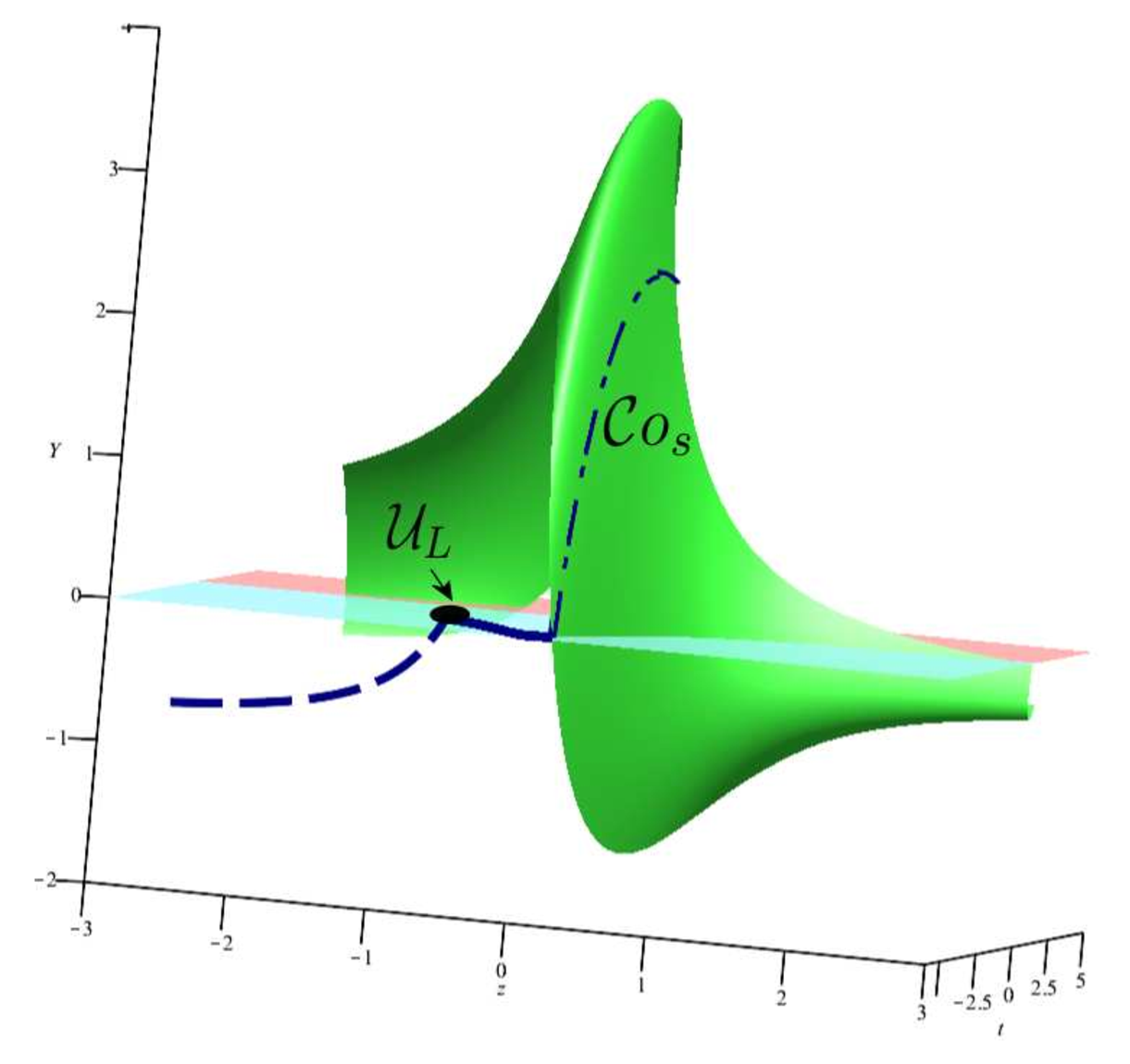}
\caption{}
\end{subfigure}
\caption[]{
(a)~The golden intermediate surface
of the slow rarefaction curve $\mathcal{R}_{s}$;
see the caption of Fig.~\ref{fig:lenta1}.
(b)~The green intermediate surface
of the slow composite curve $\mathcal{C}o_{s}$.
}
\label{fig:saturated1}
\end{center}
\end{figure}

\begin{figure}[htpb]
\begin{center}
\begin{subfigure}{0.48\linewidth}
\includegraphics[width=\linewidth]{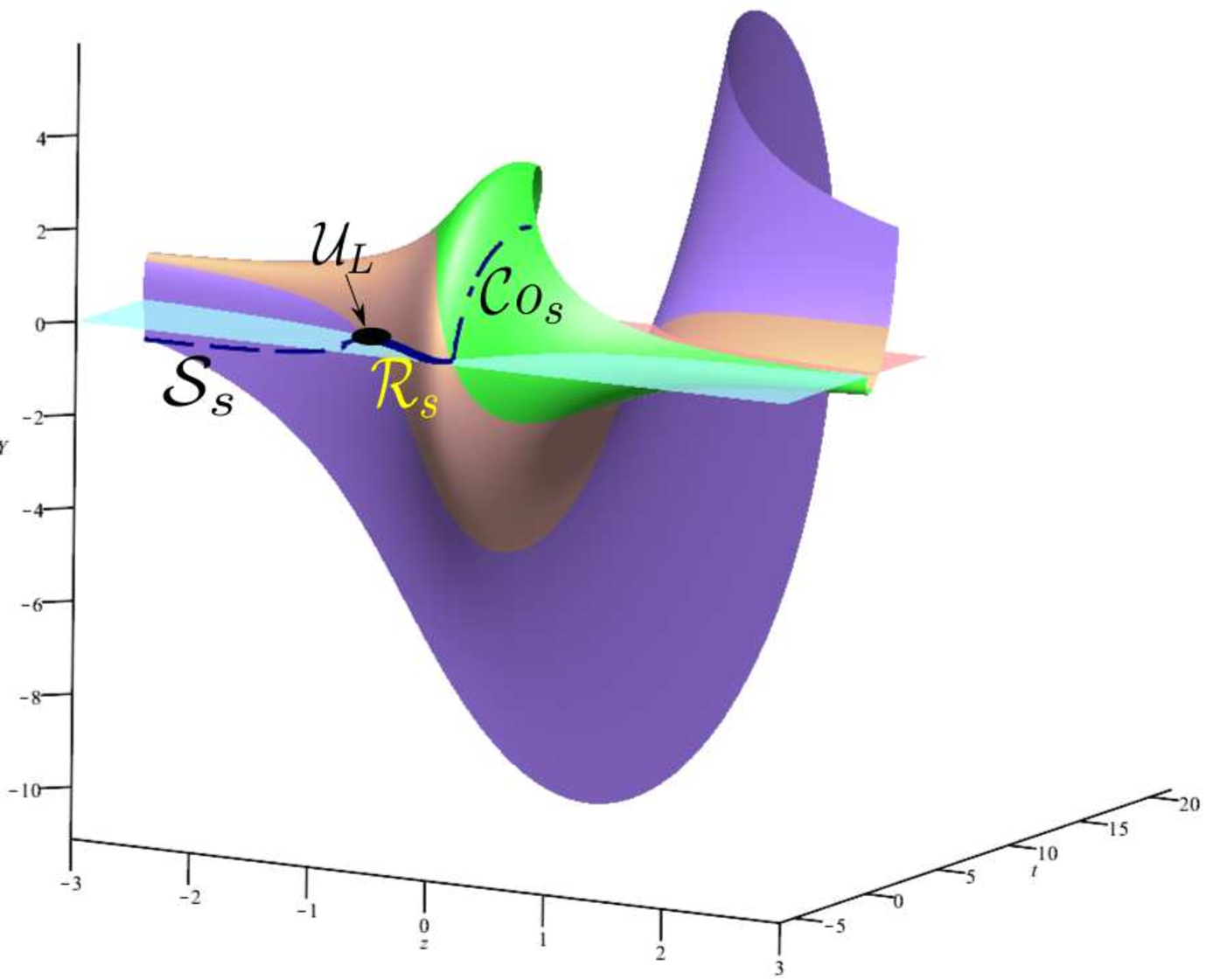}
\caption{}
\end{subfigure}
\hfil
\begin{subfigure}{0.48\linewidth}
\includegraphics[width=\linewidth]{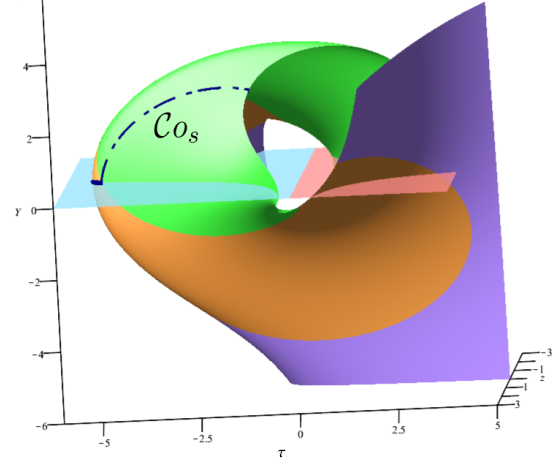}
\caption{}
\end{subfigure}
\caption[]{
Complete intermediate surface for $\mathcal{U}_{L}$,
viewed from two perspectives.
Notice that the different parts of the intermediate surface
generated by the different arcs
($\mathcal{S}_{s}$, $\mathcal{R}_{s}$, and $\mathcal{C}o_{s}$)
in the wave curve fit together smoothly,
as explained in Remark~\ref{rem:snugly}.
}
\label{fig:saturated3}
\end{center}
\end{figure}

In Figs.~\ref{fig:lenta1}-\ref{fig:saturated3},
we give an example of a slow wave curve and its intermediate surface.
We first draw the portions of the intermediate surfaces
associated with individual arcs in the wave curve;
then we draw the complete intermediate surface.
The portion for the shock curve is shown in Fig.~\ref{fig:lenta1}b;
the portion for the rarefaction curve is shown in Fig.~\ref{fig:saturated1}a;
and the portion for the composite curve is shown in Fig.~\ref{fig:saturated1}b.
The complete intermediate surface is shown in Fig.~\ref{fig:saturated3}.

\section{Riemann Solutions}
\label{sec:Riemann_Solutions}

This section describes a procedure for constructing solutions
of Eq.~\eqref{eq:cons-law}
that consist of rarefaction, shock, and composite waves
separated by constant states,
assuming that both $W_{L}$ and $W_{R}$ lie in
the strictly hyperbolic region in state space.
This procedure makes essential use of
wave curves and the intermediate surface within $\mathcal{W}$.
We examine several cases,
but do not establish existence or uniqueness of solutions of Riemann problems
for the particular model being studied.
Rather,
our purpose is to illustrate how visualizing the solution procedure
within the wave manifold $\mathcal{W}$ can aid in proving such results.

\subsection{Algorithm to find the Riemann solutions}
\label{algorithm}
We first \emph{raise} $W_L$ and $W_R$ from state space to $\mathcal{W}$,
obtaining $\mathcal{U}_{L}\in\mathcal{C}_{s}$
and $\mathcal{U}_{R}\in\mathcal{C}_{f}$, as follows.
First focus on $W_L$.
As $\mathcal{U}_L$ belongs to $\mathcal{C}$,
it has coordinates $Y=0$ and $X=z\,Y=0$.
By Eq.~\eqref{transf1},
$U=u_L$ and $V=v_L$,
so that by Eq.~\eqref{transf2},
$\widetilde{U}=b_1\,u_L$ and $V_1 = \widetilde{V} + c = v_L + c$.
The equation $G = 0$, with $G$ given by Eq.~\eqref{G},
is a quadratic equation for $z$.
Because $W_L$ lies in the hyperbolic region,
this equation has two distinct roots.
The value for $\tau$ corresponding to each root $z$
is determined using Eq.~\eqref{eq:tdef}.
By Prop.~\ref{th:slow_fast}, the two values of $\tau$ have opposite signs.
We choose the root $z_L$ so the corresponding $\tau_L$ is negative.
Thus, we obtain the coordinates $Y = 0$, $\tau = \tau_L$, and $z = z_L$
of $\mathcal{U}_L$.
Similarly, we find the coordinates $Y = 0$, $\tau = \tau_R$, and $z = z_R$
of $\mathcal{U}_R$, where $\tau_R$ is positive.

After obtaining $\mathcal{U}_{L}$ and $\mathcal{U}_{R}$,
we take the following steps to obtain the RPS.

\begin{enumerate}
\item Construct the slow wave curve starting from $\mathcal{U}_{L}$,
as in Sec.~\ref{sec:swc}.

\item Construct the intermediate surface,
as in Sec.~\ref{sec:Wave_Groups_and_the_Intermediate_Surface}.

\item Construct the fast wave curve starting from $\mathcal{U}_{R}$.

\item Form the reflection of the fast wave curve.

\item Find an intersection point $\mathcal{U}_{M,f}$
between the intermediate surface and the reflected fast wave curve.

\item Find a point $\mathcal{U}_{M,s}$
within the slow wave curve for which the Hugoniot$'$
curve contains $\mathcal{U}_{M,f}$.

\item From the slow wave group corresponding to $\mathcal{U}_{M,s}$
and the fast wave group corresponding to $\mathcal{U}_{M,f}$
construct the scale-invariant solution of Eq.~\eqref{eq:cons-law}.
\end{enumerate}

This algorithm was proposed in \cite{AEMP10}.
However, that paper omits the step in which the fast wave curve is reflected.
(In the second and third paragraphs of Sec.~6, p.~385, of Ref~\cite{AEMP10}
there are some mistakes.)
We correct this omission in the present paper.

\subsection{Examples of solutions}\label{sec:examples}
\begin{figure}[htpb]
\begin{center}
\begin{subfigure}{0.45\linewidth}
\includegraphics[width=\linewidth]{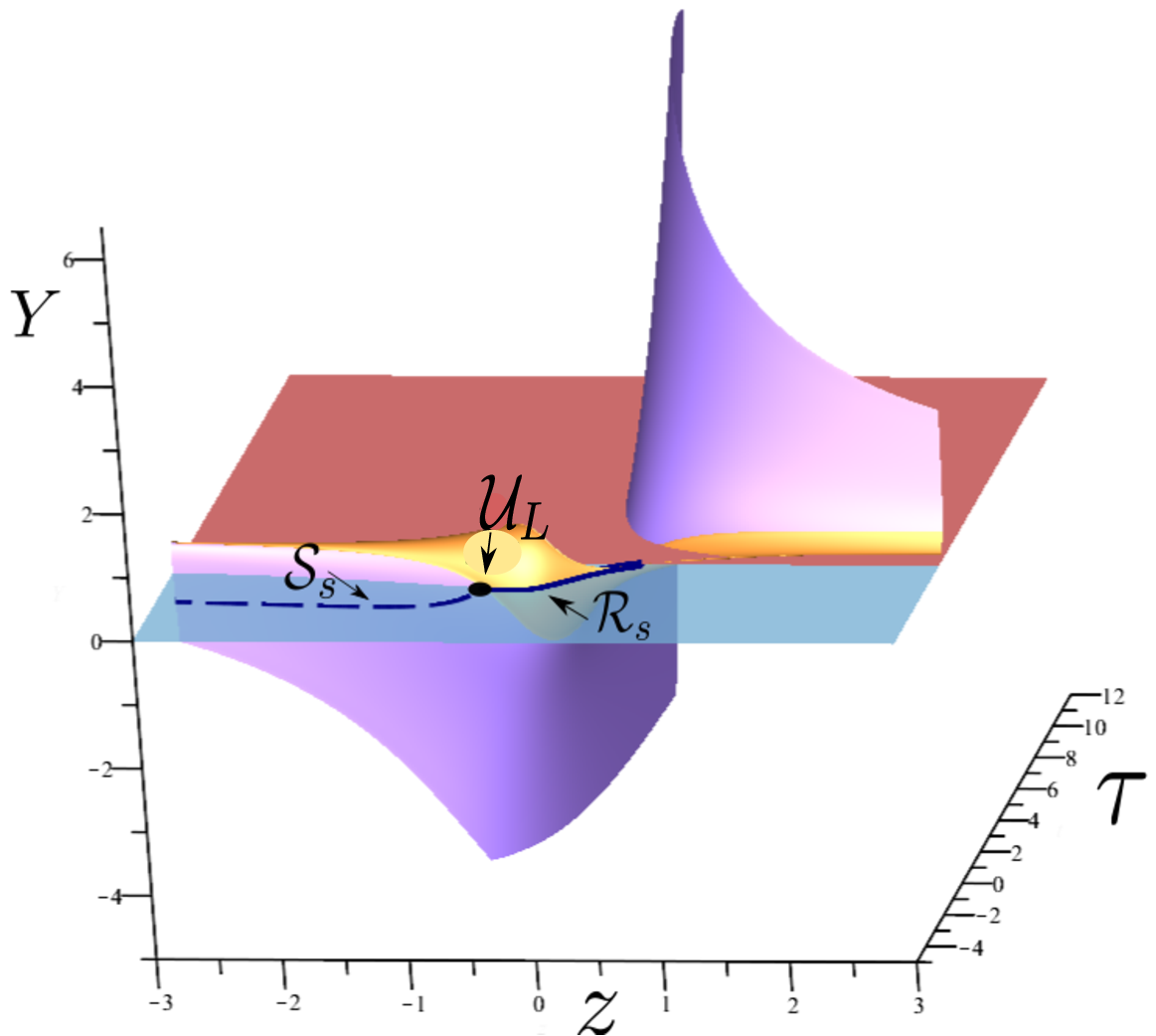}
\caption{}
\end{subfigure}
\hfil
\begin{subfigure}{0.45\linewidth}
\includegraphics[width=\linewidth]{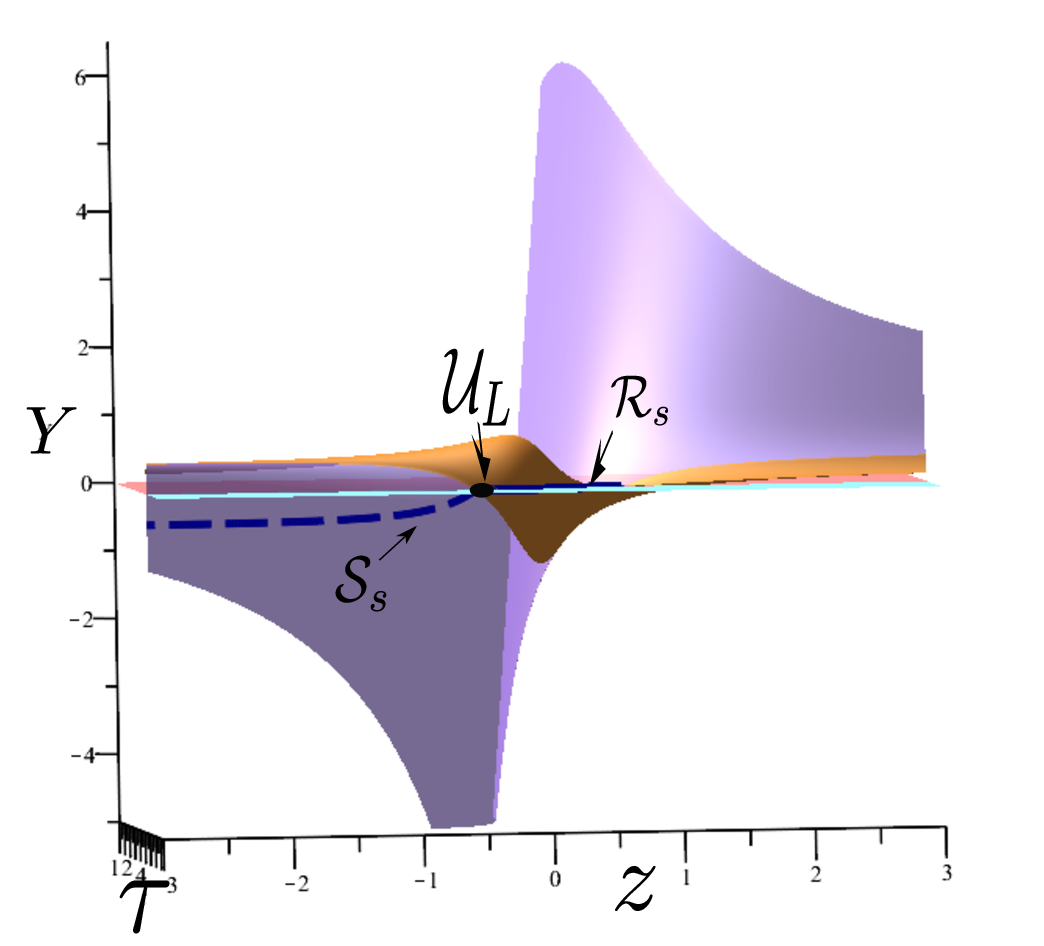}
\caption{}
\end{subfigure}
\caption[]{
The slow wave curve and its intermediate surface
for $\mathcal{U}_{L}$ in Region~$I$,
viewed from two perspectives.
Here $z_{L}=-0.5$ and $\tau_{L}=-1.5$.
The slow wave curve (dark blue)
is $\mathcal{S}_{s}^\incr\mathcal{U}_{L}\mathcal{R}_{s}^\incr$.
The intermediate surfaces for $\mathcal{S}_{s}$ and $\mathcal{R}_{s}$
are colored light purple and gold, respectively.
}
\label{fig:sol3}
\end{center}
\end{figure}

\medskip\noindent
\textbf{Example 1.}
First we choose $W_{L}=(u_{L},v_{L})$ such that
the corresponding point $\mathcal{U}_{L}=(z_L, \tau_L)$
lies in Region~$I$ of the subdivision of $\mathcal{C}_s$.
In this example, we take $u_{L}=-0.04$ and $v_{L}=-0.55$;
then $z_{L}=-0.5$ and $\tau_{L}=-1.5$.
(See Sec.~\ref{algorithm}.)
The slow wave curve of $\mathcal{U}_{L}$ is
$\mathcal{S}_{s}^\incr\mathcal{U}_{L}\mathcal{R}_{s}^\incr$
because the rarefaction curve ends at the coincidence curve, $\tau=0$.
Fig.~\ref{fig:sol3} depicts the slow wave curve, colored dark blue;
the dashed curve is the shock curve $\mathcal{S}_s$
and the solid curve is the rarefaction curve $\mathcal{R}_s$.
This figure also shows the intermediate surface for the slow wave curve,
formed out of Hugoniot$'$ curves through points on $\mathcal{S}_s$
(the light purple surface) and $\mathcal{R}_s$ (the gold surface).

\begin{figure}[htpb]
\begin{center}
\begin{subfigure}{0.40\linewidth}
\includegraphics[width=\linewidth]{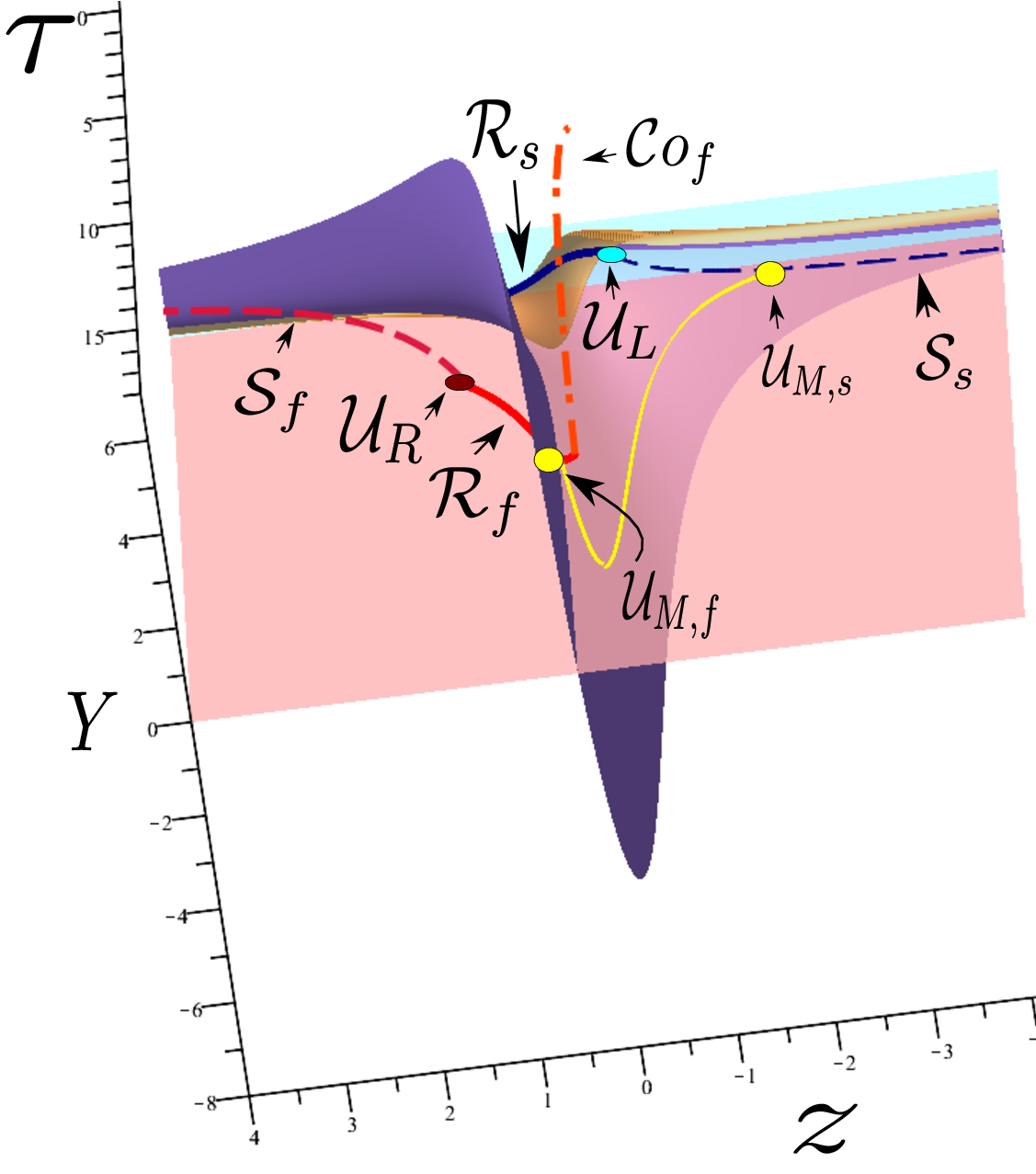}
\caption{}
\end{subfigure}
\hfil
\begin{subfigure}{0.40\linewidth}
\includegraphics[width=\linewidth]{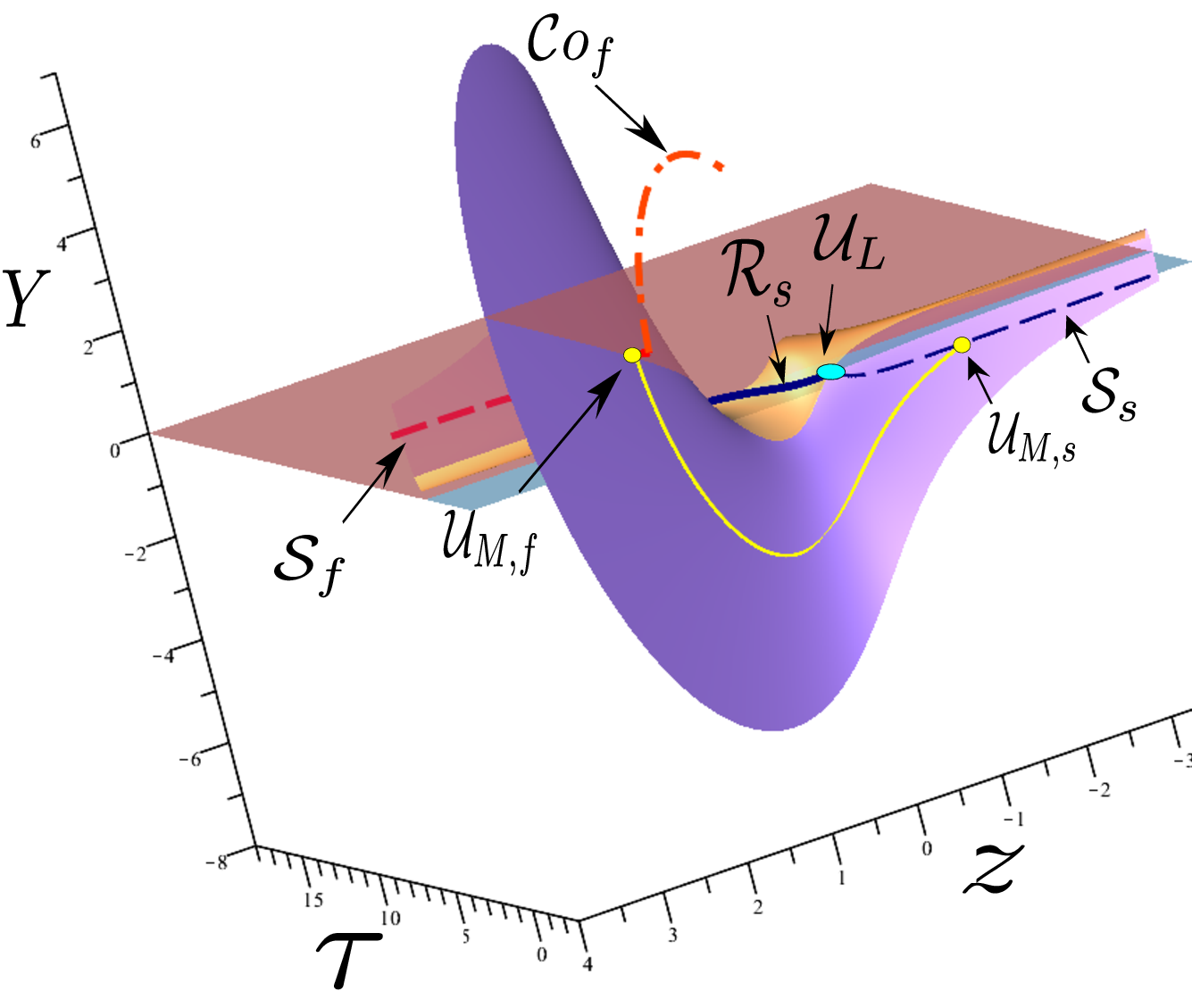}
\caption{}
\end{subfigure}
\caption[]{
Example~1 of a Riemann problem solved using wave curves in $\mathcal{W}$,
shown from two vantage points.
The reflected fast wave curve from $\mathcal{U}_{R}$
is $\mathcal{S}_{f}^\decr\mathcal{U}_{R}\mathcal{R}_{f}^\decr
\mathcal{C}o_{f}^\incr$.
Its reflection is depicted as the red curve.
In this example, the reflection of $\mathcal{R}_f$,
which coincides with itself,
intersects the intermediate surface of $\mathcal{S}_{s}$
(colored light purple) at $\mathcal{U}_{M,f}$.
A Hugoniot$'$ curve (yellow) connects $\mathcal{U}_{M,f}$
to a point $\mathcal{U}_{M,s}$ within $\mathcal{S}_s$.
The RPS consists of the slow shock wave $\mathcal{U}_{M,s}$
and the fast rarefaction wave extending from
$\mathcal{U}_R$ to $\mathcal{U}_{M,f}$.
}
\label{fig:sol2}
\end{center}
\end{figure}

As an illustration of the procedure for solving Riemann problems,
we take $W_R=(u_R, v_R)=(0.125, 3.5)$,
for which the corresponding point in $\mathcal{C}$ is
$\mathcal{U}_{R}=(z_R, \tau_R)=(1,4)$.
The fast wave curve from $\mathcal{U}_{R}$ has the form
$\mathcal{S}_f^\decr\,\mathcal{U}_{R}
\,\mathcal{R}_f^\decr\,\mathcal{C}o_f^\incr$.
Its reflection is drawn as the red curve in Fig.~\ref{fig:sol2}.
The dashed curve is the shock curve $\mathcal{S}_f$,
the solid curve is the rarefaction curve $\mathcal{R}_f$,
and the dash-dotted curve is the composite curve $\mathcal{C}o_f$.

To find the RPS,
we intersect the reflected fast wave curve with the intermediate surface.
In this example, $\mathcal{R}_f$, which is its own reflection,
intersects the intermediate surface of $\mathcal{S}_{s}$
(the light purple surface) at the point $\mathcal{U}_{M,f}$.
Because $\mathcal{U}_{M,f}$ belongs to
the intermediate surface of $\mathcal{S}_s$,
there is a point $\mathcal{U}_{M,s} \in \mathcal{S}_s$
that lies on the same Hugoniot$'$ curve (colored yellow)
as $\mathcal{U}_{M,f}$,
meaning that $\mathcal{U}_{M,s}$ and $\mathcal{U}_{M,f}$
have the same right state.
The RPS for the points $\mathcal{U}_L$ and $\mathcal{U_R}$
consists of the slow shock wave $\mathcal{U}_{M,s}$
and the fast rarefaction wave extending from
$\mathcal{U}_R$ to $\mathcal{U}_{M,f}$.

\begin{figure}[htpb]
\begin{center}
\begin{subfigure}{0.45\linewidth}
\includegraphics[width=\linewidth]{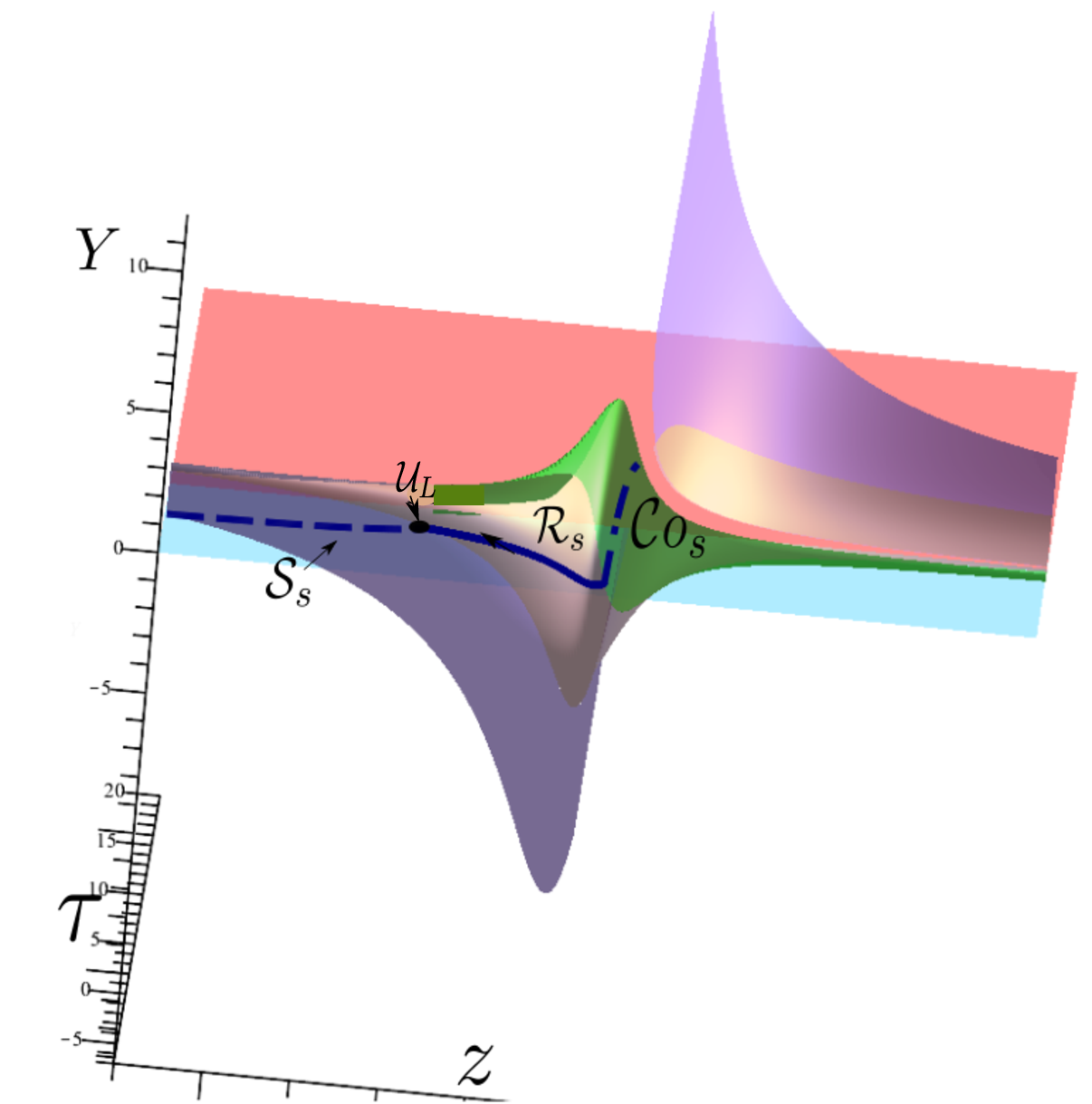}
\caption{}
\end{subfigure}
\hfil
\begin{subfigure}{0.45\linewidth}
\includegraphics[width=\linewidth]{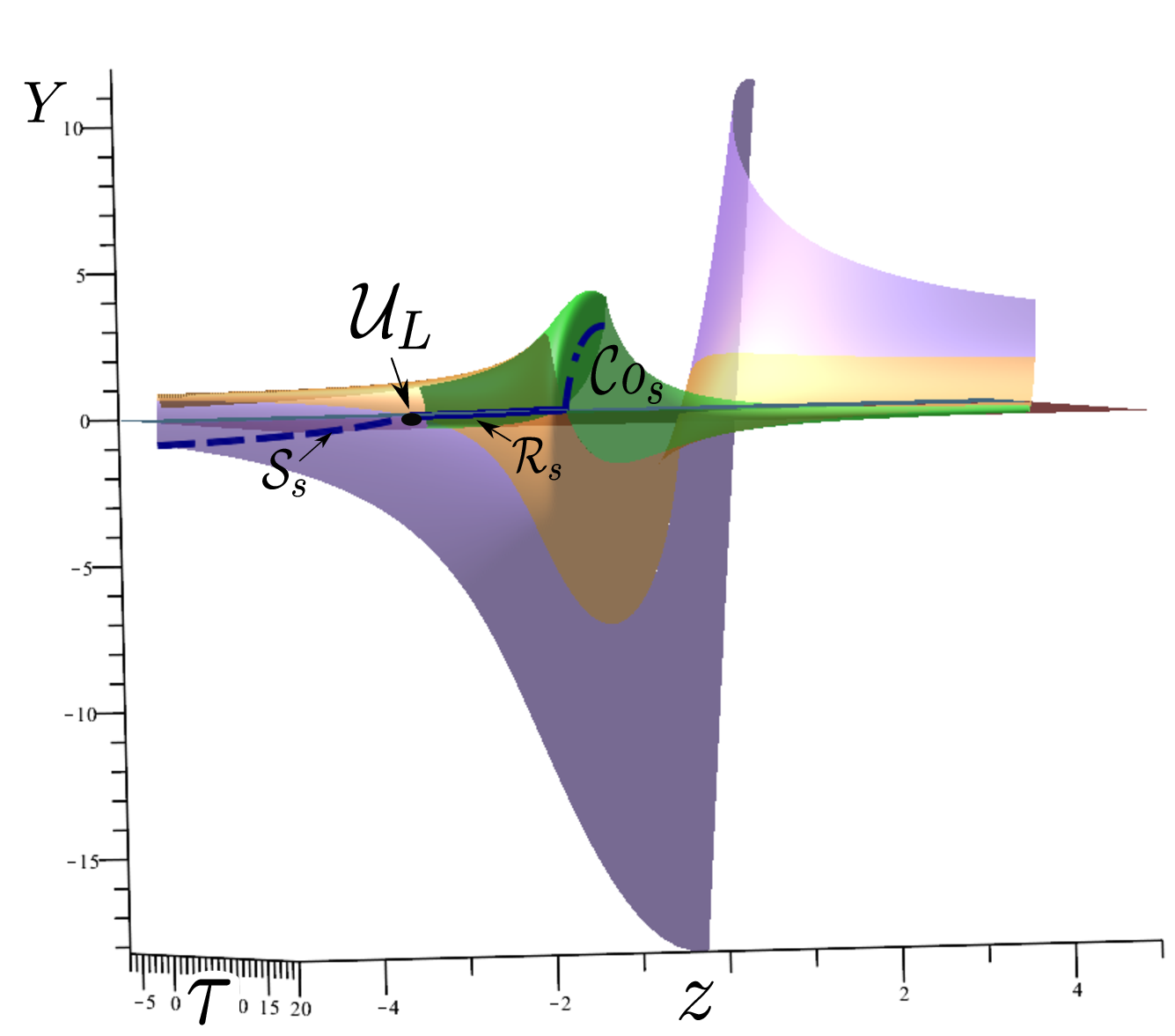}
\caption{}
\end{subfigure}
\caption[]{
The slow wave curve and its intermediate surface
for $\mathcal{U}_{L}$ in Region~$II$,
viewed from two perspectives.
The coordinates for $\mathcal{U}_L$ are
$z_{L}=-2.0$ and $\tau_{L}=-2.0$.
The slow wave curve of $\mathcal{U}_{L}$ is
$\mathcal{S}_{s}^\incr\mathcal{U}_{L}\mathcal{R}_{s}^\incr\mathcal{C}o_s^\decr$.
The portions of the intermediate surface generated by $\mathcal{S}_{s}$,
$\mathcal{R}_{s}$, and $\mathcal{C}o_f$
are colored light purple, gold, and green.
}
\label{fig:ex2-1}
\end{center}
\end{figure}

\medskip\noindent
\textbf{Example 2.} For this example,
we take $u_{L}=-0.85$ and $v_{L}=3.2$,
so that $z_{L}=-2.0$ and $\tau_{L}=-2.0$;
therefore, $\mathcal{U}_L$ belongs to Region~$II$ in $\mathcal{C}_{s}$.
The slow wave curve of $\mathcal{U}_{L}$ is
$\mathcal{S}_{s}^\incr\mathcal{U}_{L}\mathcal{R}_{s}^\incr\mathcal{C}o_s^\decr$.
Fig.~\ref{fig:ex2-1} shows this wave curve (dark blue)
and its intermediate surface in $\mathcal{W}$.
The light purple, gold, and green surfaces arise from
$\mathcal{S}_{s}$, $\mathcal{R}_{s}$, and $\mathcal{C}o_f$, respectively.

\begin{figure}[htpb]
\begin{center}
\begin{subfigure}{0.42\linewidth}
\includegraphics[width=\linewidth]{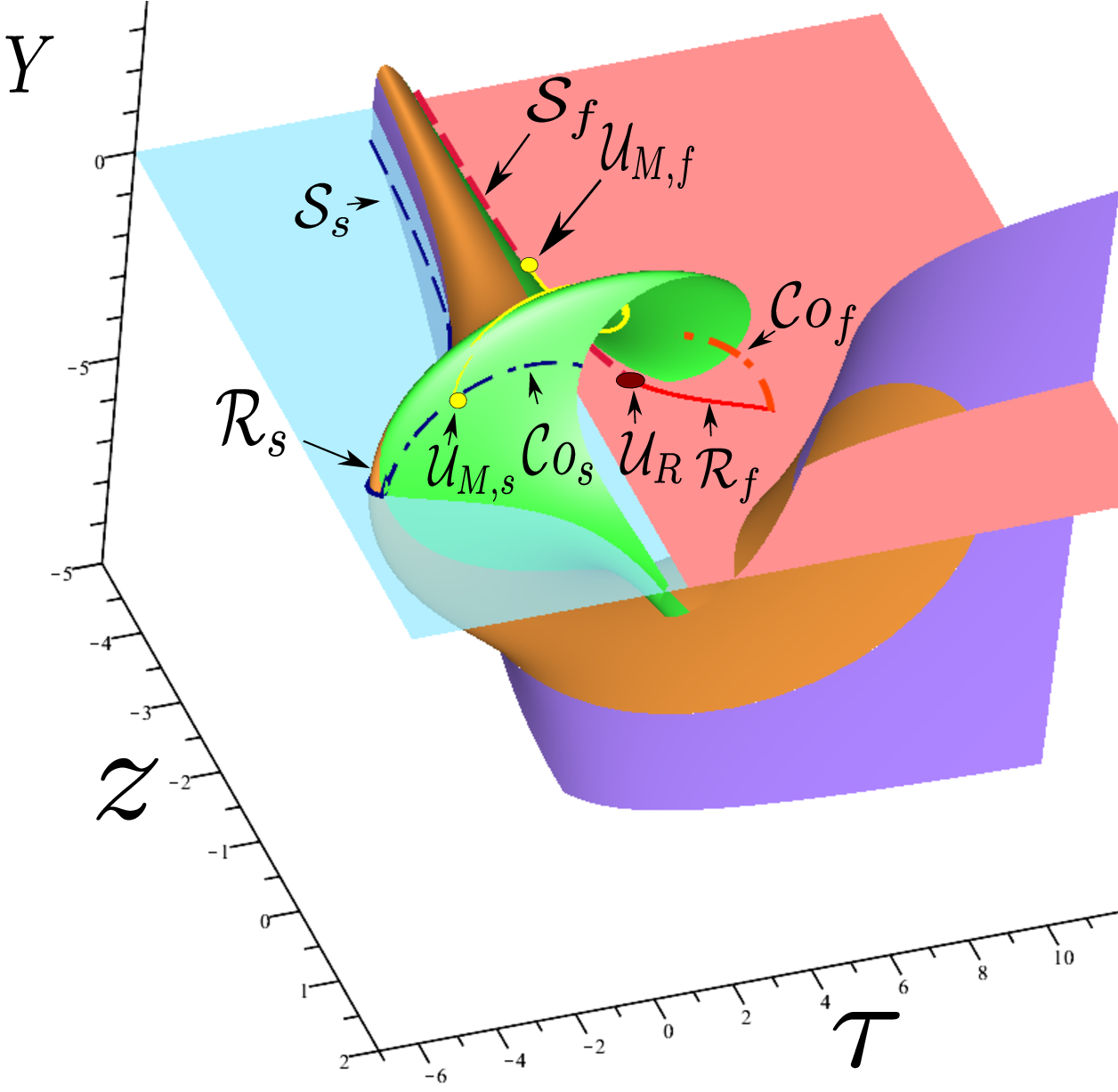}
\caption{}
\end{subfigure}
\hfil
\begin{subfigure}{0.42\linewidth}
\includegraphics[width=\linewidth]{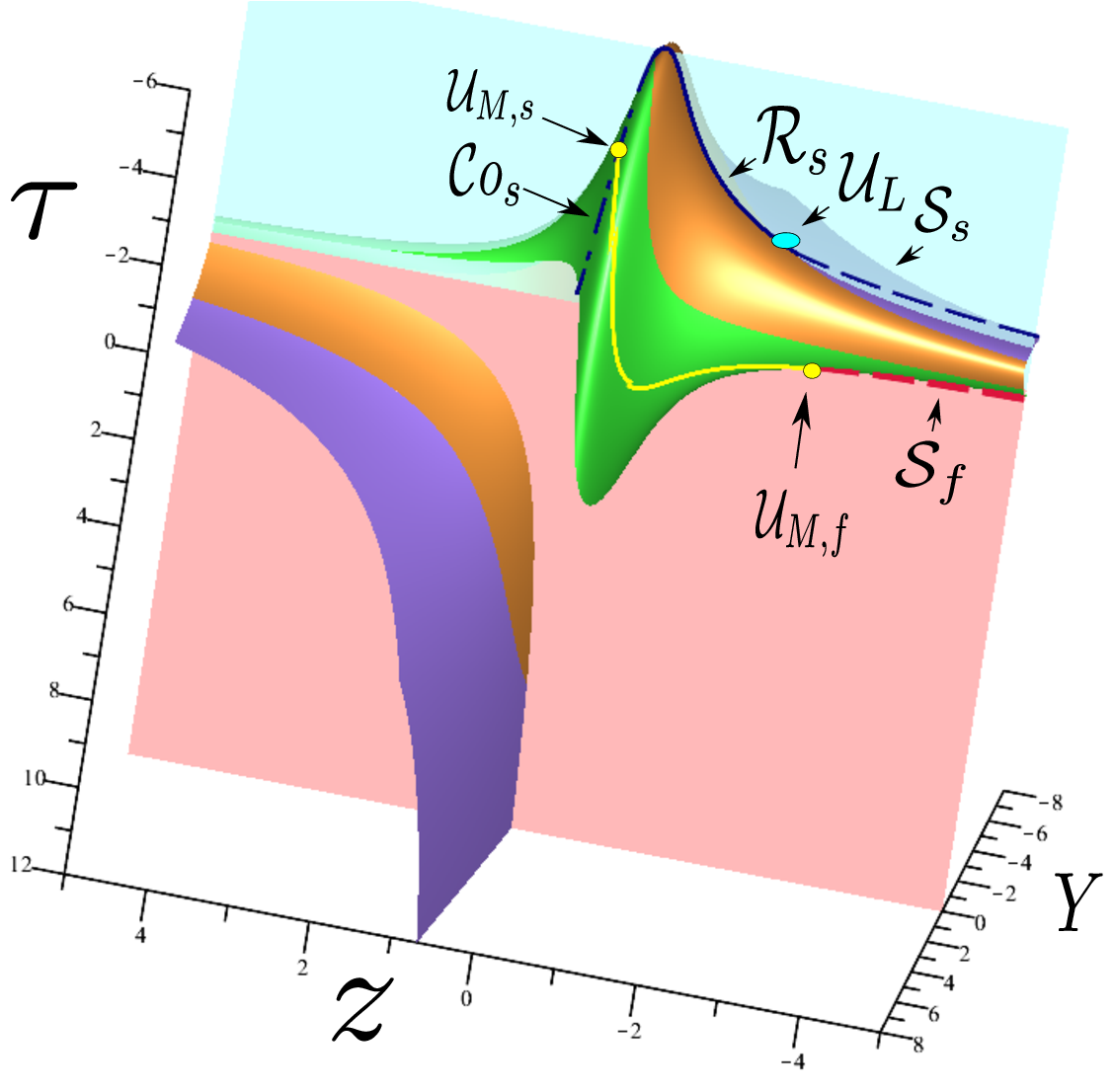}
\caption{}
\end{subfigure}
\caption[]{
Example~2 of a Riemann problem solved using wave curves in $\mathcal{W}$,
shown from two vantage points.
The fast wave curve from $\mathcal{U}_{R}$ is
$\mathcal{S}_{f}^\decr\mathcal{U}_{R}\mathcal{R}_{f}^\decr
\mathcal{C}o_{f}^\incr$; its reflection is the curve colored red.
The reflected shock curve
crosses the intermediate surface generated by $\mathcal{C}o_{s}$ (green)
at a point $\mathcal{U}_{M,f}$.
A Hugoniot$'$ curve (yellow) connects $\mathcal{U}_{M,f}$
to a point $\mathcal{U}_{M,s}$ within $\mathcal{C}o_s$.
The RPS consists of the slow composite wave
associated with $\mathcal{U}_{M,s}$
and the reflection of the fast shock wave $\mathcal{U}_{M,f}$.
}
\label{fig:ex2-3}
\end{center}
\end{figure}

With $W_R=(1.6,7.2)$, the corresponding point in $\mathcal{W}$
is $\mathcal{U}_{R}=(2.0,4.0)$.
The fast wave curve from $\mathcal{U}_{R}$,
$\mathcal{S}_{f}^\decr\mathcal{U}_{R}\mathcal{R}_{f}^\decr
\mathcal{C}o_{f}^\incr$, is reflected and drawn in red in Fig.~\ref{fig:ex2-3}.
As it happens, the reflection of $\mathcal{S}_{f}$ crosses
the intermediate surface of $\mathcal{C}o_{s}$ (the green surface)
at the point $\mathcal{U}_{M,f}$,
which is connected by a Hugoniot$'$ curve (yellow)
to a point $\mathcal{U}_{M,s} \in \mathcal{C}o_s$.
Therefore, the RPS consists of the slow composite wave
associated with $\mathcal{U}_{M,s}$
(consisting of a portion of $\mathcal{R}_s$
together with sonic slow shock wave $\mathcal{U}_{M,s}$)
together with the reflection of the fast shock wave $\mathcal{U}_{M,f}$.

\begin{figure}[htpb]
\begin{center}
\begin{subfigure}{0.45\linewidth}
\includegraphics[width=\linewidth]{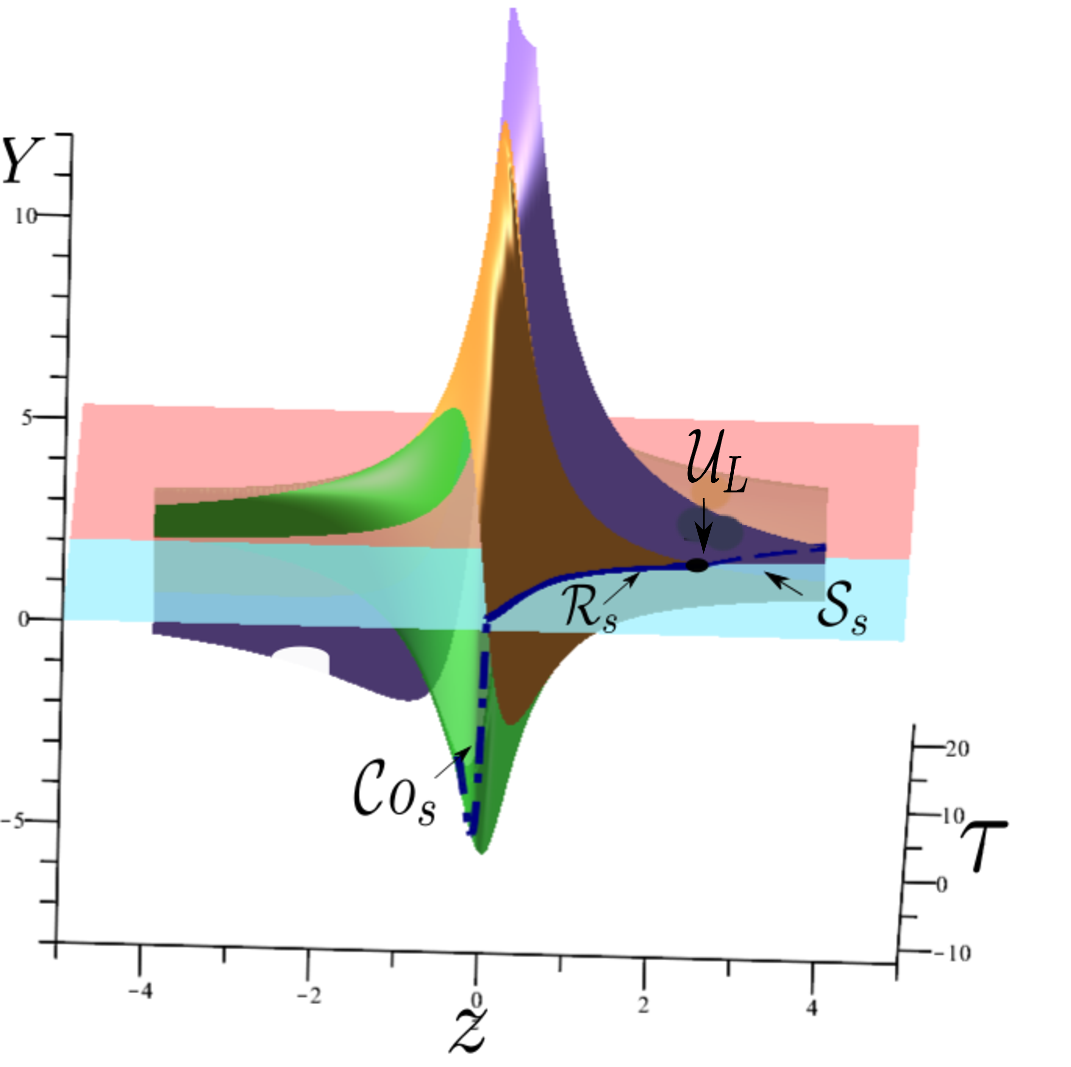}
\caption{}
\end{subfigure}
\hfil
\begin{subfigure}{0.45\linewidth}
\includegraphics[width=\linewidth]{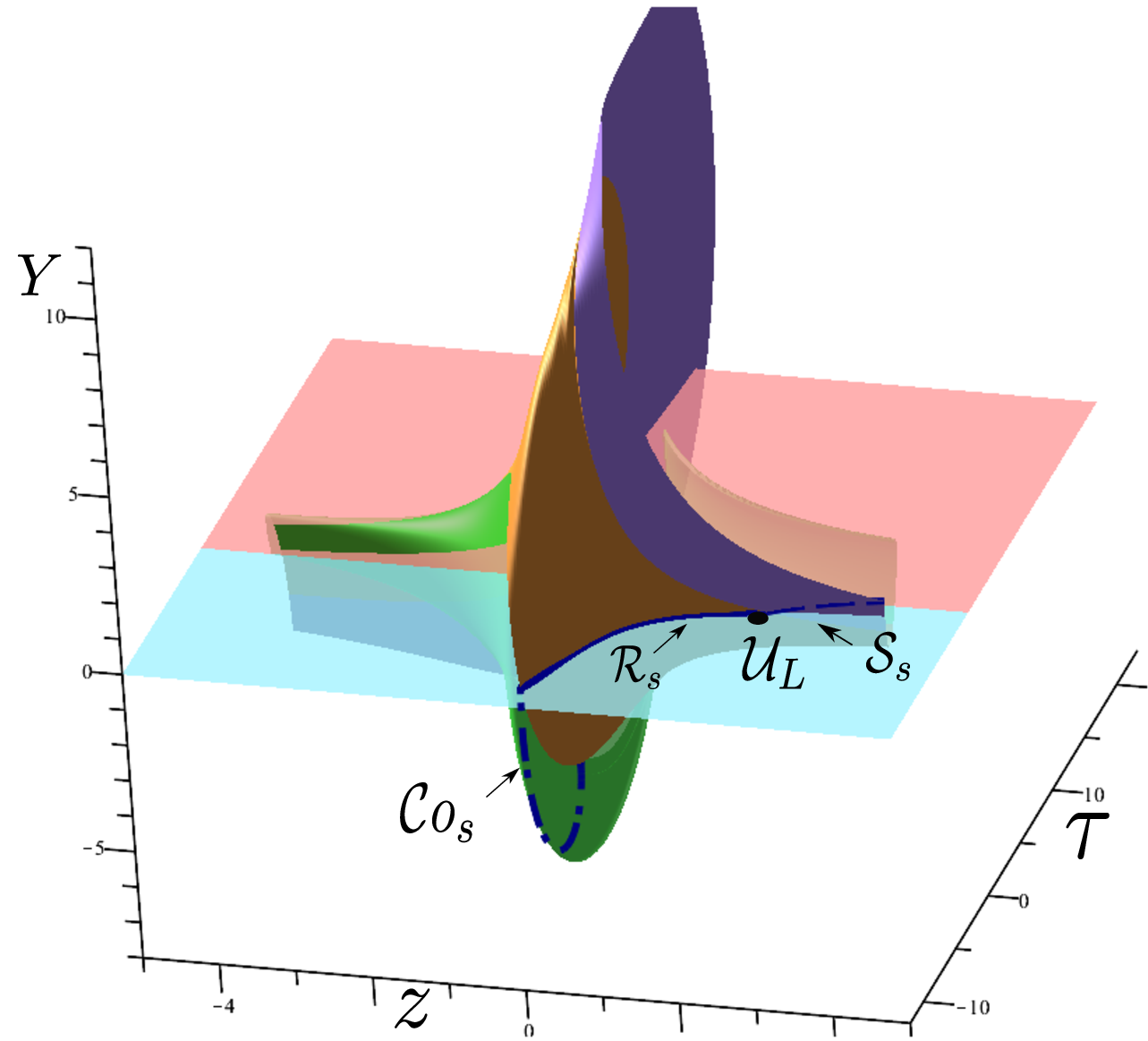}
\caption{}
\end{subfigure}
\caption[]{
The slow wave curve and its intermediate surface
for $\mathcal{U}_{L}$ in Region~$III$,
viewed from two perspectives.
The coordinates for $\mathcal{U}_{L}$
are $z_{L}=2.5$ and $\tau_{L}=-1.5$.
The slow wave curve of $\mathcal{U}_{L}$ is
$\mathcal{S}_{s}^\incr\mathcal{U}_{L}\mathcal{R}_{s}^\incr
\mathcal{C}o_s^\decr\mathcal{S}_{s}^\decr$,
but the last shock arc is not shown.
\label{fig:ex3-1}
}
\end{center}
\end{figure}

\medskip\noindent
\textbf{Example 3.}
As a final example, we take $W_{L}=(-0.898,-4.612)$;
its corresponding point $\mathcal{U}_L=(2.5,-1.5)$ belongs to Region~$III$.
The slow wave curve of $\mathcal{U}_{L}$ is
$\mathcal{S}_{s}^\incr\mathcal{U}_{L}\mathcal{R}_{s}^\incr
\mathcal{C}o_s^\decr\mathcal{S}_{s}^\decr$,
but for clarity we omit the last shock arc from the figures.
Figure~\ref{fig:ex3-1} shows the slow wave curve
and its intermediate surface.

We illustrate the construction of the RPS
with $W_R=(1.413,6.2)$,
\emph{i.e.}, $\mathcal{U}_{R}=(2.0,3.5)$.
The fast wave curve from $\mathcal{U}_{R}$ is
$\mathcal{S}_{f}^\decr\mathcal{U}_{R}\mathcal{R}_{f}^\decr
\mathcal{C}o_{f}^\incr$;
its reflection is the red curve drawn in Fig.~\ref{fig:ex3-1}.
The reflected fast composite arc crosses the intermediate surface
generated by $\mathcal{R}_{s}$ (the gold surface)
at a point $\mathcal{U}_{M,f}$,
and this point is connected to a point $\mathcal{U}_{M,s}$
on the slow rarefaction curve.
Therefore, the RPS consists of a slow rarefaction wave
extending from $\mathcal{U}_{L}$ to $\mathcal{U}_{M,s}$
and the composite wave associated with
the reflection of the sonic shock wave $\mathcal{U}_{M,f}$.

\begin{figure}[htpb]
\begin{center}
\begin{subfigure}{0.42\linewidth}
\includegraphics[width=\linewidth]{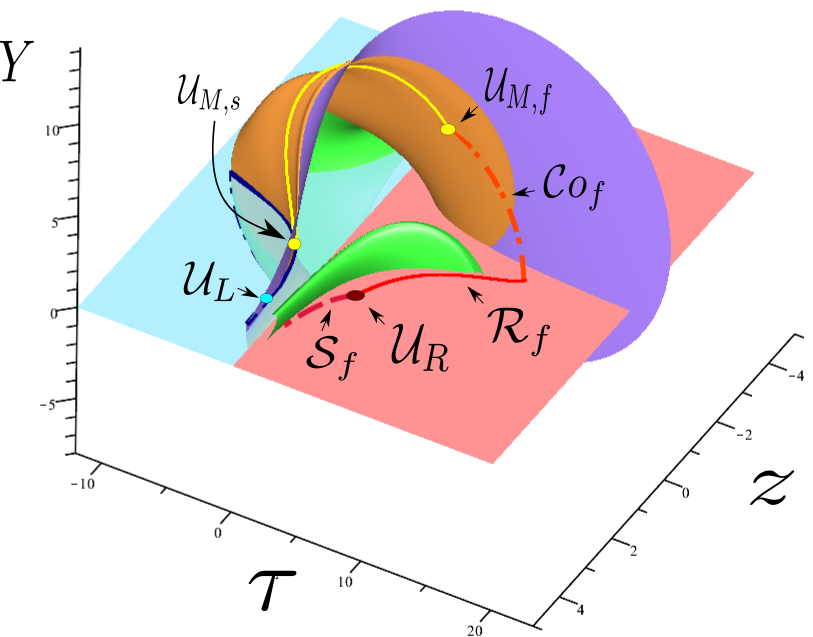}
\caption{}
\end{subfigure}
\hfil
\begin{subfigure}{0.42\linewidth}
\includegraphics[width=\linewidth]{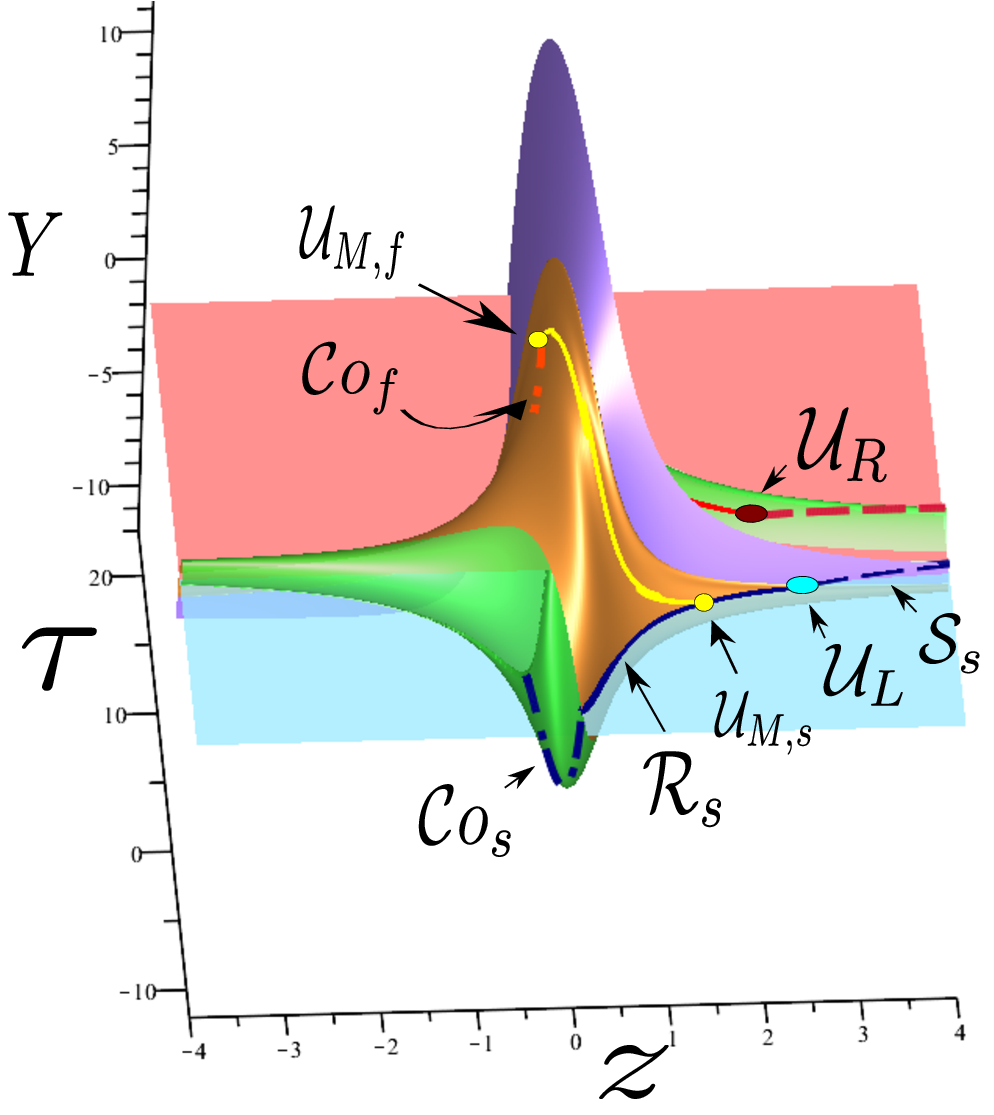}
\caption{}
\end{subfigure}
\caption[]{
Example~3 of a Riemann problem solved using wave curves in $\mathcal{W}$.
The fast wave curve from $\mathcal{U}_{R}$ is
$\mathcal{S}_{f}^\decr\mathcal{U}_{R}\mathcal{R}_{f}^\decr
\mathcal{C}o_{f}^\incr$,
and its reflection is depicted in red.
The intermediate surface
for $\mathcal{R}_{s}$ (the gold surface)
is crossed by $\mathcal{C}o_{f}$
at the point $\mathcal{U}_{M,f}$,
which is connected to a point $\mathcal{U}_{M,s} \in \mathcal{R}_s$
by a Hugoniot$'$ curve.
The RPS consists of the slow rarefaction wave
from $\mathcal{U}_L$ to $\mathcal{U}_{M,s}$
and the fast composite wave associated with
the reflection of $\mathcal{U}_{M,f}$.
}
\label{fig:ex3-3}
\end{center}
\end{figure}

\section{Concluding Remarks}
\label{sec:Concluding_Remarks}

In this study, we employ the wave manifold and utilize topological tools to
address the Riemann problem. We also analyze key structures inherent in a
hyperbolic system defined by two equations featuring an elliptic region. The
wave manifold proves to be an effective tool for tackling Riemann Problems
stemming from systems governed by quadratic flux functions. Additionally, the
construction and intersection of wave curves with important surfaces are
powerful methodologies for solving these intricate problems.

The graphical representation of important curves and surfaces within the
manifold provides valuable insights into the nature of shock waves,
rarefaction, and composite waves, aiding in formulating appropriate solutions.   The advantage of this approach is that in the three-di\-men\-sio\-nal wave manifold, Hugoniot curves yield a one-dimensional foliation in contrast with the usual cobweb of Hugoniot curves in state-space. This work extends the topological perspective to solving Riemann problems
and illustrates the procedure with a few examples.


\appendix

\section*{Appendices}
\label{sec:Appendices}

\section{Proofs of Results}
\label{sec:Proofs_of_Results}

This Appendix presents proofs of the results stated in the main text.

\subsection{Characterizing \texorpdfstring{$\mathcal{C}_{s}$}{C\_s}
and \texorpdfstring{$\mathcal{C}_{f}$}{C\_f}}
\label{subsec:Characterizing}
In Sec.~\ref{sec:sonlifs},
we defined and characterized $Son'_{s}$ and $Son'_{f}$.
As presented in Sec.~\ref{sec:lax}, shock curves either start on
$\mathcal{C}_{s}$ or on $Son'_{s}$.
So we need to characterize $\mathcal{C}_{s}$ and $\mathcal{C}_{f}$.

The speed $\sigma$ along a Hugoniot curve
of a point $(z_{0},\tau_{0},Y_{0})$ is obtained
by substituting $u_{0}(z_{0},\tau_{0},Y_{0})$
and $v_{0}(z_{0},\tau_{0},Y_{0})$,
as obtained from Eq.~\eqref{eq:parahugpont},
into the expression $\sigma=\sigma(u_{0}, v_{0},z)$
given by \eqref{eq:speed}.
We obtain
\begin{equation}
\sigma_{hug} = \frac{s_{p3}z^{3}-s_{p2}z^{2}-s_{p1}z+s_{p0}}
{2b_{1}\vartheta_0[z^{2}\theta_-+1]},
\label{eq:velhug}
\end{equation}
where $s_{p3}=b_{1}\theta_+\left\{Y_{0}\vartheta_0
+2c[1+z_{0}\tau_{0}\vartheta_0]\right\}$,
$s_{p2}=\theta_+[b_{1}z_{0}Y_{0}\vartheta_0+2c(2z_{0}
+\tau_{0}z_{0}^{4}-\tau_{0})]$,\break
$s_{p1}=b_{1}\left\{Y_{0}\vartheta_0
+2c[z_{0}\tau_{0}\vartheta_0-2z_{0}^{2}-1]\right\}$,
$s_{p0}=b_{1}z_{0}Y_{0}\vartheta_0+2c(2z_{0}+\tau_{0}z_{0}^{4}-\tau_{0})$,
and $\vartheta_0 =z_{0}^{2}+1$, $\theta_+=b_1+1$, $\theta_-=b_1-1$.

The equation for the shock speed $\sigma$ along a Hugoniot curve of a point
$(z_{0},\tau_{0},0)$ on the characteristic surface (characteristic speed) is
obtained by substituting $Y_{0}=0$ into \eqref{eq:velhug}:
\begin{equation}
\sigma_{ch}=\frac{s_{ch3}z^{3}- s_{ch2}z^{2}
- s_{ch1} z+s_{ch0}}{b_{1}\vartheta_0[\theta_-z^{2}+1]},
\label{eq:velcar}
\end{equation}
where $s_{ch3}=b_{1}c[\tau_{0}z_{0}\vartheta_0+1]\theta_+$,
$s_{ch2}=c(\tau_{0}z_{0}^{4} +2z_{0}-\tau_{0})\theta_+$,
$s_{ch1}=b_{1}c(\tau_{0}z_{0}\vartheta_0-2z_{0}^{2}-1)$,
and $s_{ch0}=c(\tau_{0}z_{0}^{4} +2z_{0}-\tau_{0})$.

\subsection{}
\begin{proof}[Proof of Prop.~\ref{th:slow_fast}]
A Hugoniot curve of a point on the characteristic has parametric equations
\begin{equation}\label{eq:hugch}
z=z,\ \tau=\displaystyle
\frac{D_{0}z^{3}+E_{0}z^{2}+F_{0}z+G_{0}}
{c\left(z_{0}^{2}+1\right)[(b_{1}-1)z^{4}+b_{1}z^{2}+1]},
\ Y=\displaystyle
\frac{A_{0}z^{2}+B_{0}z+C_{0}}{(z_{0}^{2}+1)[(b_{1}-1)z^{2}+1]},
\end{equation}
where $A_{0}$,\ $B_{0}$,\ $C_{0}$,\ $D_{0}$,\ $E_{0}$,
\ $F_{0}$ and $G_{0}$ are obtained from \eqref{eq:hugponto} taking $Y_{0}=0$.

Solving $Y=0$ in \eqref{eq:hugponto} we get
$z_{1}=-[\tau_{0}(1+z_{0}^{2})-z_{0}]/[\tau_{0}z_{0}(1+z_{0}^{2})+1]$ and
$z_{0}$, coordinates of $z$ at the intersection points and $\mathcal{C}$. It
is easy to see that $z_{0}=z_{1}$ if and only if $\tau_{0}=0$, characterizing
the \emph{coincidence curve}.

Replacing $z$ by $z_{1}$ in the expression of $\tau$ in Eqs.~\eqref{eq:hugch},
we get the expression of $\tau_{1}$. A straightforward computation shows that
$\tau_{0}\cdot \tau_{1}<0$, so, an intersection point lies on $\mathcal{C}_{s}$
and the other one on $\mathcal{C}_{f}$.

Let $\sigma_{ch0}$ and $\sigma_{ch1}$ be the characteristic speeds for
$z=z_{0}$ and $z=z_{1}$. Straightforward computations give
$\sigma_{ch0}-\sigma_{ch1}=c(z_{0}^{2}+1)\tau_{0}$. It follows that $
\sigma_{ch0}>\sigma_{ch1}$ if and only if $\tau_{0}>0$. So, the value of
$\sigma$ at the point on $\mathcal{C}_{f}$ is larger than the value of $\sigma$
at the point on $\mathcal{C}_{s}$.
\end{proof}

\subsection{}
\begin{figure}[htpb]
\begin{center}
\begin{subfigure}{0.48\linewidth}
\includegraphics[width=\linewidth]{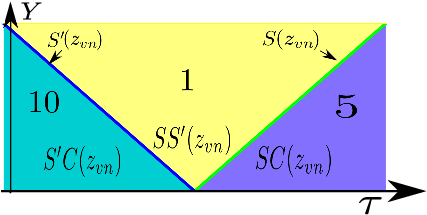}
\caption{}
\end{subfigure}
\hfil
\begin{subfigure}{0.48\linewidth}
\includegraphics[width=\linewidth]{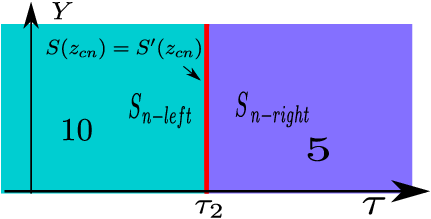}
\caption{}
\end{subfigure}
\begin{subfigure}{0.48\linewidth}
\includegraphics[width=\linewidth]{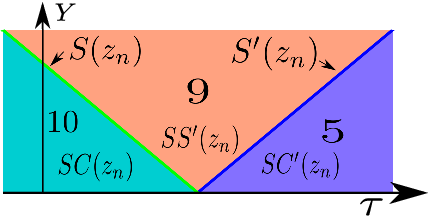}
\caption{}
\end{subfigure}
\hfil
\begin{subfigure}{0.48\linewidth}
\includegraphics[width=\linewidth]{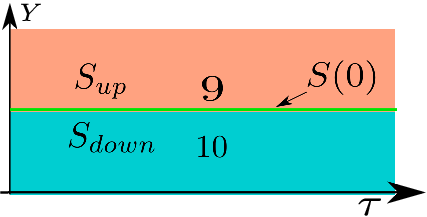}
\caption{}
\end{subfigure}
\begin{subfigure}{0.48\linewidth}
\includegraphics[width=\linewidth]{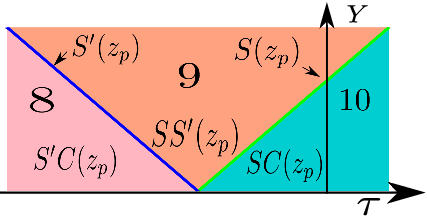}
\caption{}
\end{subfigure}
\hfil
\begin{subfigure}{0.48\linewidth}
\includegraphics[width=\linewidth]{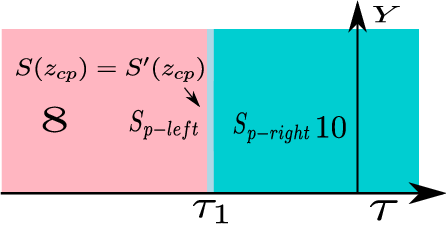}
\caption{}
\end{subfigure}
\begin{subfigure}{0.48\linewidth}
\includegraphics[width=\linewidth]{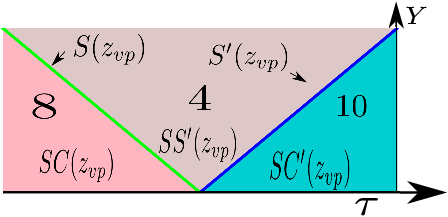}
\caption{}
\end{subfigure}
\caption[]{
Constant-$z$ slices of $\mathcal{W}$.
The numbered regions are as indicated in Fig.~\ref{fig:02}.
(a)~$\Pi_{z_{{vn}}} \cap \mathcal{W}$.
(b)~$\Pi_{z_{{cn}}} \cap \mathcal{W}$.
Here $\tau_2$ corresponds to $z_{{cn}}$ through Eq.~\eqref{eq:tdoublesonic},
so that the vertical red line is the double sonic locus.
(c)~$\Pi_{z_{{n}}} \cap \mathcal{W}$.
(d)~$\Pi_{0} \cap \mathcal{W}$.
(e)~$\Pi_{z_{{p}}} \cap \mathcal{W}$.
(f)~$\Pi_{z_{{cp}}} \cap \mathcal{W}$.
Here $\tau_1$ corresponds $z_{{cp}}$ through Eq.~\eqref{eq:tdoublesonic},
so that the vertical line in light blue is the double sonic locus.
(g)~$\Pi_{z_{{vp}}} \cap \mathcal{W}$ at $z_{{vp}}$.
}
\label{fig:cor}
\end{center}
\end{figure}

\begin{figure}[htpb]
\begin{center}
\begin{subfigure}{0.48\linewidth}
\includegraphics[width=\linewidth]{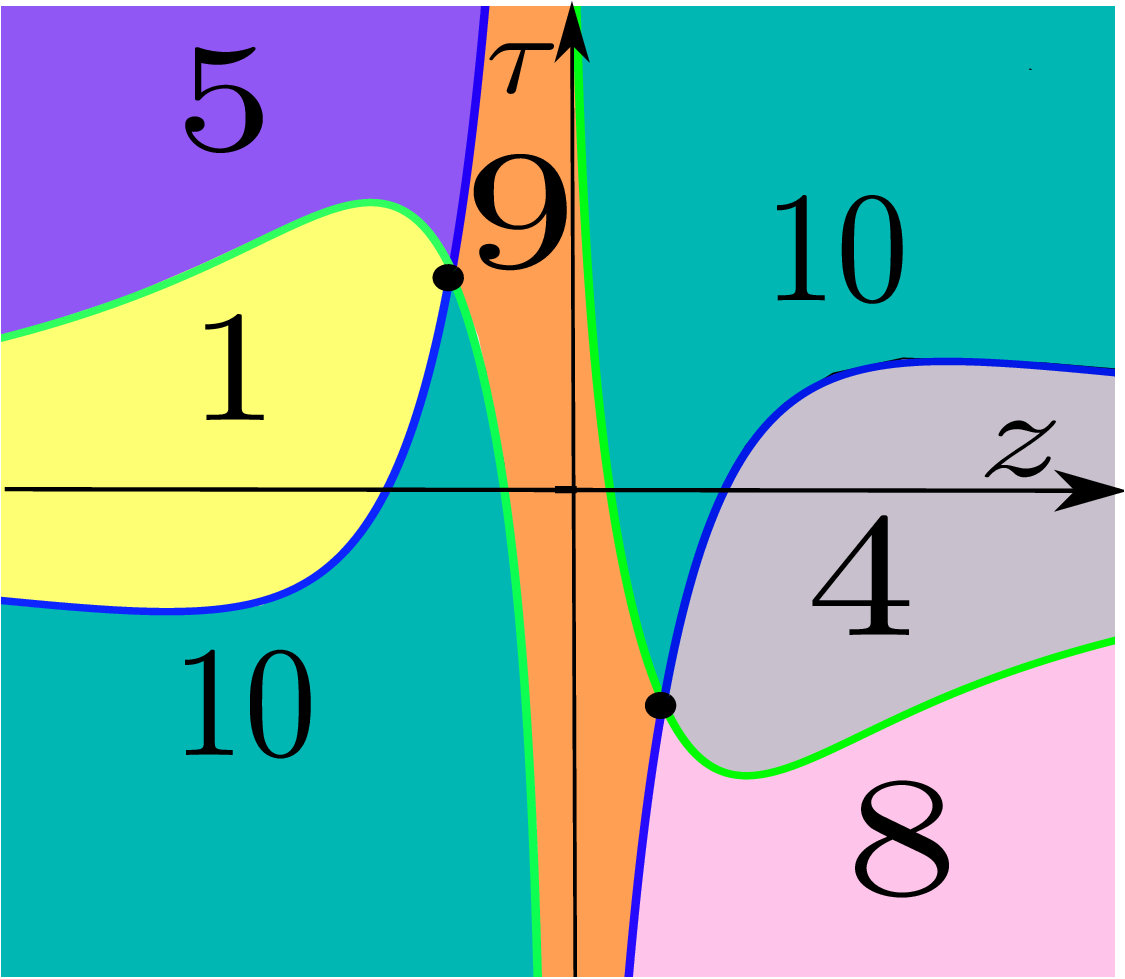}
\caption{}
\end{subfigure}
\begin{subfigure}{0.48\linewidth}
\includegraphics[width=\linewidth]{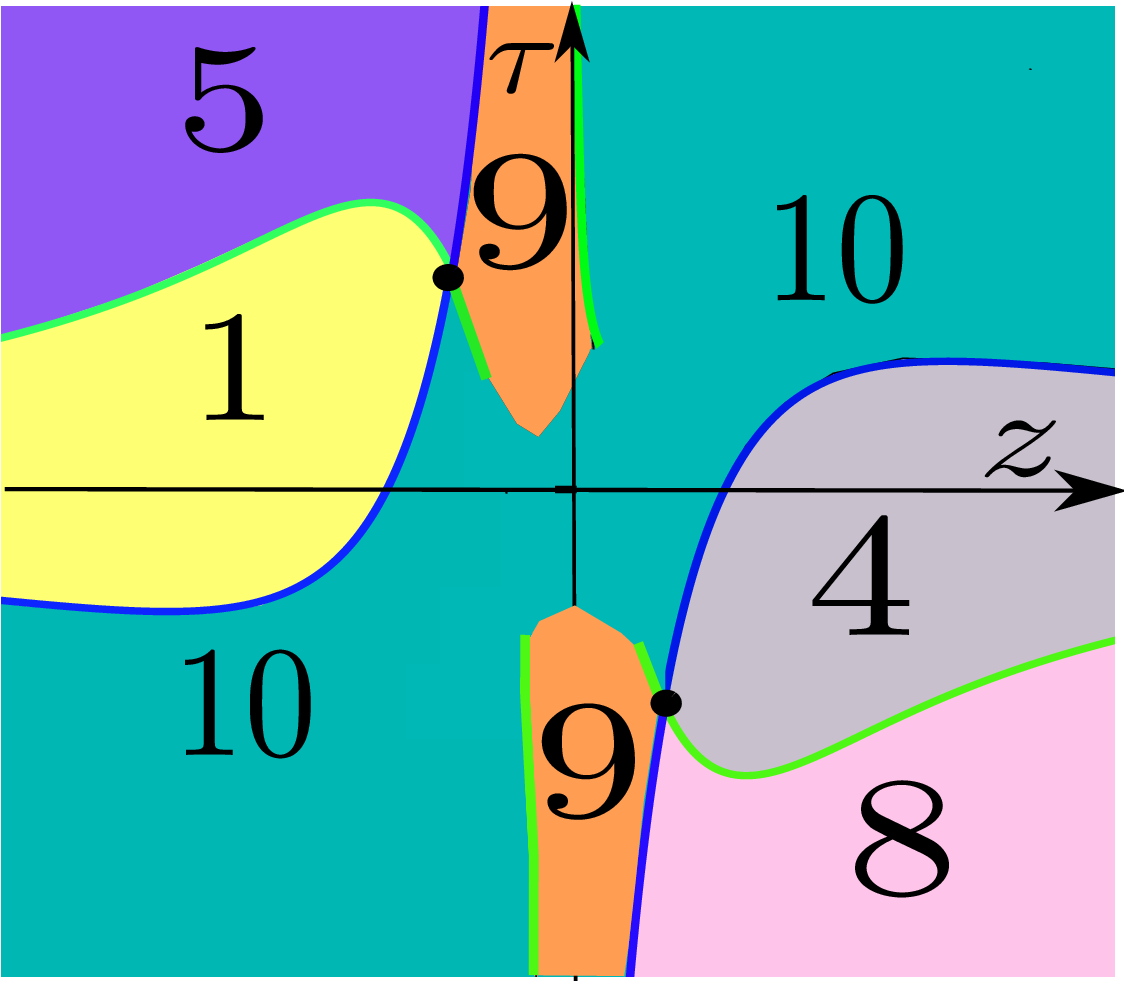}
\caption{}
\end{subfigure}
\caption{(a)~This figure illustrates
a cross-section parallel to the $z-\tau$ plane,
with $Y$ held constant at a positive value.
The cross-section for $Y$ is slightly larger than $Y^*$,
where $Y^*$ is the maximum value of coordinate $Y$ for $SCC'$.
There are six regions ($1$, $4$, $5$, $8$, $9$, $10$).
(b)~This figure represents the same cross-section,
albeit with $Y$ being smaller than $Y^*$. Notice that region 10 is connected.
}
\label{fig:cor4}
\end{center}
\end{figure}

\begin{proof}[Proof of Prop.~\ref{proposition:divideM}]
Since $\mathcal{C}$ is just the horizontal plane $Y=0$ and the reflection $(z,
\tau, Y) \to (z, \tau, -Y)$ interchanges $Son$ and $Son'$, it is sufficient to
divide the upper half space $Y>0$, and we have a symmetric division in the
$Y<0$ half space. From now on, we focus on the $Y>0$ half space,
$\mathcal{W}^+$, and show that $\mathcal{W}^+$ minus the union of $S$ and $S'$
has six connected components.

We slice $\mathcal{W}^+$ by the $\Pi_z$ planes ($z=\text{constant}$) and see
how each $\Pi_z$ is divided into 2-dimensional regions by $\mathcal{C}$, $Son$,
and $Son'$ and see how these 2-dimensional regions generate 3-dimensional
regions as $z$ varies.

In $\Pi_z$, we utilize the natural $\tau$ and $Y$ axes inherited from
$(z, \tau, Y)$ space.
The intersections of $\Pi_z$ with $\mathcal{C}$, $Son$, and $Son'$ are,
respectively, the $\tau$ axis (horizontal) and two half lines $S(z)$
and $S'(z)$ with equations,
\begin{equation}
m(z)Y + n(z)\tau + p(z) = 0
\end{equation}
and its reflection
\begin{equation}
-m(z)Y + n(z)\tau + p(z) = 0,
\end{equation}
where $m(z)=(b_1+1)z^4+b_1z^2-1$, $n(z)=-( (b_1+1)z^4+(b_1+4)z^2+3 )cz$,
and $p(z)=2c((b_1-1)z^2-1 )$.
These are just equations \eqref{eq:son} and \eqref{eq:son'}.

It is easy to see that $S(z)$ is vertical if $z^2=1/(b_1+1)$,
with solution given by
\begin{equation}
z_{crit1}=z_{cp}= 1/\sqrt{b_{1}+1} \quad\text{and}\quad
z_{crit2}=z_{cn}=-1/\sqrt{b_{1}+1} \label{doublecon}
\end{equation}
and horizontal if $z=0$ (actually $S(0)$ has equation $Y=2c$ and $S'(0)$ has
equation $Y=-2c$, which places it out of $\mathcal{W}^+$).

The $\tau$ corresponding to $z_{crit1}$ and $z_{crit2}$ are
\begin{equation}
\tau_{1}=-b_{1} \sqrt{b_{1}+1}/[2(b_{1}+2)]
\quad\text{and}\quad
\tau_{2}= b_{1} \sqrt{b_{1}+1}/[2(b_{1}+2)].\label{eq:tdoublesonic}
\end{equation}

Varying $z$, we see that the angular coefficient of $S(z)$ is positive if
$z<z_{cn}$ and changes sign as $z$ crosses $z_{cn}$, $0$, and $z_{cp}$, where
$z_{cp}$ and $z_{cn}$ are defined in \eqref{doublecon}. The intersection
point of $S(z)$ with the $\tau$ axis has positive coordinate $\tau$ if $z<0$
and negative coordinate $\tau$ if $z>0$. This $\tau$ coordinate tends to
infinity if $z$ tends to $0^-$ and tends to $-\infty$ if $z$ tends to $0^+$.

We make seven representative slices:
\begin{align*}
 &z=z_{{vn}}<z_{cn};\; z=z_{ {cn}};\; z_{cn}<z=z_n<0;\; z=0;\\
 &0<z=z_p<z_{cp};\; z=z_{{cp}};\; z_{cp}<z=z_{ {vp}}.
\end{align*}
This notation is supposed to be mnemonic: ${vn}$
stands for very negative, ${cn}$ for critical negative,
${n}$ for negative, ${p}$ for positive, ${cp}$ for critical positive,
${vp}$ for very positive. Since the angular coefficient of $S'(z)$
is minus the angular coefficient of $S(z)$,
and $S(z)$ and $S'(z)$ intersect on the $\tau$ axis,
we have the following pictures for the seven slices mentioned above.

Let us call $SS'(z)$ the 2D region of $\Pi_z$ bounded by the lines $S(z)$
and $S'(z)$. Call $SC(z)$ the region of $\Pi_z$ bounded by $S(z)$
and the $\tau$ axis, and call $S'C(z)$ the region bounded by $S'(z)$
and the $\tau$ axis. Boundaries not included.

Notice that for $z=z_{{cn}}$, $S(z_{{cn}})=S'(z_{{cn}})$,
and there are only two regions,
which are called $SC_{{n}-\text{left}}$ and $SC_{{n}-\text{right}}$.
They could just as well be called
$S'C_{\text{left}}$ and $S'C_{\text{right}}$.
In the same way, for $z=z_{{cp}}$,
there are only two regions $SC_{p-\text{left}}$ and $SC_{p-\text{right}}$.
Also, for $z=0$ there are only the 2 regions:
``above $S(0)$" and ``between $S(0)$ and the $\tau$ axis".
Let us call them $S_{\text{up}}$ and $S_{\text{down}}$.

As $z$ varies, each 2D region defined above generates a 3D region:
\begin{align*}
&SS'_{{vn}} = \cup (SS'(z), z<z_{cn}), \\
&SS'_{{n}} = \cup (SS'(z), z_{cn}<z<0), \\
&SS'_{{p}} = \cup (SS'(z), 0<z<z_{cp}), \\
&SS'_{{vp}} = \cup (SS'(z), z_{cp}<z),\\
&SC_{{vn}} = \cup (SC(z), z<z_{cn}), \\
&SC_{{n}} = \cup (SC(z), z_{cn}<z<0), \\
&SC_{{p}} = \cup (SC(z), 0<z<z_{cp}), \\
&SC_{{vp}} = \cup (SC(z), z_{cp}<z).
\end{align*}
In the same way, we define
$S'C_{{vn}}$, $S'C_{{n}}$, $S'C_{{p}}$, and $S'C_{{vp}}$.

Let us check that some of these 12 regions connect, generating the
decomposition of $\mathcal{W}^+ \setminus Son \cup Son'$ into six
connected components. We do this by taking a point in each of the twelve
regions and seeing if it can be connected by a continuous path to one of its
neighbor regions without crossing its boundary.

Drawing a path connecting any two points in the same 3D region is easy. Let us
see what happens when $z$ moves across $z_{{cn}}$, $z_{{cp}}$, or $0$.

We have three different cases, according to whether
we begin with a point in $SC(z_{{vn}})$, $SS'(z_{{vn}})$, or $S'C(z_{{vn}})$. As
$z$ approaches $z_{{cn}}$, $SS'(z)$ becomes thinner and disappears at
$z_{{cn}}$. So, if the beginning point is in $SS'(z_{{vn}})$, the region
$SS'_{{vn}}$ will not connect to any other. If the beginning point is in
$SC(z_{{vn}})$, the path will arrive at $SC(z_{{cn}})=SC_{{n}-\text{right}}$ and
proceed into $S'C_{{n}}$ (recall that at $z_{{cn}}$, there is a position switch
between the half lines $SC(z)$ and $S'C(z)$). In the same fashion, if the
beginning point is in $S'C(z_{{vn}})$, the path enters $SC_{{n}}$, via
$SC_{{n}-\text{left}}$.
Schematically we have, see Figs.~\ref{fig:cor}-\ref{fig:cor4}:
\begin{align*}
&SC_{{vn}} \to S_{{n}-\text{right}} \to S'C_{{n}}, \\
&SS'_{{vn}} \to \text{nothing},\\
&S'C_{{vn}} \to S_{{n}-\text{left}} \to SC_{{n}}.
\end{align*}

Let us continue from a point in $\Pi_{z_{{n}}}$. As $z$ goes to $0^-$, the
intersection point of the half-lines $S(z)$ and $S'(z)$ goes to infinity, and
$S(z)$ becomes the line $Y=2$. So the region $SC$ connects to
$S_{\text{down}}$, and the region $SS'_{{n}}$ connects to $S_{\text{up}}$.
Also, the half line $S'(z)$ is pushed out of the picture, so the region
$S'C_{{n}}$ does not connect to anything, and we have the diagram:
\begin{align*}
&SC_{{vn}} \to S_{{n}-\text{right}} \to S'C_{{n}} \to \text{nothing}, \\
&SS'_{{vn}} \to \text{nothing} \quad SS'_{{n}} \to S_{\text{up}}, \\
&S'C_{{vn}} \to S_{{n}-\text{left}} \to SC_{{n}} \to S_{\text{down}}.
\end{align*}

We can continue to $z_{{vp}}$ and complete the diagram, where we count six
connected components: 2 on the first line, three on the 2nd line, and one on
the 3rd line. The ``nothing'' separates the components:
\begin{align*}
&SC_{{vn}} \to S_{{n}-\text{right}} \to S'C_{{n}} \to \text{nothing}
\quad SC_{{p}} \to S_{p-\text{right}} \to S'C_{{vp}}, \\
&SS'_{{vn}} \to \text{nothing}
\quad SS'_{{n}} \to S_{\text{up}} \to SS'_{{p}} \to \text{nothing}
\quad SS'_{{vp}}, \\
&S'C_{{vn}} \to S_{{n}-\text{left}} \to SC_{{n}}
\to S_{\text{down}} \to S'C_{{p}} \to S_{p-\text{left}} \to S'C_{{vp}}.
\end{align*}

In Figs.~\ref{fig:cor}--\ref{fig:cor4} we illustrate the connections,
where colors indicate connected regions.
These regions are the same regions
shown in Fig.~\ref{fig:02}.
\end{proof}

\subsection{}
\begin{proof}[Proof of Prop.~\ref{th:Tf}]
We recall its statement:
\emph{The surface $SCC$ is tangent to $\mathcal{C}$ along $\mathcal{E}$
and to $Son'$ along $ECC\,'$.}

From Eqs.~\eqref{eq:hugponto} with $\tau_0 = 0$ and $Y_0 = 0$ we obtain
\begin{equation}
\begin{split}
z&=z, \\
\tau&=2\,\dfrac{-cz^{3}-2cz_0z^{2}+c(2+2z_{0}^{2}-b_{1}z_0^2)z-2cz_0}
{\vartheta_{0} c[(b_{1}-1)z^{4}+b_{1}z^{2}+1]}, \\
Y&=2\displaystyle \frac{-cz^{2}+2cz_0z +cz_0}
{\vartheta_{0}[(b_{1}-1)z^{2}+1]},
\end{split}
\label{eq:aphugponto}
\end{equation}
where $\vartheta_0=z_0^2+1$,
which can be thought as parametric equations
for the $SCC$ in parameters $z$ and $z_0$.
If we differentiate in $z$ and $z_0$,
we obtain vectors $V_1$ and $V_2$.
To verify the tangency with $\mathcal{C}$,
it is simple to check that for $z=z_0$,
the $Y$ components of $V_1$ and $V_2$ are zero.

To see tangency to $Son'$,
$V_1$ and $V_2$ must be orthogonal to $\operatorname{grad}(son')$
at the intersection of $SCC$ and $Son'$,
recall that $son'=0$ is the equation of the sonic\,$'$ surface \eqref{eq:son'}.
A straightforward calculation to find the dot products
$V_1\cdot \operatorname{grad}(son')$
and $V_2\cdot \operatorname{grad}(son')$
and verify that they are both zero.

First, we solve equations \eqref{eq:aphugponto} and \eqref{eq:son'} in $\tau$
and $Y$, obtaining a relation between $z$ and $z_0$,
\begin{equation}\label{eq:A9new}
z_{0}=\frac{(b_{1}+1)z^{2}+1}{2z}.
\end{equation}
Then, we replace $z_0$ by this value and $\tau$ and $Y$ by their values from
\eqref{eq:aphugponto} in $V_1\cdot \operatorname{grad}(son')$ and $V_2\cdot
\operatorname{grad}(son')$. This results in expressions in $z$, which, when
simplified, turns out to be zero. Computations are long, so we used
Maple\textsuperscript{\textregistered}.
\end{proof}

If we substitute in \eqref{eq:aphugponto}, $z_0$ by its value from
\eqref{eq:A9new} we get the following parametrization for $ECC'$

\begin{equation}
\left \{
\begin{array}{l}
z=z,
\tau=-\dfrac{(b_{1}+2)z[(b_{1}-1)z^{2}+1)]}
{\vartheta[\theta_+^{2}z^{4}+2(b_{1}+3)z^{2}+1]},
Y=-\dfrac{2c[(b_{1}-1)z^{2}+1]}{\theta_+^{2}z^{4}+2(b_{1}+3)z^{2}+1}
\end{array}
\right\},
\label{eq:parsonli}
\end{equation}
where $\vartheta=z^2+1$ and $\theta_+=b_1+1$.

\subsection{}
\label{sec:propson}
\begin{proof}[Proof of Prop. \ref{propson}]
For a point $(z_{0},\tau'_{0},Y_{0})$ on $Son'$,
we obtain from Eq.~\eqref{eq:son'} that
\begin{equation}\label{eq:velponto}
\tau'_{0}=\frac{Y_{0}\,(z_{0}^{2}+1)\,[(b_{1}+1)\,z_{0}^{2}-1]
-2\,c\,[(b_{1}-1)\,z_{0}^{2}+1]}
{2\,c\,z_{0}\,(z_{0}^{2}+1)\,[(b_{1}+1)\,z_{0}^{2}+3]}.
\end{equation}
The shock speed at this point, $\sigma_{son'}(z_{0}, \tau_{0}',Y_{0})$,
is obtained putting $\tau_{0}=\tau_{0}'$ into Eq.~\eqref{eq:speed},

\begin{equation}
\sigma_{son'}(z_{0},\tau_{0}',Y_{0})
=\frac{[(b_{1}+1)z_{0}^{2}-1]^{2}Y_{0}+6c(b_{1}+1)z_{0}^{2}+2c}
{2b_{1}z_{0}[(b_{1}+1)z_{0}^{2}+3]}.
\label{eq:speedson'}
\end{equation}

Parametric equations for the Hugoniot curve of $(z_{0},\tau_{0}',Y_{0})$
are obtained changing $\tau_{0}$ into $\tau_{0}'$ in Eqs.~\eqref{eq:hugponto},
giving $(z, \tau_{hug}(z),Y_{hug}(z))$.
Solving $Y_{hug}=0$, we obtain the $z$ coordinates of
the intersection points of Hugoniot curve with $\mathcal{C}$:
\begin{equation}
\text{$z_{\mathcal{C}1}=\frac{(b_{1}+1)z_{0}^{2}+1}{2z_{0}}$
and $z_{\mathcal{C}2}
=-\frac{2\,(Y_{0}-c)\,z_{0}}{(b_{1}+1)\,Y_{0}\,z_{0}^{2}+Y_{0}+2c}$}.
\end{equation}
Substituting $z_{\mathcal{C}1}$ and $z_{\mathcal{C}2}$
into the expression for $\tau_{hug}(z)$,
we obtain the $\tau$ coordinates of the intersection points:
\begin{align}
\tau_{\mathcal{C}1}&=\frac{2z_{0}\left\{
[(b_{1}+1)^{2}z_{0}^{4}+2(b_{1}+3)z_{0}^{2}+1]Y_{0}
+2c[(b_{1}-1)z_{0}^{2}+1]\right\}}
{c[(b_{1}+1)z_{0}^{2}+3][(b_{1}+1)^{2}z_{0}^{4}+2(b_{1}+3)z_{0}^{2}+1]} \\
\tau_{\mathcal{C}2}&=-\frac{(A_{\tau_{\mathcal{C}2}}\,Y_{0}
+B_{\tau_{\mathcal{C}2}})
[(\theta_+\,z_0^2+1)Y_{0}+2c]^{2}}
{2\,c\,z_{0}\,(\theta_+\,z_0^2+3)(C_{\tau_{\mathcal{C}2}}\,Y_{0}^{2}
+ D_{\tau_{\mathcal{C}2}}\,Y_{0} + E_{\tau_{\mathcal{C}2}})},
\label{tauc2}
\end{align}
where
\begin{align}
A_{\tau_{\mathcal{C}2}}&=\theta_+\,z_{0}^{4}+2\,(b_{1}+3)\,z_{0}^{2}+1,\label{eq:A.16} \\
B_{\tau_{\mathcal{C}2}}&=2\,c\,(\theta_-\,z_0^2+1),\label{eq:A.17}  \\
C_{\tau_{\mathcal{C}2}}&=\theta_+^{2}\,z_{0}^{4}+2\,(b_{1}+3)\,z_{0}^{2}+1,\\
D_{\tau_{\mathcal{C}2}}&=4\,c\,(\theta_-\,z_0^2+1), \\
E_{\tau_{\mathcal{C}2}}&=4\,c^{2}\,\vartheta_0.
\end{align}
Therefore, the intersection points of the Hugoniot curve with $\mathcal{C}$
are $(z_{\mathcal{C}1}, \tau_{\mathcal{C}1},0)$
and $(z_{\mathcal{C}2},\tau_{\mathcal{C}2},0)$.

Substituting these values into Eq.~\eqref{eq:speed} we obtain, respectively,
\begin{align*}
\sigma_{\mathcal{C}1}=
&\dfrac{\{(b_{1}+1)^{3}z_{0}^{4}+2[(b_{1}+1)^{2}-2]z_{0}^{2}+(b_{1}+1)\}Y_{0}}
{2z_{0}b_{1}[(b_{1}+1)z_{0}^{2}+3]}\\
&+\dfrac{2c[(b_{1}+1)^{2}z_{0}^{2}+2z_{0}^{2}+(b_{1}+1)]}
{2z_{0}b_{1}[(b_{1}+1)z_{0}^{2}+3]}
\end{align*}
and
\begin{equation}
\sigma_{\mathcal{C}2}
=\frac{[(b_{1}+1)z_{0}^{2}-1]^{2}Y_{0}+6c(b_{1}+1)z_{0}^{2}+2c}
{2b_{1}z_{0}[(b_{1}+1)z_{0}^{2}+3]}.
\end{equation}

Simple inspection shows that
$\sigma_{\mathcal{C}2}=\sigma_{son'}(z_{0},\tau_{0}',Y_{0})$.
Therefore, the coordinates $z_{\mathcal{C}2}$ and $\tau_{\mathcal{C}2}$
are the ones we seek.

According to Prop.~\ref{th:slow_fast},
we study the sign of $\tau_{\mathcal{C}2}$.
If $\tau_{\mathcal{C}2}>0$,
the point $(z_{0},\tau_{0}',Y_{0})$ is in $Son'_{f}$;
otherwise, it is in $Son'_{s}$.

These computations are not valid for $z_{0}=0$;
another coordinate system must be used
to ensure the validity of our results,
for instance, $Z=1/z$, $T=z\tau$, and $X=z\,Y$.
\end{proof}

\subsection{}
\begin{proof}[Proof of Prop.~\ref{proposition:son'sson'f}]
To understand the properties of the system, we need to examine the sign of
$\tau_{\mathcal{C}2}$, as defined in Eq.~\eqref{tauc2}. The denominator of this
expression can be factored into the product of $2cz_{0}[(b_{1}+1)z_{0}^{2}+3]$
and a quadratic polynomial in $Y_{0}$ with negative discriminant. Since the
coefficient of $Y_{0}^{2}$ is positive, the denominator has the same sign as
$z_{0}$, and this sign does not change since these computations are not valid
for $z_{0}=0$. The numerator is of the form
$-(A_{\tau_{\mathcal{C}2}}(z_{0})Y_{0}+B_{\tau_{\mathcal{C}2}}(z_{0}))$
multiplied by a positive term, where $A_{\tau_{\mathcal{C}2}}$ and
$B_{\tau_{\mathcal{C}2}}$ are both positive and are given by the expressions after 
\eqref{tauc2}.

If $z_{0}>0$, then the point lies on $Son'_{s}$ if
$A_{\tau_{\mathcal{C}2}}Y+B_{\tau_{\mathcal{C}2}}>0$ and in $Son'_{f}$ if
$A_{\tau_{\mathcal{C}2}}Y+B_{\tau_{\mathcal{C}2}}<0$. If $z_{0}<0$, then the
point lies on $Son'_{s}$ if
$A_{\tau_{\mathcal{C}2}}Y+B_{\tau_{\mathcal{C}2}}<0$ and in $Son'_{f}$ if
$A_{\tau_{\mathcal{C}2}}Y+B_{\tau_{\mathcal{C}2}}>0$.
The curve $A_{\tau_{\mathcal{C}2}}Y+B_{\tau_{\mathcal{C}2}}=0$
is the projection of the $ECC\,'$ curve onto the $(z,Y)$-plane
(the third equation in \eqref{eq:parsonli}).

These computations are not valid for $z = 0$
because $z$ appears in the denominator of \eqref{tauc2}.
To check what happens near $z=0$,
we use the projection of the curve $ECC\,'$
onto the $(z,\tau)$-plane.
We follow the five steps at the beginning of Sec.~\ref{sec:sonlifs}.

Given a point $(z_1, \tau_1,Y'_1)$ on $Son'$,
we obtain from Eq.~\eqref{eq:son'} that
\begin{equation}
Y'_1=\frac{2c[(b_1+1)\tau_1z_1^5+(b_1+4)\tau_1z_1^3+(b_1-1)z_1^2
+ 3\tau_1z_1+1]}{(z_1^2+1)[(b_1+1)z_1^2-1]}.
\end{equation}
We observe that $z$ does not appear as a factor in the denominator of $Y'_1$.
These calculations are not valid for $z_1^2=1/(b_1+1)$.
We are interested in values of $z_1$
such that $z_{crit2}< z_1 <z_{crit1}$.
The shock speed at this point is obtained by changing $(Y,\tau,Z)$
into $(z_1,\tau_1,Y'_1)$ in Eq.~\eqref{eq:speed}:
\begin{equation}
\sigma_{\tau z}(z_1,\tau_1,Y'_1)
=\frac{c[\tau_1(b_1+1)z_1^4+b_1\tau_1z_1^2+(b_1+2)z_1-\tau_1]}{b_1(z_1^2+1)}.
\label{eq:speedtauz}
\end{equation}

Parametric equations for the Hugoniot curve of $(z_1,\tau_1,Y'_1)$
are obtained changing $z_0$, $\tau_0$, $Y_0$ into $z_1$, $\tau_1$, $Y'_1$
in Eqs.~\eqref{eq:hugponto} giving $(z, \tau_{Hug}(z),Y_{Hug}(z))$.
Solutions of $Y_{Hug}=0$ give the $z$ coordinates of the intersection points
of Hugoniot curve with $\mathcal{C}$,
$z_{\mathcal{C}3}=[(b_{1}+1)z_{1}^{2}+1]/2z_{1}$
and $z_{\mathcal{C}4}=-[2(z_1^2+1)\tau_1z_1-z_1^2+1]
/\{(b_{1}+1)\tau_1z_1^4+(b_1+2)\tau_1z_1^2+b_1z_1+\tau_1\}$.

Substituting $z_{\mathcal{C}3}$ and $z_{\mathcal{C}4}$
into the expression of $\tau_{Hug}$,
we obtain the $\tau$ coordinates of the intersection points:
\begin{equation}
\tau_{\mathcal{C}3}=\frac{4z_1^2[\eta_1(z_1)\tau_1+z_1\eta_2(z_1)]}
{[(b_1+1)z_1^2-1]\eta_1(z_1)}
\end{equation}
\begin{equation}
\tau_{\mathcal{C}4}=\frac{-(\eta_3(z_1)
\tau_1+b_1z_1)^2[\eta_1(z_1)\tau_1
+\eta_2(z_1)]}{(z_1^2+1)((b_1+1)z_1^2-1)
((z_1^2+1)\eta_1(z_1)\tau_1^2+2(z_1^2+1)
\eta_2(z_1)+\eta_4(z_1))},
\end{equation}
where $\eta_1(z_1)=(z_1^2+1)[(b_1+1)^2z_1^4+2(b_1+3)z_1^2+1]$,
$\eta_2(z_1)=(b_1+2)z_1[(b_1-1)z_1^2+1]$,
$\eta_3(z_1)=(b_1+1)z_1^4+(b_1+2)z_1^2+1$,
and $\eta_4(z_1)=z_1^4+(b_1^2-2)z_1^2+1$.

So the intersection points of the Hugoniot curve of $(z_1,\tau_1,Y'_1)$
are $\mathcal{C}_3=(\tau_{\mathcal{C}3},z_{\mathcal{C}3},0)$
and $\mathcal{C}_4=(\tau_{\mathcal{C}4},z_{\mathcal{C}4},0)$.
Substituting $\mathcal{C}_3$ and $\mathcal{C}_4$
into Eq.~\eqref{eq:speed} we get $\sigma_{\mathcal{C}3}$ and
\begin{equation}
\sigma_{\mathcal{C}4}=\frac{c[(b_1+1)\tau_1 z_1^{4}+b_1\tau_1 z_1^{2}
+(b_1+2)z_1-\tau_1]}{b_1(z_1^2+1)}.
\end{equation}

Straightforward computations show that
$\sigma_{\mathcal{C}3}-\sigma_{\tau z}\neq 0$.
A simple inspection shows that $ \sigma_{\mathcal{C}4}=\sigma_{\tau z}$.
So, we need to study the sign of $ \sigma_{\mathcal{C}4}$.

The denominator of $\sigma_{\mathcal{C}4}$ is the product of a positive term,
a quadratic polynomial in $\tau_1$ with negative discriminant,
so it is positive.
It also involves a quadratic polynomial in $z_1$,
which has roots $z_{crit1}$ and $z_{crit2}$ as in \eqref{doublecon}.
Since we are interested in $z_{crit1} < z_1 < z_{crit2}$,
the denominator of $\sigma_{\mathcal{C}4}$ remains negative.
Therefore, the sign of $\sigma_{\mathcal{C}4}$ is the same as the
sign of its numerator. Therefore, the sign of $\sigma_{\mathcal{C}4}$
is the same as the sign of its numerator,
which is of the form $\eta_1(z_1)\tau_1+\eta_2(z_1)$
multiplied by a positive term.

Hence, $\sigma_{\mathcal{C}4}$ has the same sign as
$\eta_1(z_1) \tau_1+\eta_2(z_1)$.
Therefore, if $\eta_1(z_1) \tau_1+\eta_2(z_1) > 0$,
then the point lies in $Son'_f$;
otherwise, it lies in $Son'_s$.
Solving for $\tau_1$, we obtain the second equation in \eqref{eq:parsonli},
which is the equation of the projection of $ECC\,'$ onto the $(z,\tau)$-plane.
It follows that the curve $ECC\,'$ subdivides $Son'_{s}$ and $Son'_{f}$
in a neighborhood of $z=0$.
\end{proof}

\section{Normal Forms and Coordinates}
\label{sec:Normal_Forms_and_Coordinates}

\subsection{Normal form for the flux}
\label{subsec:Normal_form_for_the_flux}
This paper, like many previous ones,
studies systems of two conservation laws \eqref{eq:cons-law} for $F=(f,g)$
with $f(u,v)=v^{2}/2+(b_{1}+1)\,u^{2}/2+a_{1}\,u+a_{2}\,v$
and $g(u,v)=u\,v-b_2\,v^2/2+a_{3}\,u+a_{4}\,v$.
Here, we focus on the so-called symmetric case, $b_{2}=0$.

In \cite{schaeffer87},
Eq.~\eqref{eq:cons-law} was studied under the hypothesis that
$DF$ is strictly hyperbolic (\emph{i.e.}, has distinct eigenvalues)
everywhere except at $(0,0)$,
where $DF(0,0)=I_{2\times 2}$,
\emph{i.e.}, $(0,0)$ is an umbilic point.
This led to the study of quadratic $F$,
and it was shown that it is sufficient to take $F=\operatorname{grad}(C)$
with $C=a\,u^{3}/3+b\,u\,v+u\,v^{2}$.

In \cite{palmeira88},
the study of $F=\operatorname{grad}(C) + \text{\emph{linear terms}}$
such that $F$ is no longer a gradient was begun.
The addition of linear terms replaced the umbilic point with an elliptic region.
Before adding linear terms, a simple change of coordinates was made so that
$C$ became $uv^{2}/2+[(b_{1}+1)u^{3}]/6-b_{2}v^{3}/6$.
This change of coordinates was chosen to simplify the differential equation
of the eigenspaces of $DF$. Subsequent papers kept this choice of $F$.

\subsection{Coordinates}
\label{subsec:Coordinates}
The coordinates $\widetilde U$, $\widetilde V$, $X$, $Y$, and $Z$
were introduced in \cite{marpal94b}.
In these coordinates, $\mathcal{W}$ was given by the equations
$(1-Z^{2})\widetilde V-Z\widetilde U+c=0$ and $Y=ZX$.
Both of the natural coordinate choices,
$(\widetilde U, Z, X)$ and $(\widetilde V, Z, X)$, were used.

In the present paper, we are treating symmetric Case~IV
in the classification of \cite{schaeffer87}.
In this case, the secondary bifurcation line is contained in plane $Z=0$;
to stay away from it, we use the coordinate $z=1/Z$.
For convenience, we replace $\widetilde V$ by $V_1 = \widetilde V + a_2$.
Thus, the equations of $\mathcal{W}$ become
$(z^{2}-1)\,V_{1}-z\,\widetilde U+c=0$ and $X=z\,Y$.
The solution set is the characteristic manifold $\mathcal{C}$.

Each slice of $\mathcal{C}$ for fixed $z$
is a horizontal line in $(\widetilde U, V_{1}, z)$-space.
In the $(\widetilde U, V_{1})$-plane,
this line is spanned by the vector $(z^{2}-1,z)$.
The coincidence curve is the singular set of the projection
$(\widetilde U, V_{1},z) \mapsto (\widetilde U, V_{1},0)$,
as restricted to $\mathcal{C}$.
Putting $h(\widetilde U, V_{1}, z) = (z^{2}-1)\,V_{1}-z\,\widetilde U+c$,
we obtain a parametrization of
the coincidence curve by solving the linear system
\begin{equation}
\text{$h(\widetilde U, V_{1},z)=0$ and $\frac{dh}{dz}(\widetilde U, V_{1},z)=0$}
\end{equation}
for $\widetilde U$, $V_{1}$;
the result is $\widetilde U=2\,c\,z/(z^{2}+1)$ and $V_{1}=c/(z^{2}+1)$.
To define the coordinate $\tau$,
we take the rules of $\mathcal{C}$ to originate at the coincidence curve.
For fixed $z$, the coordinate $\tau$ measures distance
from the coincidence curve along the rule:
$\widetilde U = 2cz/(z^{2}+1) + c\,\tau\,(z^{2}-1)$
and $V_{1}=c/(z^{2}+1) +c\,\tau\,z$.

\subsection*{Acknowledgments}
Computations were assisted by
the symbolic manipulation software package
Maple\textsuperscript{\textregistered}
and results were checked by the ELI interactive solver \cite{ELISoft}.
The work of M.~M. L\'{o}pez-Flores was supported by IMPA PCI program
under grant 300758/2022-7.
The work of D. Marchesin was supported by CAPES grant 88881.156518/2017-01,
by CNPq under grants 405366/2021-3, 306566/2019-2,
and by FAPERJ under grants E-26llf/210.738/2014,
E-26/202.764/2017, and E-26/201.159/2021.

\addcontentsline{toc}{section}{References}

\bibliographystyle{amsplain}
\providecommand{\bysame}{\leavevmode\hbox to3em{\hrulefill}\thinspace}
\providecommand{\MR}{\relax\ifhmode\unskip\space\fi MR }
\providecommand{\MRhref}[2]{%
 \href{http://www.ams.org/mathscinet-getitem?mr=#1}{#2}
}
\providecommand{\href}[2]{#2}

\end{document}